\documentclass[11pt,reqno]{amsart}

\usepackage{times}
\usepackage[T1]{fontenc}
\usepackage{mathrsfs}
\usepackage{latexsym}
\usepackage[dvips]{graphics}
\usepackage{graphicx}
\usepackage{fancyhdr,amssymb}
\usepackage{epsfig}
\usepackage{amsmath,amsfonts,amsthm,amssymb,amscd}
\input amssym.def
\input amssym.tex
\usepackage{color}
\usepackage{hyperref}
\usepackage{url}
\usepackage{breakurl}
\newcommand{\bburl}[1]{\textcolor{blue}{\url{#1}}}

\usepackage[margin=1in]{geometry}
\usepackage{mathtools}

\newcommand\be{\begin{equation}}
\newcommand\ee{\end{equation}}
\newcommand\bea{\begin{eqnarray}}
\newcommand\eea{\end{eqnarray}}
\newcommand\bi{\begin{itemize}}
\newcommand\ei{\end{itemize}}
\newcommand\ben{\begin{enumerate}}
\newcommand\een{\end{enumerate}}
\newcommand\bc{\begin{center}}
\newcommand\ec{\end{center}}
\newcommand\ba{\begin{array}}
\newcommand\ea{\end{array}}

\newcommand{\R}{\ensuremath{\mathbb{R}}}

\newcommand{\Z}{\ensuremath{\mathbb{Z}}}
\newcommand{\Q}{\mathbb{Q}}
\newcommand{\N}{\mathbb{N}}
\newcommand{\F}{\mathbb{F}}

\newtheorem{thm}{Theorem}[section]
\newtheorem{conj}[thm]{Conjecture}

\newcommand{\ncr}[2]{{#1 \choose #2}}
\newcommand{\twocase}[5]{#1 \begin{cases} #2 & \text{#3}\\ #4
&\text{#5} \end{cases}   }

\newcommand{\diag}{diag}

\newcommand{\fixed}[2][1]{%
  \begingroup
  \spaceskip=#1\fontdimen2\font minus \fontdimen4\font
  \xspaceskip=0pt\relax
  #2%
  \endgroup
}

\begin{document}


\title{Combinatorial and Additive Number Theory Problem Sessions: '09--'24}

\author[Miller]{Steven J. Miller}
\email{\textcolor{blue}{\href{mailto:sjm1@williams.edu}{sjm1@williams.edu}},  \textcolor{blue}{\href{Steven.Miller.MC.96@aya.yale.edu}{Steven.Miller.MC.96@aya.yale.edu}}}
\address{Department of Mathematics and Statistics, Williams College, Williamstown, MA 01267}



\maketitle




\begin{abstract} These notes are a summary of the problem session discussions at various CANT (Combinatorial and Additive Number Theory Conferences). Currently they include all years from 2009 through 2018  (inclusive); the goal is to supplement this file each year. These additions will include the problem session notes from that year, and occasionally discussions on progress on previous problems. If you are interested in pursuing any of these problems and want additional information as to progress, please email the author.

For more information, visit the conference homepage at \begin{center}\bburl{http://www.theoryofnumbers.com/}\end{center} or email either the typist at \textcolor{blue}{\href{mailto:sjm1@williams.edu}{sjm1@williams.edu}} or \textcolor{blue}{\href{Steven.Miller.MC.96@aya.yale.edu}{Steven.Miller.MC.96@aya.yale.edu}}, or the organizer at \textcolor{blue}{\href{mailto:melvyn.nathanson@lehman.cuny.edu}{melvyn.nathanson@lehman.cuny.edu}}. \\

\textbf{\textcolor{red}{Warning: Many of these notes were LaTeX-ed in real-time by Steven J. Miller; all errors should be attributed solely to him.}}

\end{abstract}

\thanks{The typist was supported by NSF grants DMS0600848, DMS0600848 and DMS1265673. These notes would not have been possible without the help of the participants, especially Huixi Li, Zack McGuirk, Kevin O'Bryant, Steven Senger, my students, especially Olivia Beckwith, Alan Chang, Ginny Hogan, Jared D. Lichtman, Jasmine Powell, Ryan Ronan, Maddie Weinstein and of course the organizer Mel Nathanson.}

\tableofcontents


\newpage
\setcounter{section}{2008}
\section{CANT Problem Sessions: 2009}

\subsection{Problem Session I: Tuesday, May 26th (Chair Kevin O'Bryant)}

\subsubsection{Steven J Miller: I (sjm1@williams.edu)}

\textbf{Probability an element is in an MSTD}

Let $\gamma(k,n)$ be the probability that $k$ is in an MSTD set $A$ with $A \subset [0,n]$; see for instance the figure below

\begin{figure}

\begin{center}

\scalebox{1.00}{\includegraphics{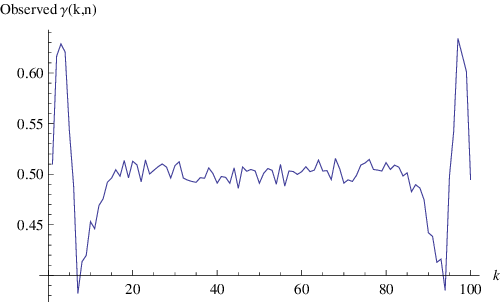}}

\caption{\small\textbf{Observed $\gamma(k,100)$, random sample 4458 MSTD sets.}}

\end{center}\end{figure}

\ \\

\textbf{Conjecture:} Fix a constant $0<\alpha<1$. Then $\lim_{n\rightarrow\infty}\gamma(k,n)=1/2$ for $\lfloor \alpha n \rfloor \le k \le n - \lfloor \alpha n \rfloor$.\\

Questions:

\bi

\item How big are the spikes? Do the sizes of the spikes tend to zero as $n \to \infty$?

\item Is the spike up equal to the spike down?

\item Study more generally $g(n) \le k \le n - g(n)$; is it sufficient for $g(n) \to \infty$ monotonically at any rate to have all $k$ in the region above having probability 1/2 of being in an MSTD set? Can we take $g(n) = \log\log\log n$, or do we need $g(n) = \alpha n$?

\item A generic MSTD set has about n/2 elements; what other properties of a generic set are inherited by an MSTD set?

\item How big are the fluctuations in the middle?

\item Do we want to look at all MSTD sets containing $1$ and $n$, or do we want to just look at all subsets of $[1,n]$ that are MSTD sets

\ei

\textbf{Note since the end of the conference: Kevin O'Bryant has observations relevant to this and other problems proposed by Miller.}

\textbf{Note added in 2014: Many of these claims were later proved by Zhao.}

\subsubsection{Steven J Miller: II (sjm1@williams.edu)}

With Dan S. and Brooke O. we constructed very dense families of MSTD sets in $[0,r]$ such that these families were $C/r^4$ of all subsets of $[0,r]$. This isn't a positive percentage of sets, but it is significantly larger than any previous family. Can one do better? Can one get a positive percentage?

\subsubsection{Peter Hegarty: I (hegarty@chalmers.se)}

Smallest size of an MSTD is 8 elements: $A =$ $\{0, 2$, $3$, $4$, $7$, $10, 12, 14\}$; remove $4$ and symmetric about $7$. If look in $\Z^2$, can construct an MSTD set of size 16 from this: take $A \times \{0,1\}$ (any set of size two would work). Can you construct MSTD sets in $\Z^2$ without going through an MSTD set in $\Z$. Need a computer to show this set $A$ was minimal (about 15 hours to find all MSTD sets of size 8, and thus see that this set $A$ is minimal). To find all MSTD sets up to isomorphism of a certain size is a finite computation, but practically impossible for 9.

\subsubsection{Peter Hegarty: II (hegarty@chalmers.se)}

\textbf{Question:} What are the possible orders of a basis for $\Z_n$? \\

Let $A \subset \Z_n$. We say $A$ is a basis of order $h$ if $hA = \Z_n$ but $(h-1)A \neq \Z_n$. $A$ is a basis of some order if and only if $(\gcd(A),n) = 1$.

$A = \{a_1,\dots,a_k\}$, $|hA| = O(k^h)$, order should be about $\log n / \log k$, so $k^h = n$. Order for a random set to be a basis, should be of logarithmic order. Can write down a very inefficient basis where need long summands to cover all of $\Z_n$. To do this, take $A = \{0,1,\dots,k-1\}$. Order of this will be essentially $\frac{n}{k-1}$.

\ \\

\textbf{Conjecture:} If the order of a basis is $\Theta(n)$ then the order must be very close to $n/k$ for some $k$. So no number between $n/2$ and $n$ can be the order. Gives gaps. See paper by Dukes and Herke. \\

Question from participants: this can't go on forever?

Answer from Peter: Can go on forever. Fix a $k$, let $n\to\infty$, the way you would phrase to make it precise: Fix a $k$. For $n \gg 0$ if the order of a basis is greater than $\frac{n}{k+1} + O(k)$ then the order must be within $O(k)$ of $n/\ell$ for some $\ell \le k$.

\textbf{Note : Since the end of the conference, Peter Hegarty has solved this problem. His result is available at \bburl{http://arxiv.org/abs/0906.5484}}

\subsubsection{Kevin O'Bryant (obryant@gmail.com)}

Take $g_0 = 0$, $g_i$ to be the least positive integer such that $\{g_0,g_1,\dots,g_i\}$ has no solutions to $5w+2x=5y+2z$. This is building a set greedily.

Let  $a_0 = 0$ and \be a_i \ = \ \lfloor \frac{5 + 7\sum_{j=0}^{i-1} a_j}{2} \rfloor. \ee Let $A$ equal the sum of distinct $a_i$'s.\\

\textbf{Conjecture:} $G = \{5x+y: x \in A, 0 \le y \le 4\}$.

\ \\

Appears computationally that there is some description of this sort when one number is at least twice as large as another; can replace $(5,2)$ with $(11,4)$ without trouble, but not with $(4,3)$. \\

Question from audience: Why 5 and 2?

Answer: 5, 2 smallest haven't solved and have done the most computation.

Question from audience: How many other cases investigated?

Answer: Calculated all terms up to about 100,000 if both numbers at most 12 (and can exclude cases, such as cases with common prime factors). Nice structure if one is twice the other, else irregular and nothing to say (though all irregular in the same way).

\subsubsection{Ruzsa (through Simon Griffiths through Kevin O'Bryant)}

Let $A \subset \Z$, $|A| = n$, and \be S_k \ = \ \{\sum_{a \in B} a: B \subset A, \ |B| = k\}. \nonumber\ee Note $|S_k| = |S_{n-k}|$. For example, if $A = \{1,2,4,8,\dots,2^{n-1}\}$ then $|S_{k+1}| = \ncr{n}{k+1} = \frac{n-k}{k+1} |S_k|$.\\

\textbf{Question:} $|S_{k+1}| \le \frac{n-k}{k+1} |S_k|$ whenever $k < n/2$?\\

\textbf{Theorem (Ruzsa):} Yes, when $n > \frac{k^2 + 7k}{2}$. \\

\textbf{Exercise:} $|S_{k+1}| \le \frac{n}{k+1} |S_k|$.


\subsection{Problem Session II: Wednesday, May 27th}

\subsubsection{Comments after Nathanson's Talk (Mel Nathanson:\\ melvyn.nathanson@lehman.cuny.edu)}

Paper is online at \bburl{http://arxiv.org/pdf/0811.3990}.

Take $G_i$ with $A_i$ of generators. Only one direct product, but many sets of generators that can construct from generators of individual groups. Could take direct product of generators. That's complicated. Given groups and generating sets, many ways to put together new sets of generators. Never thought about finite groups because thinking about geometric group theory. For finite groups know at some point all spheres empty.

Question: A result like this might not be true for semi-groups: bunch of things with finite spheres then empty at some point. Additive sub-model of integers, all positive integers exceeding 1000. Can you have infinite sphere, finite sphere, infinite sphere, finite sphere.... Answer: don't know. Wanted to create oscillating sets of spheres, turned out couldn't.

\subsubsection{Constructing MSTD Sets (Kevin O'Bryant, communicated to Steven Miller)}

\ \\

\textbf{Theorem:} $d_i \in \{3,4,5\}$ independent uniformly distributed, $x_1 = 4$, $x_2 = 5$, $x_i = x_{i-1} + d_i$, $A = \{1\} \cup \{0, \pm x_1, \dots, \pm x_n\}$. Then $|A + A| > |A-A|$ with probability 1.\\

Note \bea A+A &\  = \ &  (X+X) \cup (X+1) \cup \{2\} \nonumber\\ A-A &=& (X-X) \cup (X-1) \cup (1-X) \cup \{0\} \nonumber\\ &=& (X+X) \cup (X-1) \cup (1+X) \cup \{0\}, \nonumber \eea where $X$ is the set of the $x_i$'s.

\subsubsection{David Newman (davidsnewman@gmail.com)}

Suppose we have a basis for the non-negative integers, that is a set so that for any non-negative number we can find two elements of the set whose sum is this given number. If we arrange the numbers in the set in ascending order then we can cut it off at a certain point and look at the first $N$ terms of this basis. \\

\textbf{Question:} Can this beginning of a basis be extended into a minimal basis? By minimal basis I mean a basis where if you remove any element it is no longer a basis.\\

Has to have 0 and 1 as a start. I think the answer is yes is because I haven't seen the beginning of a basis I couldn't extend to a minimal basis. I have an algorithm implemented in Mathematica and in a few seconds gives a set which is a minimal basis. That's about all the info I have, other than one family of bases that I can always extend to a minimal basis.\\

Another problem (from the theory of partitions): Consider \be (1+x)(1+x^2) \cdots \ = \ \sum a_n x^n,\nonumber \ee which is the generating function for partitioning into distinct parts. Now put in minus signs: \be (1-x)(1-x^2) \cdots \ = \ 1 - x - x^2 + x^5 + x^7 + \cdots, \nonumber \ee where all coefficients are in $\{0,\pm 1\}$. Now do partitions into unrestricted parts: \be (1+x+x^2+x^3 +\cdots) (1+x^2 + x^4 + x^6 + \cdots) \cdots \ = \ \sum b_n x^n. \nonumber \ee \\

\textbf{Question:} can we change some of the signs above into minus signs so that the $b_n$'s are also in $\{0, \pm 1\}$.

\ \\

\textbf{Note since the conclusion of the conference: Peter Hegarty and David Newman have made progress on this. They are currently working on a paper: Let $h > 1$ be an integer, for any basis $A$ for $\N_0$ of order $h$ and any $n \in \N_0$ the initial segment $A \cap [0,n]$ can be extend to a basis of $A'$ of order $h$ which is also a minimal asymptotic basis of this order.}

\subsubsection{Infinitude of Primes (Steven Miller)}

Two types of proofs of the infinitude of primes, those that give lower bounds and those that don't (such as Furstenberg's topological proof). What category does $\zeta(2) = \pi^2/6 \neq \Q$ fall under? It implies there must be infinitely many primes, as this is $\zeta(2) = \prod_p (1 - p^{-2})^{-1}$; if we knew how well $\pi^2$ can be approximated by rationals, we could convert this to knowledge about spacings between primes. Unfortunately while we know the irrationality exponent for $\pi^2$ is at most 	5.441243 (Rhin and Viola, 1996), their proof uses the prime number theorem to estimate ${\rm lcm}(1,\dots,n)$. This leads to $\pi(x) \gg \log\log x /\log\log\log x$ infinitely often; actually, I can show: let $g(x) = o(x/\log x)$ then $\pi(x) \ge g(x)$. A preprint of my paper (with M. Schiffman and B. Wieland) is online at \begin{center} \bburl{http://arxiv.org/PS_cache/arxiv/pdf/0709/0709.2184v3.pdf} \end{center} and I hope to have a final, cleaned up version in  a few months. I'm looking for a proof of the finiteness of the irrationality measure of $\zeta(2)$ that doesn't assume the prime number theorem. \textbf{Note added in 2014: Miller is currently working on this with some of his students.}

\subsubsection{Kent Boklan (boklan@boole.cs.qc.edu)}

There are infinitely many primes, don't know much about twin primes. Know sum of reciprocals of twins converges by Brun's theorem. This is a hard theorem -- I want to do elementary things. How do you show there are infinitely many primes which are not twin primes. Trivial proof: There are infinitely many primes of the form $15k+7$ by Dirichlet, and not prime if add or subtract 2. But Dirichlet isn't elementary!

\subsubsection{Mel Nathanson II}

$A = \{a_1,\dots,a_k\}$ finite set of integers, $n = \sum_{i=1}^k a_i x_i$ is solvable for all $n$ if and only if ${\rm gcd}(A) = 1$. In geometric group theory, can deduce algebraic properties of the group by seeing how it acts on geometric objects. Fundamental lemma of geometric group theory says the following: $G$ is a group and acts on a set $S$ (metric space), want action to be an isometry for any fixed $g$ in the group. Acts isometrically on the metric space $S$. Suppose the space is nice (Heine-Borel, want that: any closed and bounded set is compact, call this a proper space). $G$ acts properly discontinuously on $X$ if intersection non-empty for only finitely many $g$. Example: $\Z^n$ acts on $\R^n$ by $(g,x) = g + x$. Have $G\backslash X$, send $x$ to its orbit $\langle x \rangle$. Put a quotient topology on $G\backslash X$ that makes projection map continuous. Example: $\Z^n \backslash \R^n$ is the n-torus. Let $K \subset X$ compact, for every $x \in X$ there is a $y
  \in K$ such that $gy = x$. For example, $n=1$: $\Z$ acts on $\R$ by translation, take unit interval $[0,1]$ (compact), and every number is congruent modulo 1 to something in unit interval. The fundamental lemma of geometric group theory: Group acts as isometry and properly discontinuously on proper metric space then $G$ must be finitely generated. Know nothing about if it is finitely or infinitely generated, but if acts geometrically in this nice way, that can only happen if the group is finitely generated. Proof goes by finding a compact set $K$ with exactly the property above. What we know about $K$ since group action properly discontinuous, group action under $K$ only finitely many, that is a finite set of generators.

Suppose we specialize to elementary number theory: integers acting on reals by translations, compact set $K$ such that every real number is congruent modulo 1 to an element of $K$. Then we get a finite set of generators for the group, but the group is the integers and a finite set of integers is a finite set of relatively prime integers. Get certain sets of relatively prime integers. What finite sets of integers can we get geometrically in this way? Every finite set of relatively prime set of integers can be obtained this way. Curious thing is that there is this geometric way to describe these sets. Look at lattice points in two dimensions, seems quite complicated.

Article might be: \bburl{http://arxiv.org/pdf/0901.1458}.




\subsection{Problem Session III: Thursday, May 28th}

\subsubsection{Gang Yu (yu@math.kent.edu)}

$C \subset \N$ is an infinite sequence, $C(N) = C \cap [N]$ (where $[N] = \{0,\dots,N\}$ or perhaps it starts at 1), $h \ge 2$ fixed, call $A \subset [N]$ an $h$-basis of $C(N)$ if $hA \supset C(N)$.

Trivial estimate: $|A| \ge h! |C(N)|^{1/h}$.

Interesting cases: $C$ is sparse but arithmetically nice: \bea C & \ = \ & \{n^2\} \nonumber\\ C &=& \{n^\alpha\} \nonumber\\ C &=& \{f(n)\} \nonumber\ \eea where $f$ is a degree 2 polynomial. Let \be \Gamma_h/C \ = \ \overline{\lim_{N\to\infty}} \frac{|C(N)|}{D_h(C,N)^h}, \ee where \be D_h(C,N) \ = \ \min_{A \subset [N] \atop A\ {\rm is\ an\ }h{\rm-basis\ of\ } C(N)} |A|. \ee

\ \\

\textbf{Question:} Is $\Gamma_h(C) = 0$ for polynomial $C$ (ie, degree at least 2)?\\

Audience: Is it true for any sequence?

Gang: Don't know. For $AP = C$, bounded away from 0. Specifically, $A \subset [N]$, $A+A \supset \{n^2: D \le n \le \sqrt{N}\}$, $N^{1/4} = o(|A|)$?

\textbf{After the conference it was noted: Some information available at
\bburl{http://arxiv.org/pdf/0711.1604}.}

\subsubsection{Simon Griffiths (sg332@cam.ac.uk)}

An $n$-sum of a sequence $x_{1},\dots ,x_{r}$ is a sum of the form $x_{i_{1}}+\dots +x_{i_{n}}$ where $i_{1}<\dots <i_{n}$, i.e. an element that can be obtained as the sum of an $n$-term subsequence.

EGZ: Every sequence $x_1, \dots, x_{2n-1} \in \Z_n$ has 0 as an $n$-sum.

Bollob\'{a}s-Leader: Let $x_1,\dots,x_{n+r} \in G$ and suppose $D$ is not an

$n$-sum, then you have at least $r+1$ $n$-sums.

Examples: EGZ is tight as demonstrated by the sequence of $n-1$ 0s and $n-1$ 1s; Bollob\'{a}s-Leader is tight as deomnstrated by the seq of $n-1$ 0s and $r+1$ 1s.

What about finite abelian groups more generally?

$D(G)$ is Davenport constant, the minimum $r$ where every $r$-term sequence has a non-trivial subsequence with sum $0$.  For example: not difficult to show $D(\Z_n)=n$.

Example: Let $x_{1},\dots,x_{D(G)-1}$ be a sequence in $G$ with no non-trivial subsequence summing to $0$, and adjoin $n-1$ 0s by setting $x_{D(G)},...,x_{n+D(G)-2}=0$.  Then, by an easy check, we see that this sequence, of length $n+D(G)-2$ does not have $0$ as an $n$-sum.

Gao: Every sequence $x_1,\dots, x_{n+D(G)-1}$ has 0 as an $n$-sum. \\

\textbf{Question:} EGZ is to Bollob\'{a}s-Leader as Gao is to ..?..?..\\

One Answer: A theorem of Gao and Leader.

Why do we need another answer: Both of the results, Bollob\'{a}s-Leader and Gao-Leader allow us to see the set of $n$-sums grow as the length of the underlying sequence increases.  However perhaps in the case of general abelian groups there may be a different way to see this growth - to see this growth as a growth of dimension in the sense described below.

Our approach to defining the dimension of a subset $S\subset G$ is similar to describing the dimension of a subspace via the maximum dimension of an independent subspace.  Call a sequence zero-sum-free if no non-trivial subsequence has sum $0$.  Let ${\rm dim}(S)$ equal $D(G)$ minus the minimum $r$ such that for every zero-sum-free sequence

$y_1,\dots,y _r$ there exists an $s\in S$ and a subsequence $I$ such that

$s+\sum_{I} y_i = 0$. \\

If $0\in S$ then we take the minimum to be zero.  Thus,

Examples: ${\rm dim}(G) = D(G)$. $S = G - \{0\}$ implies ${\rm dim}(S) = D(G)-1$. $S = \emptyset$ implies ${\rm dim}(S) = 0$.

\ \\

\textbf{Conjecture:} $x_1,\dots, x_{n+r}$ either 0 as an $n$-sum or ${\rm
dim}\left(\{n - {\rm SUMS}\}\right) \ge r+1$. \\



\subsection{Problem Session IV: Saturday, May 30}

\subsubsection{Urban Larsson}

How small can a {\em maximal} AP-free set be? Specifically, how large is the smallest maximal (with respect to not having 3 terms in arithmetic progression) subset of $[n]$? Set
    \[\mu(n) := \min_{\substack{A \\ \text{$A$ is 3-free}}} \left| A \cap [n] \right|.\]

Examples: The greedy subset of $\{0,1,\dots\}$ with 3-term APs is $\{ n \in {\mathbb N}$ ${\rm col}on$ $\text{base-3 expansion of $n$ has no '2's}\}$. This shows that $\mu(3^t)\leq 2^t$. A better example is the set of natural numbers whose base-4 expansions have neither '2's nor '3's. This gives $\mu(4^t)\leq 2^t$.

Each pair of elements of $A$, which is 3-free, forbids at most three other numbers from $A$, so $3\binom{|A|}{2} \geq n - |A|$, so that $|A|\geq c \sqrt{n}$. Attention to detail gives $|A|\geq \sqrt{2n/3}$. If $A$ is uniformly distributed mod 4, and u.d. in $[n]$, then many pairs will not forbid three other numbers, and this gives $|A|\geq \sqrt{420n/401}$.

I conjecture that $\mu(4^t)=2^t$, and $\mu(n)\geq \sqrt{n}$ for all $n$.

\subsubsection{Renling Jin}

Let $\underline{d}(A)=\liminf_{n\to\infty} A(n)/n$ be the lower asymptotic density of $A$, and let ${\mathcal P}$ be the set of primes. Clearly
    \[\forall A \subseteq {\mathbb N} \;\; (\underline{d}(A+{\mathcal P}) \geq f(\underline{d}(A))\]
for $f(x)=x$. What is the right $f$?

Using Pl\"{u}nnecke and $\underline{d}(3{\mathcal P}) = 1$ (due to Easterman, van der Corput, and possibly others independently) we get $f(x)=x^{2/3}$.

Audience: Can replace ${\mathcal P}$ with any $h$-basis and still have $f(x)=x^{1-1/h}$.

Note that Erd\H{o}s proved the existence of $A$ with $A(n)\sim \log n$ and $A+{\mathcal P} \sim {\mathbb N}$, so the primes are not a typical basis.

Audience: Can $x^{2/3}$ be improved assuming the Goldbach conjecture? Answer: Goldbach gives only $\sigma(4{\mathcal P})=1$, so no.

\subsubsection{Mel Nathanson}

Clarifying earlier problem. We say that two points in ${\mathbb R}^n$ are congruent if their difference is in ${\mathbb Z}^n$. Suppose that $k\subseteq{\mathbb R}^n$ is compact and for each $x\in {\mathbb R}$ there is a $y\in K$ such that $x\equiv y$.

\ \\

{\bf Theorem:} $A:= (K-K)\cap {\mathbb Z}^n$ is a finite set and generates the additive group ${\mathbb Z}^n$.

\ \\

For $n=1$, there is a $K$ that will give any set of generators that contains 0 and is symmetric about 0. For $n=2$, which sets of generators arise in this fashion? Specifically, is there a $K$ (compact and hitting every residue class modulo 1) such that $A \subseteq \{(x,y) {\rm col}on xy=0\}$? Even more specifically, is there a $K$ with $A=\{(0,0),(\pm1,0),(0,\pm1)\}$?

\ \\

\textbf{Note since the end of the conference: Renling Jin has solved this problem. Independently, and by different methods, Mario Szegedy has obtained a partial solution.}



\subsection{Speakers and Participants Lists}

\subsubsection{Speakers}

\bi

\item Hoi H. Nguyen <hoi@math.rutgers.edu>

\item Mariah E. Hamel <mhamel@math.uga.edu>

\item Steven J. Miller <Steven.J.Miller@williams.edu>

\item Alex Iosevich <iosevich@gmail.com>

\item Benjamin Weiss <blweiss@umich.edu>

\item Brooke Orosz <borosz@gc.cuny.edu>

\item Charles Helou <cxh22@psu.edu>

\item Craig Spencer <craigvspencer@gmail.com>

\item Gang Yu <yu@math.kent.edu>

\item Jaewoo Lee <jlee1729@hotmail.com>

\item Jonathan Sondow <jsondow@alumni.princeton.edu>

\item Julia Wolf <julia.wolf@cantab.net>

\item Kent Boklan <boklan@boole.cs.qc.edu>

\item Kevin O'Bryant <obryant@gmail.com>

\item Lan Nguyen <ltng@umich.edu>

\item Le Thai Hoang <leth@math.ucla.edu>

\item Li Guo <liguo@andromeda.rutgers.edu>

\item Mario Szegedy <szegedy@cs.rutgers.edu>

\item Mei-Chu Chang <changmeichu03@gmail.com>

\item Mel Nathanson <melvyn.nathanson@lehman.cuny.edu>

\item Mohamed El Bachraoui <MElbachraoui@uaeu.ac.ae>

\item Neil Lyall <lyall@math.uga.edu>

\item Peter Hegarty <hegarty@chalmers.se>

\item Renling Jin <jinr@cofc.edu>

\item Rishi Nath <rnath@york.cuny.edu>

\item Simon Griffiths <sg332@cam.ac.uk>

\item Urban Larsson <urban.larsson@yahoo.se>

\ei

\subsubsection{Participants}

\noindent

Adon Amira, New York

Bela Bajnok, Gettysburg College

Kent Boklan, Queens College (CUNY)

Richard T. Bumby, Rutgers University - New Brunswick

Mei-Chu Chang, University of California-Riverside

David Chudnovsky, Polynechnic University - NYU

Bianca De Souza, CUNY Graduate Center

Mohamed El Bachraoui, United Arab Emirates University

Simon Griffiths, University of Montreal

Li Guo, Rutgers University-Newark

Mariah Hamel, University of Georgia

Piper Harris, Princeton University

Derrick Hart, Rutgers University - New Brunswick

Tanya Haxhoviq, CUNY Graduate Center

Peter Hegarty, Chalmers University of Technology and University of Gothenburg

Charles Helou, Penn State Brandywine

Le Thai Hoang, UCLA

Alex Iosevich, University of Missouri

Renling Jin, College of Charleston

Nathan Kaplan, Harvard University

Walter O. Kravec, Stevens Institute of Technology

Shanta Laishram, University of Waterloo

Urban Larsson, Chalmers University of Technology and University of Gothenburg

Jaewoo Lee, Borough of Manhattan Community College (CUNY)

Ines Legatheaux Martins

Karl Levy, CUNY Graduate Center

Neil Lyall, University of Georgia

Steven J. Miller, Williams College

Rishi Nath, York College (CUNY)

Mel Nathanson, Lehman College (CUNY)

David Newman, New York

Hoi H. Nguyen, Rutgers University-New Brunswick

Lan Nguyen, University of Michigan

Kevin O'Bryant, College of Staten Island (CUNY)

Brooke Orosz, Essex County College

Gina-Louise Santamaria, Montclair State University

Steven Senger, University of Missouri - Columbia

Satyanand Singh, CUNY Graduate Center

Jonathan Sondow, New York

Craig Spencer, Institute for Advanced Study

Jacob Steinhardt

Mario Szegedy, Rutgers University-New Brunswick

Jonathan Wang

Benjamin Weiss, University of Michigan

Julia Wolf, Rutgers University-New Brunswick

Thomas Wright, Johns Hopkins University

Gang Yu, Kent State University

\newpage

\section{CANT Problem Sessions: 2010}

\subsection{Problem Session III: Friday, May 28th (Chair)}

\subsubsection{Nathanson: Classical Problems in Additive Number Theory}

N. G. de Bruijn had two papers:

\ben

\item \emph{On bases for the set of integers}, 1949.

\item \emph{On number systems}, 1956.

\een

Very few references to these papers. The second paper: he stated and solved a problem; in the first he stated a problem but neither he nor others could solve. Lately, however, these have become of interest to people in harmonic analysis.

These are related to the idea of complementing sets. Given a finite set $A$, can you find an infinite set $B$ such that $A \oplus B = \Z$? De Bruijn considered a slightly different problem, but in the same spirit. Given a family of sets $\{A_i\}_{i \in I}$ with $I = \N$ or $\{1,\dots,n\}$, we are interested in sets with the property that $\N_0 = \oplus_{i\in I} A_i$; in other words, every non-negative integer is of the form $\sum_{i\in I} a_i$ and $a_i \neq 0$ only finitely often. De Bruijn calls this a \emph{British number system}. Years ago 12 pence in a shilling, .... The British number system (pence, shillings, pounds) is the motivation for notation. If you have 835 pence that is 3 pounds, 9 shillings and 7 pence. The British number system is based on 12 and 20.

Using 12 pence is 1 shilling and 20 shillings is 1 pound. Take $\{g_i\}_{i \in I}$, $g_i \ge 2$, $G_0 - 1$, $G_1 = g_1$, $G_2 = g_1g_2$, $G_i = g_1g_2\cdots g_n$, $$A_n\ =\ G_{n-1} \ast [0,1,2,\dots, g_{i-1}) \ = \ G_{i-1} \ast [0, g_i),$$  and \bea A_1 & \ = \ & \{0, 1, 2, \dots, g_{1-1}\}  \nonumber\\ A_2 &=& G_1 \ast \{0,1,\dots, g_{2-1}\} \nonumber\\ A_3 &=& G_2 \ast \{0,1,2, \dots, g_{3-1}\} \nonumber\\ A_{n+1} &=& G_n \ast \N_0, \nonumber \eea where $$d \ast A \  = \ \{d a: a \in A\}.$$ Are there other sets? Yes. Let
$\{A_i\}_{i\in I}$ and $I = \cup_{j \in J} I_i$ with $I_j \cap I_{j'} = \emptyset$ for $j \neq j'$, $B_j = \sum_{i \in I_j} A_i$, $\{B_j\}_{j\in J}$. Comes down to choosing sequence of $g$'s to be prime numbers to get indecomposable sequence.

Consider a set $B$ of integers such that every $n\in \Z$ has a unique representation in the form $$n \ = \ \sum_{b \in B} \epsilon_b b$$ where $\epsilon b \in \{0,1\}$ and $\epsilon_b = 1$ finitely often. Let $A_i = \{0,b_i\}$, $\oplus_{i=1}^\infty A_i = \Z$.

Take set of powers of 2: $\{2^i\}_{i=0}^\infty$: get all non-negative integers. Suppose we look at $\{\epsilon_i 2^i\}_{i=0}^\infty$ where $\epsilon_i \in \{\pm 1\}$. Need infinitely many $+1$s and $-1$s to be a basis. No other condition necessary. Works if there are infinitely many $+1$s and infinitely many $-1$s. What infinite sets $B$ have this property? Want subset sums to give each integer once and only once. First thing can say is that we better not have everything even, so must have at least one odd integer in the set. Then de Bruijn proves something clever: not only at least one odd integer, but at most one odd integer. Was a conjecture of someone else, de Bruijn proves this conjecture.

Think about this for a minute. Exactly one odd number. If you are going to represent an even number it cannot have that odd number, and thus if divide all even numbers by 2 get another system of this form, so one of these and only one of these is divisible by 2 and not 4. By induction, see for every power of 2 there is one and only one number $x$ in this set such that $2^7||x$. We can thus write $b_i = d_n 2^i$ with $d_i$ odd. So $B$ comes from a sequence of odd numbers. Let's call this sequence of odd numbers $\{d_i\}_{i=1}^\infty$ \emph{okay}; sequence of odd numbers. In other words, it is an additive basis. Just restated the problem -- what sequences of odd integers are okay?

No one knows what sequences of odd numbers are okay. De Bruijn proved the following: suppose $\{d_1, d_2, d_3, \dots\}$ is an okay sequence; this is an okay sequence if and only if $\{d_2, d_3, d_4, \dots\}$ is okay. Can throw off any bunch -- do again. Can really screw around with an okay sequence -- can chop at any point, any finite sequence of garbage in the beginning. This is an interesting problem. I went this morning to MathSciNet to see what papers have referenced this paper of de Bruijn. There was a gap of about 50 years, but now relevant for something in harmonic analysis (they can't solve this problem, but it is in the same spirit as something they are interested in).

\subsubsection{Schnirelman}

When did additive number theory start? In 1930s Schnirelman proved that every even number is the sum of a large number of primes; he did this by proving a theorem about sumsets. Before this the results were beautiful (Fermat, Lagrange, Gauss, Hardy, Ramanujan, Littlewood, Vinogradov); Schnirelman had a general theorem about integers.Let $$A + B \ = \ \{a + b: a \in A, \ \ b \in B\}.$$ Counting function $$A(n) \ = \ \sum_{a \in A \atop 1 \le a \le n} 1.$$ Let's say $0 \in A\cap B$ to be safe. Defining $$\delta(A) \ = \ \inf \frac{A(n)}{n}.$$ Schnirelman proved $$\delta(A+B) \ \ge \ \delta(A) + \delta(B) - \delta(A)\delta(B).$$

Norwegians are funny -- go off to the mountains and come down with a great theorem. In WWI de Bruijn goes up to mountains and invents a sieve  method that allows him to prove things about Goldbach and Twin Primes. No one could understand the paper. Landau couldn't understand it, didn't try. Schnirelman did understand and used it to get his results, which made the result / method fashionable. Now people studied de Bruijn's paper. Landau's exposition in one of the seminal journals became the standard exposition for the de Bruijn sieve. Same thing happened with Selberg. WWII started, he was captured by Germans, released if promised not to stay in Oslo, went to family home and proved results on zeros of $\zeta(s)$.

Could also look at $$\delta_L(A) \ = \ \lim_{n\to \infty} \frac{A(n)}{n}.$$ Say $A \sim B$ if there is an $N$ such that for all $n \ge N$ we have $n \in A$ if and only if $n \in B$. Embarrassment: if every element is even then all sums even, must be careful. If $\delta_L(hA) > 0$ for some $h$ then if $d = \gcd(A)$ and $0 \in A$ there there is an $h_0$ such that $h_0 A \sim d \ast N_0$. First time appears is in a paper with John C. M. Nash (the son).

Let $0\in A$ and $d_L(hA) = 0$ for all $h \ge 1$. Assume $$A \ \subseteq \ 2A \ \subseteq \ 3A \ \subseteq \ 4A \ \subseteq \ \cdots \ \subseteq \ hA \ \subseteq \ \cdots.$$ If any set has positive density then get all integers from some point onward. Maybe in this case some infinite case appears. The question is: take a set of non-negative integers containing 0 such that all of these sets have asymptotic density zero. Get an increasing sequence. As add to itself more and more times, does any structure appear? Can you say something that interests your friends mathematically about this? Is there anything that must happen?

\subsubsection{Alex Kontorovich}

Not convinced problem is difficult, but we haven't made progress. The question is the additive energy in ${\rm SL}(2,\Z)$. This means that we take elements $\gamma_1, \gamma_2, \gamma_3, \gamma_4$ in ${\rm SL}(2,\Z)$ in a ball $B_N$ and want to know how many there are such that $\gamma_1+\gamma_2 = \gamma_3 + \gamma_4\}$. The number of points in a ball (4 variables, 1 quadratic equation) gives $cN^2$.

We have 16 variables (unknowns), 4 quadratic equations (the determinants equaling 1) and 4 linear equations. Want an upper bound of the form $\ll N^{4+\epsilon}$. Have a trivial lower bound of $N^4$. Easy thing to prove is $N^5$ for the following reason. Let $$\twocase{\eta(\omega) \ = \ \#\{\gamma_1 \gamma_2 \in B_N({\rm SL}_2): \gamma_1 -\gamma_2 = \omega\} \ = \ }{N^2}{if $\omega = 0$}{N^{1+\epsilon}}{if $\omega \neq 0.$}$$ What about  non-trivial bounds? See \begin{center} \bburl{http://arxiv.org/pdf/1310.7190v1.pdf} \end{center} for more on this problem.

\subsubsection{Peter Hegarty}

This is a problem on Phase Transitions inspired by Hannah Alpert's talk. Let $G$ be an Abelian group, $A$ a set of generators, everything infinite. Have $S(1)$, $S(2)$, $S(3)$, $\dots$, $S(r)$, $\dots$, where the sequence is $\infty, \infty, \dots, \infty, (n), 0, 0, \dots$ (where we may or may not have the $n$ term). We should be able to compare the sizes of infinite, i.e., their measure. Suppose $G$ is a compact Abelian group, such as the circle, and let $A$ be a measurable set. Want to look at Lebesgue (or Haar) measure of the sets: $\mu(S(1))$, $\mu(S(2))$, et cetera. The sequence should be unimodal (regular).

David Neumann looked at something similar. For finite groups looking at the sizes, did a lot of computations with different groups and generating sets. Did find an example where it wasn't the case, but typically do have unimodality. Hegarty conjectured that for any finite group (not necessarily Abelian) can always find a set of generators such that the sequence is unimodal.


\subsection{Problem Session IV: Saturday, May 29th}

\subsubsection{Peter Hegarty}

Let $A \subset \N$, $$r(A,n) \ = \ \#\{(a_1,a_2): a_1 + a_2 = n\}.$$ What sequences of non-negative integers can be asymptotic representation functions? Of course there are restrictions if start from 0. Obviously only one way to represent 0 (0+0). Given a sequence of numbers, want the sequence to equal $r(A,n)$ starting at some point. Assuming Erd\"os-Turan, cannot be bounded and simultaneously not have infinitely many zeros.

Comment from Nathanson: Matter of choice whether take $r(A,n)$ or the function $$\mathfrak{r}(A,n) \ = \ \#\{(a_1,a_2): a_1 + a_2 = n, a_1 \le a_2\}.$$

More generally, say $|S| = \infty$ and $S \subseteq G$, $A$ is an asymptotic basis (of order 2) for $S$ if $S {\subset \atop \sim} A + A$ (up to finite sets). Let $$r_A(s) \ = \ \#\{(a_1,a_2): a_1+a_2 = s\}.$$ We have $r_A:S \to \N \cup \{\infty\}$. Hardest problem is what we had earlier.


\subsubsection{Mel Nathanson}

Erd\"os-Renyi Method: Let $$\Omega \ = \ \{{\rm all\ sequences\ of\ non-negative\ integers}\}.$$ Let $0 \le p(n) \le 1$ for $n =0, 1, 2, \dots$. Then there exists a probability measure $P_r$ on $\Omega$ such that $$P_n(E_n) \ := \ P_n(\{A \in \Omega: n \in A\}) \ = \ p(n)$$ then the events $E_n$ are independent.

If choose $p(n)$ to be something like a logarithm over a power of $n$, say $\alpha\frac{\log^\beta n}{n^\gamma}$ with $1/3 < \gamma \le 1/2$ -- want a result that doesn't use any probability. If put this probability measure on the sequence of integers, then if $A \subset \N_0$ with $A+A \sim \N_0$ and $S(n) = \{a \in A:n-a \in A\}$, then for $m\neq n$ we have $$|S(m) \cap S(n)| \ \le \ \frac{2}{3\gamma - 1}$$ for all but finitely many pairs of integers.

Below is an example of where this result was used. An asymptotic basis means every number from some point onward can be written as $a+a'$ with $a,a'\in A$. An asymptotic basis $A$ is minimal if no proper subset of $A$ is an asymptotic basis. This means we have the set of integers with the property that if throw away any number then all of a sudden infinitely many numbers cannot be represented. Came up in an attempt to construct a counter-example to the Erd\"os-Turan conjecture. Not every asymptotic basis contains a minimal basis. There is a theorem that says that if $r_A(n) \to \infty$ and $|S_A(m)\cap S_A(n)| = O(1)$ then $A$ contains a minimal asymptotic basis.

Theorem: If have a sequence with $r_A(n) > c\log n$ for some $c > 1 / \log(4/3) \approx 3.47606$ and $n \ge n_0$ then $A$ contains a minimal asymptotic basis.


\subsection{Speaker List}

\bi

\item Hannah Alpert, University of Chicago
\item Paul Baginski, Universite Claude Bernard Lyon, France
\item Gautami Bhowmik, Universite Lille, France
\item Kent Boklan, Queens College (CUNY)
\item Mei-Chu Chang, University of California-Riverside
\item Scott Chapman, Sam Houston State University
\item Brian Cook, University of British Columbia
\item David Covert, University of Missouri
\item Aviezri Fraenkel, Weizmann Institute of Science, Israel
\item John Friedlander, University of Toronto
\item John Griesmer, University of British Columbia
\item Sinan Gunturk, Courant Institute, NYU
\item Peter Hegarty, Chalmers University of Technology and University of Gothenburg
\item Charles Helou, Penn State Brandywine
\item Alex Iosevich, University of Missouri
\item Renling Jin, College of Charleston
\item William J. Keith, Drexel University
\item Alex Kontorovich, Institute for Advanced Study
\item Brandt Kronholm, SUNY at Albany
\item Urban Larsson, Chalmers University of Technology and University of Gothenburg
\item Jaewoo Lee, Borough of Manhattan Community College (CUNY)
\item Zeljka Ljujic, CUNY Graduate Center
\item Neil Lyall, University of Georgia
\item Steven J. Miller, Williams College
\item Rishi Nath, York College (CUNY)
\item Mel Nathanson, Lehman College (CUNY)
\item Lan Nguyen, University of Michigan
\item Kevin O'Bryant, College of Staten Island (CUNY)
\item Alex Rice, University of Georgia
\item Steven Senger, University of Missouri
\item Jonathan Sondow, New York
\item John Steinberger, Institute for Theoretical Computer Science, Tsinghua University

\ei

\newpage

\section{CANT Problem Sessions: 2011}

\subsection{Problem Sessions}

There was an issue with my computer and the original file was lost for 2011; the items below are restored from earlier copies, though I have lost who spoke on what day and thus have run this as one entry.

\subsubsection{Seva Lev}

\noindent \textbf{Problem:} Let $A \subset \mathbb{F}_2^n$, $p \in \mathbb{F}_2[x_j]_{j=1}^n$, and for all $a, b \in A$ if $a\neq b$ then $p(a+b)=0$. Does this imply that $p(0) = 0$?\\

For example, if $A = \{a,b\}$ then $p(a+b)=0$ does not imply $p(0) = 0$.

If $A$ is large and the degree of $p$ is small, what is true? For a given $p$, how large must $|A|$ be for this to be true? We have the following: \bea \deg P & \ \ \ \ & {\rm need} \nonumber\\ 0 & & |A| \ge 2 \nonumber\\ 1 & & |A| \ge 3 \nonumber\\ 2 & & |A| \ge n+3 \nonumber\\ 3 & & |A| \ge 2n \nonumber\\ \le \left(frac12 o(1)\right)n & & ???. \nonumber \eea


\subsubsection{Giorgis Petridis}

P-R: $D_2 \ge 1$ implies that there exists $v_0$ vertex disjointed paths of length 2 in $G$.\\

\noindent \textbf{Problem:} What can be said when $D_2 \ge k \in \Z$?\\

Guess: there exist $v_0$ vertex disjoint trees in $G$ each having at least $k_i$ vertices in $V_i$. Note: there is an example which shows that one cannot hope to prove this guess using max flow - min out. Guess confirmed in $k=|V_0|=2$ by Petridis.

\subsubsection{Mel Nathanson}

Believe the following is an unsolved problem by Hamidoune (he proposed it and no one has solved it):\\

\noindent \textbf{Problem:} Let $G$ be a torsion free group, $G \neq \{e\}$. Let $S$ be a finite subset of $G$, $e \in S$, $$\kappa_k(S) \ = \ \min \left\{ |XS| - |X|: {\rm finite sets}\ X \subset G, |X| \ge k\right\}.$$ Hamidoune conjectured that there is an $A \subset G$ with $|AS|-|A| = \kappa_k(S)$ and $|A|=k$. \\

True for $k=1$, unknown for $k \ge 2$. It is true for ordered groups. As every free abelian group of finite rank can be ordered, true here. In general for $k=2$ still unknown.


\subsubsection{Matthew DeVos}

\noindent \textbf{Problem:} Let $G$ be a multiplicative group, $S \subset G$ a finite set, and set $$\Pi(S) \ = \ \{s_1\cdots s_k: s_i \in S, s_i = s_j \Longleftrightarrow i=j\} \cup\{1\}.$$ Not allowed to use an element multiple times. Conjecture: there is a $c>0$ such that for every group $G$ and set $S \subset G$ there exists $H \subset G$ with $|\Pi(S)| \ge |H| + c|H| \cdot |S\setminus H|^2$.\\

True with $c = 1/64$ when $G$ is abelian.

\subsubsection{David Newman}

\noindent \textbf{Problem:} How many partitions are there where no frequency is used more than once?\\

For example, the partitions of 4 are $\{4\}$, $\{3,1\}$, $\{2,2\}$, $\{2,1,1\}$ and $\{1,1,1,1\}$. The ones that are okay are all but $\{3,1\}$. The problem here is that the two decompositions each occur just once: we have one 3 and one 1.


\subsubsection{Steven J. Miller, Sean Pegado, Luc Robinson}

\noindent \textbf{Problem:} For each positive integer $k$, consider all $A$ such that $|kA+kA| > |kA - kA|$ and $1 \in A$ (for normalization purposes). Let $C_k$ be the smallest of the largest elements of such $A$'s. What can you say about the growth of $C_k$?\\

$C_1 = 15$, $C_2 = 31$, $\dots$.


\subsubsection{Speaker unremembered}

\noindent \textbf{Problem:} Assume that you have $A, B$ in a general group and $|AB| < \alpha |A|$ and $|AbB| \le \alpha|A|$ for all $b\in B$. Does there exist an absolute $c$ such that $X \subset A$ then $|X B^h| \le \alpha^{ch} |X|$?\\

Rusza showed that if you have $|A+B_j| \le \alpha_j |A|$ for $j = 1, 2$, then there is an $X$ such that $|X+B_1+B_2| \le \alpha_1 \alpha_2 |X|$.\\

\noindent \textbf{Problem:} Is there a prescription for $X$ given that Rusza's theorem shows the existence of $X$.


\subsubsection{Speaker unremembered}

\noindent \textbf{Problem:} Let $\mathcal{B}$ be a partition of $n$. Consider the partition where $c_1 +\cdots + c_k = n$, and $1^{d_1} 2^{d_2} \cdots n^{d_n}$. The $d_i$'s are the number of the $c_j$'s and $d_1 + \cdots + d_n = m$.  Consider $$\sum_{\mathcal{B} \in \mathcal{P}(n)} \ncr{n}{c_1,\dots,c_k} \ncr{n+m+1}{n+1,d_1,\dots,d_n} \left[\frac1{m+n+1}\right].$$ What can be said?\\

Try putting in an $r^n$ and summing over $n$. Maybe this is a holomorphic part of a non-holomorphic Maass form.


\subsubsection{Peter Hegarty}

\noindent \textbf{Problem:} Consider the least residue of $n$ modulo $q$, denoted $[n]_q$, which is in $\{-q/2, \dots, q/2\}$. Want a function from $\pi: \{1,\dots,27\}$ to itself (a permutation, so 1-1) with the property that given any $a,b,c$ not all equal with $|[a+c-2b]_{27}| \le 1$ then $|[\pi(a)+\pi(c)-2\pi(b)]_{27}| \ge 2$.\\

Motivation: replace 27 with $n$, ..., have a permutation avoiding a progression. Conjecture that a permutation of $\Z_n$ exists for every $n$ sufficiently large.


\subsubsection{Speaker unremembered}

\noindent \textbf{Problem:} Define $h: \{1,\dots,N\} \to \Z/N\Z$; call it a partial homomorphism if  it a bijection such that whenever $a,b, ab \in \{1,\dots,N\}$ then $h(ab) = h(a) + h(b) \bmod N$. Does such a function exist for all $N$? \\

Have built by hand for all $N$ up to 64?


\subsubsection{Steven Senger}

The basic idea is that an additive shift will destroy multiplicative structure. Given a large, finite set, $A \subset \mathbb{N}$, suppose that $|AA|$ = $n$. We know that there exists no generalized geometric progression, $G$, of length $c_1n$, such that $| (AA+1) \cap G| \geq c_2n$, where $c_1$ and $c_2$ do not depend on $n$. The question is, given the same conditions on $A$, do there exist sets $E, F \subset \mathbb{N}$, such that the following hold for $c_3, c_4$ independent of $n$, and $\delta > 0$:

\begin{itemize}
\item $|E|,|F| \geq n^\delta$
\item $|EF| = c_3n$
\item $|(AA+1) \cap EF| \geq c_4n$
\end{itemize}
Even partial results would be interesting to me. Also, considering the problem over $\mathbb{R}$ would be interesting to me.


\subsubsection{Urban Larsson}
2 pile Nim can be described as the set of moves on a chessboard made by a rook, moving only down and left. Players take turns moving the rook, and the person to move it to the lower-left corner is the winner. The set of legal moves is defined to be
$$\{(0,x),(x,0) \}.$$
In this case, the positions which guarantee victory following perfect play, or {\it p-positions} are along the diagonal. That is, the player who consistently moves the rook to the diagonal will eventually win.

In Wythoff Nim, the piece is replaced by a queen, and the diagonal move is added. The set of legal moves for Wythoff Nim is
$$\{ (0,x),(x,0),(x,x) \}.$$
This game has p-positions close to the lines of slope $\phi$ and $\phi^{-1}$, where, $\phi$ denotes the golden ratio. For example, the points $(\lfloor \phi x\rfloor, \lfloor \phi^2 x \rfloor)$ are p-positions in Wythoff Nim.

Now, adjoin the multiples of the last possible p-positions from Wythoff Nim which are not in Wythoff Nim, namely the multiples of the knight's move. The legal moves of the new game are
$$\{0,x),(x,0),(x,x),(x,2x),(2x,x)\}.$$
The p-positions for this game appear to split along lines of slopes nearly 2.25 and 1.43. Why?


\subsubsection{Thomas Chartier}
Let $n,k \in \mathbb{N}$, and $p=nk+1$ be prime. Exclude 1 and 2. Fixing $n$ does there exist a $k$ such that
$$1^k, 2^k, 3^k, \dots, n^k$$
are distinct mod $p$?
The conjecture is that such a $k$ exists for every non-trivial $n$.


\subsubsection{Mel Nathanson}
Recall the classical sum-product problem of Erd\H os. Given a large set of positive integers, $A \subset \mathbb{N}$, either the set of sums or the set of products should be large. The conjecture is that, for such an $A$, with $c$ independent of $n$, for any $\epsilon >0$,
$$\max\{|A+A|,|AA|\}\geq cn^{2-\epsilon}.$$

\subsection{Speaker List}

\begin{enumerate}

\item Paul Baginski, Universite Claude Bernard Lyon, France
\item Mei-Chu Chang, University of California-Riverside
\item Scott Chapman, Sam Houston State University
\item Jonathan Cutler, Montclair State University
\item Matthew DeVos, Simon Fraser University
\item Aviezri Fraenkel, Weizmann Institute of Science, Israel
\item Peter Hegarty, University of Gothenburg, Sweden
\item Charles Helou, Penn State Brandywine
\item Jerry Hu, University of Houston - Victoria
\item Alex Iosevich, University of Missouri
\item Renling Jin, College of Charleston
\item Nathan Kaplan, Harvard University
\item Mizan R. Khan, Eastern Connecticut State University
\item Omar Kihel, Brock University, Canada
\item Alex Kontorovich, SUNY at Stony Brook
\item Urban Larsson, University of Gothemburg, Sweden
\item Thai Hoang Le, Institute for Advanced Study
\item Vsevolod Lev, University of Haifa, Israel
\item Zeljka Ljujic, CUNY Graduate Center
\item Neil Lyall, University of Georgia
\item Steven J. Miller, Williams College
\item Rishi Nath, York College (CUNY)
\item Mel Nathanson, Lehman College (CUNY)
\item Hoi H. Nguyen, University of Pennsylvania
\item Lan Nguyen
\item Sean Pegado, Williams College
\item Giorgis Petridis, University of Cambridge
\item Luc Robinson, WIlliams College
\item Steve Senger, University of Missouri
\item Jonathan Sondow, New York
\end{enumerate}

\newpage

\section{CANT Problem Sessions: 2012}


\subsection{Problem Session I: Tuesday, May 22nd (Chair Renling Jin)}

\begin{itemize}

\item \emph{From Renling Jin, jinr@cofc.edu:} Define a subset of the natural numbers $B$ to be an essential component if for all $A \subset \N$, $\sigma(A+B) > \sigma(A)$ if $0 < \sigma(A) < 1$. $B$ is an extraordinary component if $$\liminf_{\sigma(A) \to 0}\frac{\sigma(A+B)}{\sigma(A)} \ = \ \infty.$$ Here $$\sigma(A) \ = \ \inf_{x \ge 1} \frac{A(x)}{x}.$$

    Ruzsa conjectured that every essential component is an extraordinary component.

    What are the essential components we know? If $$B \ = \ \{k^2: k \in \N\}$$ then $$\sigma(A+B) \ \ge \ \sigma(A)^{1 - 1/4}$$ since $B$ is a basis of order four. We get $$\frac{\sigma(A+B)}{\sigma(A)} \ \ge \ \frac{1}{\sqrt[4]{\sigma(A)}}.$$ Similar for cubes or $k$-powers.

\item \emph{From Steven J. Miller, sjm1@williams.edu:} We say a set $A$ is a More Sums Than Differences Set, or an MSTD set, if $|A+A| > |A-A|$, where \bea A+A & \ = \ & \{a_i + a_j: a_i, a_j \in A\} \nonumber\\ A-A & \ = \ & \{a_i - a_j: a_i, a_j \in A\}. \nonumber\ \eea As addition is commutative and subtraction is not, it's expected that `most' sets are difference dominated; however, Martin and O'Bryant proved that a positive percentage of sets are sum-dominated. There are explicit constructions of infinite families of sum-dominant sets. Initially the best result was a density of $n^c 2^{n/2} / 2^n$, then $1/n^4$ (or $1/n^2$), and now the record is $1/n$ (where our sets $A$ are chosen uniformly from subsets of $\{0, 1, \dots, n-1\}$). Can you find an `explicit' family that is a positive percentage.

\item \emph{From Urban Larsson, urban.larsson@yahoo.se:} Let $A = \{0, 1, 3, 4, \dots\}$ for a set that avoids arithmetic progression, thought to be best set to avoid arithmetic progression but not (comes from a greedy construction). Equivalence with a base 3 construction: $A = \{0, 1, 10, 11, 100, 101, \dots\}$ gives $A((3^n+1)/2) = 2^n$, where $A(n) = \#\{i\in A\mid i < n\}$. Hence, for all $n$, $A(n) < Cn^{\log 2/\log 3}\approx n^{2/3}$. Study impartial heap games. Is it possible to find a game such that the P and N-positions correspond to the numbers in this construction? (A position is in N if and only if the first player wins.) In some sense such that: \begin{center}\begin{tabular}{ccccccc}
                                                     P& P &N  & P & P & N & N \\
                                                     \hline
                                                     0& 1 & 2 &3  & 4 & 5 & 6 \\
                                                  \end{tabular}\end{center}
We rather use three heaps of sizes in three-term arithmetic progression. A legal move is to erase the largest pile and then to announce one of the smaller piles as the new largest pile. Notation $(x,y)$, where $x$ is the number of tokens in the smallest heap and $y$ in the second smallest. In the table below, the first entry is the outcome, the second is the position, the third is the Grundy value, and the fourth are the options.
\bea P\ \ \  (0,1) & & 0 \ \ \ \nonumber\\
 N\ \ \  (0,2) & & 1 \ \ \ (0,1) \nonumber\\
 P\ \ \  (0,3) & & 0 \ \ \ (1,2)\nonumber\\
 N\ \ \  (0,4) & & 0 \ \ \ (0,2), \ (2,3) \nonumber\\
 N\ \ \  (1,2) & & 1 \ \ \ (0,1)  \nonumber\\
 N\ \ \  (1,3) & & 0 \ \ \ (1,2), \nonumber\\
 P\ \ \  (1,4) & & 0 \ \ \ (0,2), \ (2,3) \nonumber\\
 N\ \ \  (2,3) & & 2 \ \ \ (0,1), \ (1,2) \nonumber\\
 N\ \ \  (2,4) & & 3 \ \ \ (1,2), \ (0,2), \ (2,3)\nonumber\\
 N\ \ \  (3,4) & & 0 \ \ \ (1,2), \ (0,2), \ (2,3)\nonumber \eea

The P positions (Grundy value 0) have both lower heap sizes in the set $A$. The N positions have Grundy values $> 0$, defined as the \emph{minimal exclusive} of the Grundy values of the options. What are they? Is it possible to extend the game by \emph{adjoining moves} to obtain limsup$A(n)/n^{\log 2/\log 3}~=~\infty$? The game generalizes to $k$-term arithmetic progressions and the Sidon-condition for example.

How do we extend such games? We need a general definition for the family of games. A \emph{ruleset} is a set of finite sets of positive integers. From a \emph{position} consisting of a set $S$ of non-negative integers, choose one of the numbers $s\in S$ and a set $M$ of numbers from the given ruleset. The next position, which is a set of nonnegative numbers, is $\{s-m\mid m\in M\}$, provided max$M\le s$. We get a recursive definition of the set $A$ which determines the P-positions for a given $M$. A position $S$ is in P if and only if $S\subset A$. That is $S$ is in N iff $S\cap A\neq\emptyset$. In this sense we can abuse notation and regard $A$ as the set of ``P-positions''. A game extension of $M$ is $M\cup M'$, for $M'$ a set of finite sets of nonnegative numbers. For our game the set $M$ is $M=\{\{d,2d\}\mid d>0\}$. One first example of a game extension is $M = \{\{d,2d\}\mid d>0\}\cup \{\{1\}\}$. Question: does the set $A$ become less dense for this game than for our original AP-avoiding game?

\end{itemize}


\subsection{Problem Session II: Wednesday, May 23rd (Chair Steven J Miller)}

\begin{itemize}

\item \emph{From Steven J Miller, sjm1@williams.edu:} We investigated  in \hfill \\ \bburl{http://arxiv.org/pdf/1109.4700v2.pdf}) properties of $|A+A|$ and $A+A$ as $A$ varies uniformly over all subsets of $\{0, 1, \dots, n-1\}$. How does the behavior change if we change the probability of choosing various $A$'s (see for example my work with Peter Hegarty: \hfill \\ \bburl{http://arxiv.org/pdf/0707.3417v5}).

    Another related problem is to `clean-up' the formula we have for the variance. This involves sums of products of Fibonacci numbers -- can the answer be simplified?

    What about the expected values of $2kA$ versus $kA - kA$.

\item \emph{From Ryan Ronan, ryan.p.ronan@gmail.com:} Earlier today I discussed joint work on generalized Ramanujan primes, \hfill \\ \bburl{http://arxiv.org/pdf/1108.0475}. One natural question is whether or not for each prime $p$ there is some constant $c_p$ such that $p$ is a $c_p$-Ramanujan prime.

    Another question is the distribution of $c$-Ramanujan primes among the primes, in particular the length of runs of these and non-these. It can take awhile for the limiting behavior of primes to set in; it's dangerous to make conjectures based on small sized data sets. Are the calculations here sufficiently far enough down the number line to have hit the limiting behavior? For a related question, perhaps the Cramer model is not the right model to use to build predictions, and instead we should use a modified sieve to construct `random primes'. It would be worthwhile to do so and see what happens / what the predictions are.

\item \emph{From Steven Senger, senger@math.udel.edu:} Have a subset $A$ of a finite field $\mathbb{F}_q$ satisfying for all $\epsilon$ and $\delta$ positive  (1) $|A|\ |AA| \ge q^{3/2 + \epsilon}$, (2) $|AA| \le q^{1-\delta}$. For all generalized geometric progressions $G$ with $|G| \approx |AA|$ we have $|(AA+1)\setminus G| \ge q^{\delta}$. Can reduce the size constraint (1)? Can we increase the size of $|(AA+1)\setminus G| \ge q^{\delta}$?

\item \emph{From Kevin O'Bryant, obryant@gmail.com:} How far out can you go $\{x_1$, $x_2$, $x_3$, $x_4$, $\cdots\}$ such that the first four are in the first four quadrant, the first nine in the first nine subdivisions ($3 \times 3$), the first 16 in the first $4 \times 4$ and so on.... We know this can't go on forever, violates Schmidt.

    The discrepancy of the sequence $\{x_i\}$ is $${\rm Disc}(\{x_i\}_{i=1}^d) \ = \ \sup_R \left|\frac{\#\{x_i \in R\}}{d} - A(R) \right|.$$ We have ${\rm Disc}(\{x_i\}_{i=1}^d) \ge C \frac{\log d}{d}$. If we spread the points too well, the discrepancy gets very low.

    Let me rephrase -- I strongly believe that this logarithmic factor will kill this arrangement.

\end{itemize}


\subsection{Problem Session III: Thursday, May 24th (Chair Alex Iosevich)}

\begin{itemize}

\item \emph{From Jerry Hu, HuJ@uhv.edu:}

This problem is related to Nathanson's talk ``The Calkin-Wilf tree and a forest of linear fractional transformations'' from Tuesday.  Recall the form of the Calkin-Wilf tree, where we have:
\begin{eqnarray}
   &  \frac{a}{b} &     \nonumber \\
    \swarrow & \ & \searrow  \nonumber \\
    \frac{a}{a+b} \ \ \ \ \ \ & \ & \ \ \ \ \ \ \frac{a+b}{b}  \nonumber \\
\swarrow \ \ \ \ \searrow \ \ \ \ \  & \ & \ \ \ \ \ \swarrow \ \ \ \ \searrow \nonumber \\
\frac{a}{2a+b} \ \ \ \ \frac{2a+b}{a+b} & \ & \frac{a+b}{a+2b} \ \ \ \ \frac{a+2b}{b}.\nonumber
\end{eqnarray}
When $a$ and $b$ are both initialized as 1, each positive rational number appears on the tree exactly once.

The question is: how can we generalize this?  More specifically, do there exist other trees of the form
\begin{eqnarray}
  &z&     \nonumber \\
    \swarrow & \ & \searrow  \nonumber \\
     L(z) & \ & R(z) \nonumber
\end{eqnarray}
in which every positive rational number appears exactly once?  Can we find all, or any, nontrivial functional pairs $L, R$ such that this condition holds?

\item \emph{From  Nathan Kaplan, nathanckaplan@gmail.com: } Here is a problem about counting lines among points in
$\mathbb{F}_3^n$.  I will give two different kinds of motivation for
why someone might be interested in this.

The card game SET is played with $81$ distinct cards, each of which
has four attributes (number, color, shading, and shape), where each
attribute has three possibilities.  We can identify a card with a
$4$-tuple $(x_1, x_2, x_3, x_4)$, where each $x_i \in \mathbb{F}_3$.
The game is played by collecting sets.  A set is a collection of three
cards $(x,y,z)$ such that for each of the four attributes each card is
the same or all three cards are different.  It is equivalent that the
vectors in $\mathbb{F}_3^4$ represented by our three cards take the
form $(x,y,-(x+y))$, or equivalently, $(x,x+d,x+2d)$.  Therefore, we
see that what we are looking for is a three term arithmetic
progression in $\mathbb{F}_3^n$.  In $\mathbb{F}_3^n$ a $3$-term AP is
equivalent to a line.  A set of vectors with no $3$-term AP is called
a cap set.  The cap set problem asks, ``What is the maximum size of a
cap set in $\mathbb{F}_3^n$?''.  This problem is very hard and has
been well-studied.  Exact answers are known only for $n\le 6$.  We
note that for $n=3$ the cap set problem is equivalent to asking for
the maximum number of SET cards one can have so that there is no set
among them.  The answer to this is $20$ and an argument is given in
the paper \emph{The Card Game Set} by Benjamin Davis and Diane
MacLagan.

There is a related problem motivated by SET which does not seem to
have appeared in the literature.  The game is usually played by
dealing out $12$ cards.  We know that it is possible to have no sets
at all, but we could ask for the largest number of sets which could
occur among $12$ cards.  I can show that this is $14$, but the
argument is sort of ad hoc and not so satisfying.  I have not found
anything written before about the following question.  What is the
maximum number of lines that $m$ points in $\mathbb{F}_3^n$ can
contain?  Note that any two points determine a unique line, so if a
set contains many lines, then it determines few lines.  Equivalently
we could ask for the minimum number of lines determined by $m$ points
in $\mathbb{F}_3^n$.  This question is very general and includes the
cap set problem as a subcase.  This is because the number of lines
contained in a subset of $\mathbb{F}_3^n$ determines the number of
lines contained in its complement, so if we know the maximum number of
lines among any collection of $m$ points for all $m$, then we also
know the minimum number of lines among $m$ points.

Here is the actual problem I am asking.  In the argument for the
maximum number of lines among $12$ points in $\mathbb{F}_3^4$ is $14$,
it is clear that the maximum number of lines among $12$ points in
$\mathbb{F}_3^n$ is $14$ for any $n\ge 3$.  That is, if we want lots
of lines, the best thing that we can do is to put our points into the
smallest possible dimensional subspace that can contain them.

\begin{conj}
Fix $m\ge 0$ and let $d = \lceil \log_3(m) \rceil$.  For any $n\ge d$,
the maximum number of lines contained among $m$ points in
$\mathbb{F}_3^n$ is equal to the maximum number of lines contained
among $m$ points in $\mathbb{F}_3^d$.
\end{conj}

I think that this is probably true and that the proof for it is
probably easy.  One could also ask similar questions for
$\mathbb{F}_q^n$ for other $q$.

Here is some extra motivation that the cap set problem is interesting.
 Tic-Tac-Toe on a $3\times 3\times 3$ board can never end in a draw no
matter how many moves are made by each player.  This is the first case
of a more general phenomenon, the Hales-Jewett Theorem.  Given $k$,
there exists a $d$ such that Tic-Tac-Toe on a $k\times \cdots k =
[k]^n$ board (where it takes $k$ in a row to win) cannot end in a draw
no matter how many times each player moves.  A more precise statement
is that for large enough $n$, either a set or its complement must
contain a combinatorial line.  I won't define exactly what a
combinatorial line is, but it is a slightly more restrictive condition
than a Tic-Tac-Toe line, which is slightly more restrictive than the
type of line described above in the discussion of SET.

A few years ago, the initial Polymath project organized by Tim Gowers
was focused on giving a combinatorial proof of the Density
Hales-Jewett Theorem.  The only previous proof of this theorem
involved arguments from ergodic theory.  Let $c_{n,k}$ be the largest
number of points of $[k]^n$ which does not contain a combinatorial
line.  Let $c'_{n,k}$ be the largest number of points of $[k]^n$ which
does not contain a geometric line (you can think of this as a
Tic-Tac-Toe line.  These are called Moser numbers.  Finally, let
$c''_{n,k}$ be the largest number of points of $[k]^n$ without a line
of the type described above.  Clearly $c''_{n,k} \le c'_{n,k} \le
c_{n,k}$.

\begin{thm}[Density Hales-Jewett]
Fix $k\ge 1$.  Then
$$\lim_{n \to \infty} \frac{c_{n,k}}{n^k} = 0.$$
\end{thm}

This result is important in understanding the growth of cap sets.  The
Polymath project also proved the best known lower bound for $c_{n,k}$.
It is quite difficult to compute these numbers in general, even for
small $k$.  We mentioned above that $c''_{4,3} = 20$ and it is also
known that $c''_{5,3} = 45$ and that $c''_{6,3} = 112$.  This last
statement determines the maximum number of lines among $3^6-112$
points in $[3]^6$, for example.  The Polymath project also determined
more values of $c_{n,3}$ and $c'_{n,3}$ than previously known.

Since so much work has gone into understanding large subsets of
$[k]^n$ with no lines, it seems reasonable to study collections of
points which contain the largest possible number of lines.

\end{itemize}


\subsection{Problem Session IV: Friday, May 25th (Chair Kevin O'Bryant)}

The following papers are relevant for the problems proposed by Steven Miller.

\begin{itemize}

\item \bburl{http://arxiv.org/abs/1107.2718}

\item \bburl{http://arxiv.org/abs/1008.3204}

\item \bburl{http://arxiv.org/abs/1008.3202} (the gap paper referenced below is in preprint stage, but available upon request).

\item \bburl{http://www.emis.de/journals/INTEGERS/papers/j57/j57.pdf} (Hannah Alpert).

\end{itemize}

Proposed problems.

\begin{itemize}

\item \emph{From Steven J. Miller, sjm1@williams.edu:} The following problems are related to Zeckendorf decompositions. {{\rm col}or{red}Many of these are currently being studied by my summer REU students in the Williams 2012 SMALL program. If you are interested in working on these, please email me at sjm1@williams.edu.}

    \subitem $\diamond$ We know every number has a unique Zeckendorf decomposition, and appropriately localized the number of summands converges to being a Gaussian. What happens if we have a decomposition where some integers have multiple representations? What if there are some integers that have no representations? Instead of counting the total number of summands, what if you just count how many of each summand one has (so in decimal 4031 wouldn't count as $4+0+3+1$ but $1+0+1+1$).

    \subitem $\diamond$ We have formulas for the limiting distribution of gaps between summands of Fibonacci and some generalized Fibonacci sequences. Try to find formulas for general recurrence relations as a function of the coefficients of the relations. Do this for the signed Fibonacci decomposition (see Hannah Alpert's paper; can we generalize signed distributions to other recurrence relations). What about the distribution of the largest gap (that should grow with $n$ for numbers between $H_n$ and $H_{n+1}$). If we appropriately normalize it, does it have a nice limiting distribution?

\item \emph{From Mizan Khan, 	khanm@easternct.edu:}
Let $$\mathcal{H}_n:=\lbrace (x,y)\in \mathbb{Z}\times\mathbb{Z} : xy \equiv 1 ~(\text{mod } n), 1\leq x,y\leq n-1)\rbrace.$$
Consider the convex closure of $\mathcal{H}_n$-- what can we say about the number of vertices in this convex closure? Let $v(n)$ be the number of vertices. Easily, $v(n)\geq 2(\tau(n-1)-1),$ where $\tau$ is the number of positive divisors.

It is easy to see that $\lim \sup v(n)= \infty.$ Can we show that $\lim_{n\rightarrow \infty}v(n) = \infty?$

Also, consider $D(n) = v(n) - 2(\tau(n-1)-1).$ We know that $D(n)>0$ for a set of density 1 in the naturals and furthermore $D(n)=0$ on a set which is $\gg \frac{x}{\log x}$. Can we improve the second estimate?

\item \emph{From Steven Senger, senger@math.udel.edu:}
We will call a family of sets, $P_n\subset [0,1]^2$, $s$-adaptable if they satisfy the following bound:
$$\frac{1}{{n \choose 2}} \sum_{x\neq y; x,y\in P_n} |x-y|^{-s} \lesssim 1.$$
The Szemer\' edi-Trotter incidence theorem says that for a set of $n$ points and $m$ ``reasonable" curves in the plane, the number of incidences of points and curves is bounded above by
$$I \lesssim (nm)^\frac{2}{3}+n+m.$$
Can we get better incidence bounds for $s$-adaptable sets? Specifically, can we get tighter bounds in the case of $n$ points and $n$ circles centered at those points?

\item \emph{From Nathan Pflueger, pflueger@math.harvard.edu:}
Suppose $S$ is a numerical semigroup, $S\subset \mathbb{N}_+$, closed under addition, i.e., $S+S \subset S.$ Let $G:=\mathbb{N}_+\setminus S.$ Define the weight of $S$ to be $w(s) = |\lbrace (x,y)\in S\times G: 0<x<y \rbrace |.$ Define the irreducible elements of $S$ to be the minimal generators. Define the effective weight of $S$ to be $w_{eff}(s) = |\lbrace (x,y)\in S_{irred} \times G: 0<x<y \rbrace |.$ Let the genus of $S$ be $g=|G|$.

For example, $S=\langle 3,5 \rangle.$ Then $w(s) =4$, and $w_{eff}(s) = 3.$

Can we characterize the genus $g$ subgroups of largest effective weight? We believe the largest is $\approx \frac{g^2}{4},$ and in the form $\langle a, a+1, \dots, b-1,b \rangle$, where $b<2a.$

This comes from algebraic geometry. Pick a point $p$ on an algebraic curve or surface. $S=\lbrace ord_p(f): f $ is a rational function$\rbrace,$ where $ord_p(f)$ is the order of the single pole at $p$ of $f$.

\end{itemize}


\subsection{Speaker List}

\begin{itemize}

\item John Bryk, John Jay College (CUNY)
\item Mei-Chu Chang, University of California-Riverside
\item Emel Demirel, Bergen County College
\item Frederic Gilbert, Ecole Polytechnique, Paris
\item Christopher Hanusa, Queens College (CUNY)
\item Charles Helou, Penn State Brandywine
\item Jerry Hu, University of Houston - Victoria
\item Alex Iosevich, University of Rochester
\item Geoff Iyer, University of Michigan
\item Renling Jin, College of Charleston
\item Nathan Kaplan, Harvard University
\item Mizan R. Khan, Eastern Connecticut State University
\item Sandra Kingan, Brooklyn College (CUNY)
\item Alex Kontorovich, Yale University
\item Urban Larsson, Chalmers University of Technology and University of Gothenburg
\item Oleg Lazarev, Princeton University
\item Xian-Jin Li, Brigham Young University
\item Neil Lyall, University of Georgia
\item Steven J. Miller, Williams College
\item Rishi Nath, York College (CUNY)
\item Mel Nathanson, Lehman College (CUNY)
\item Kevin O'Bryant, College of Staten Island (CUNY)
\item Kerry Ojakian, St. Joseph's College, New York
\item Ryan Ronan, Cooper Union
\item Steven Senger, University of Delaware
\item Jonathan Sondow, New York
\item Liyang Zhang, Williams College
\item Wei Zhang, Columbia University
\end{itemize}

\newpage
\section{CANT Problem Sessions: 2013}

\subsection{Problem Session I: Tuesday, May 21st (Chair Steven Miller)}

\subsubsection{MSTD sets and their Generalizations}

Proposed by Steven J. Miller and expanded on by the audience: There are many problems one can ask about More Sums Than Differences sets. Here are just a few.

\begin{itemize}

\item We know that, in the uniform model, a positive percentage of the $2^n$ subsets of $\{0, 1, \dots, n-1\}$ are sum-dominant. Unfortunately these proofs are non-constructive, in that one shows with high probability almost anything thrown between two specially chosen fringes work. Early constructions of explicit families often involved tweaking arithmetic progressions (which are balanced). While these early families were often sub-exponential in terms of their relative size, work of Miller, Scheinerman and Orosz proved that one can find `explicit' families with density $1/n^2$; Zhao obtained a density of $1/n$ through the use of bidirectional ballot sequences. Can one find an explicit formula with a better density (or, dare to dream, one that is a positive percentage?).\\ \ \\

\item Continue to investigate phase transitions, and the natural of the relative size function, for more summands with different combinations of size. \textcolor{red}{This is currently being studied by students in Miller's 2013 REU at Williams.} \\ \ \\

\item Instead of looking at $A+A$ and $A-A$, choose $A$ and $B$ randomly and study $A+B$ and $A-B$ (of course, $A-B$ allow both $a-b$ and $b-a$ for $a \in A$ and $b\in B$.\\ \ \\

\item In determining if $A$ is sum-dominant or difference-dominant, it doesn't matter how much larger one is than the other. Try and find a natural weighting on the sets, try to take into account by how much one beats the other. \\ \ \\

\item Is there a set $A$ such that $|A+A| > |A-A|$ \emph{and} $|A \cdot A| > |A/A|$? If yes, can you find an explicit, infinite family? What is the density of such sets? \textcolor{blue}{Note: Miller finds this problem interesting and wants to bring this to his REU students. Anyone interested in collaborating please email sjm1@williams.edu}.\\ \ \\

\item Instead of looking at subsets of the integers or finite groups, look at subsets of $\Z^d$, intersected with different regions (say spheres, boxes). These sets have different fringe structures. How does the shape of the fringe affect the answer? We can play with the relative sizes of the length and width of a box in two dimensions, for example. \textcolor{red}{This is currently being studied by students in Miller's 2013 REU at Williams.}\\ \ \\

\item Can we say anything about MSTD sets in the continuous case? Is this related to some results on measures? What about subsets of fractals or other special objects (similar to the modular hyperbolas Amanda mentioned). \\ \ \\

\item (Mizan Khan): Speaking of Amanda's talk, the 84\% lower bound mentioned is almost surely not the true answer. What do numerical investigations suggest? What is the correct limiting behavior?

\end{itemize}

\subsubsection{Weakened Convex Functions}

Problem proposed by Seva Lev.\\ \ \\

Consider functions $f: [0,1] \to \R$ that satisfy (1) $\max\{f(0), f(1)\} \le 0$ and (2) for any $0 \le x_1 \le \cdots \le x_m \le 1$ we have $$f\left(\frac{x_1 + \cdots + x_m}{m}\right) \ \le \ \frac{f(x_1) + \cdots + f(x_m)}{m} + (x_m-x_1).$$ Note that convex functions satisfy this.

Set $$F_m(x) \ = \ \sup\{f(x): f \in \mathcal{F}_m\},$$ where $F_m \in C([0,1])$ and $F(0) = F(1) = 0 < F(x)$ for $0 < x < 1$. What is $F_m$ explicitly?

Results are known for $m \in \{2, 3, 4\}$. For such $m$ we have $$F_m(x) \ = \ \sum_{k=1}^\infty m^{-k} \min\{||m^{k-1} x||, 1/m\}.$$ When $m=2$ we have $2 \omega(x)$, where $$\omega(x) \ = \ \sum_{k=0}^\infty 2^{-k} ||2^k x||.$$ What about $m=4$?


\subsection{Problem Session II: Wednesday, May 23rd (Chair Seva Lev)}

\subsubsection{Matrices and Curves}

Problem proposed by Seva Lev. \\ \ \\ Consider an $m \times n$ matrix $A$ whose entries are 0 or 1. Consider $n$ points in the plane $\{p_1, \dots, p_n\}$, with each point corresponding to a column of $A$. If there exists $m$ curves (continuous, no self-intersection) $\{c_1, \dots, c_m\}$ with each curve corresponding to a column of $A$, such that

\bi

\item curve $c_i$ passes through point $p_j$ if $A_{i,j}=1$ and does not pass through point $p_j$ if $A_{i,j}=0$, and

\item any two curves intersect at most once,

\ei

we will call $A$ realizable by curves. \\

The following are questions we can ask: \

\bi

\item What conditions can we put on $A$ to guarantee $A$ is realizable? Note: requiring the dot product of any two rows of $A$ to be at most one does not guarantee $A$ is realizable.

\item Can you find a small matrix that is not realizable?

\item Lastly, if $A$ is realizable, does this mean $A^T$ is realizable? The speaker does not see a reason this should be true, but hasn't found a counterexample yet.

\ei

This problem might be related to planar graphs.

\subsubsection{Enumerating Points in the Plane with Polynomials}

Problem proposed by Mel Nathanson.  \\ \ \\ Consider the set $P \ = \ \{(x,y): x \geq 0, 0 \leq y \leq \alpha x\}$. Does there exist a bijective polynomial $f: P \rightarrow \mathbb{N} \cup \{ 0 \}$?

For instance, if $\alpha \ = \ 1$, then $f(x,y)=\frac{x(x+1)}{2}+y$. Notice when $y=0$, $f(x,0)$ is a triangular number. However, even when $\alpha=2/3$, it is not clear what $f$ should be or if $f$ even exists.

\subsubsection{Sumsets}

Problem posed by Dmitry Zhelezov. \\ \ \\ Let $B$ be a set such that $|B| = n$. Let $$ B \ + \ B\ \supseteq\ A \ = \ \{a_0 < \dots < a_n\},$$ where $A$ is concave ($ a_1 - a_0 > a_2 - a_1 > \dots > a_n - a_{n-1}$) or convex ($ a_1 - a_0 < a_2 - a_1 < \dots < a_n - a_{n-1}$). Is it true that $|A| = O(n^2)$? \\

Problem posed by Steven Senger. \\ \ \\ Let $A$ and $B$ be ``large'' subsets of $\mathbb{N}$ (or $\mathbb{R},$ or $ \mathbb{Z}, \dots$). Do there exists ``large'' subsets of $\mathbb{N}$, $C$ and $D$, such that $|((A \cdot B) +1) \cap (C \cdot D)| \ = \ |A \cdot B|^{1-\epsilon}$?



\subsection{Problem Session III: Thursday, May 24th (Chair Kevin O'Bryant)}

\subsubsection{3-term Geometric Progressions in Sets of Positive Density}

Problem proposed by Kevin O'Bryant.\\ \ \\

As motivation for this problem, recall Van der Waerden's Theorem: given any partition of $\N$, at least one part has an arithmetic progression of arbitrarily large length.  Similarly, we have Szemeredi's Theorem: given any set of positive density in $\N$, there exists an arithmetic progression of arbitrarily large length.  Here we are defining the density of $A \subseteq \N$ as $d(A) = \lim_{n \to \infty}  \frac{|A \cap [1,n)|}{n}$.

It is known that Van der Waerden's Theorem holds for geometric progressions as well.  We would like to consider Szemeredi's Theorem for geometric progressions, but unfortunately it is not true: the square-free integers provide a simple counter-example.  Currently, there is work being done on which densities we can obtain with no geometric progressions.

The original problem proposed in the session was: if $A \subseteq \N$ has density $1$, does $A$ have a three-term geometric progression?  After some Googling by Nathan Kaplan, a 1996 paper by Brienne Brown and Daniel M. Gordon, ``On Sequences Without Geometric Progressions'', was found which stated that if $A \subseteq \mathbb{N}$ has a density and has no 3-term geometric progressions, then the density of $A$ is bounded by .869.

The revised problem proposed is: Given a subset $A \subseteq \N$, which densities of $A$ guarantee 3-term geometric progressions?

\subsubsection{Convex Subsets of Sumsets}

Problem proposed by Dmitry Zhelezov and requested by Giorgis Petridis.\\ \ \\

We consider a variant of the Erd\"os-Newman conjecture, but replace the idea of squaring a set with sumsets.

The problem proposed is: does there exist any set $B$ with $|B| = n$ such that $B + B \supseteq A$ for some convex set $A$ with $|A| =\Omega(n^2)$?


\subsection{Problem Session IV: Friday, May 24th (Chair Renling Jin)}

\subsubsection{Sumsets} Problem proposed by Renling Jin. \\ \ \\

Let $A, B \subseteq \N$ such that $\mathrm{max} A \geq \mathrm{max} B$, $0 \ = \ \mathrm{min} A \ = \ \mathrm{min} B$, and $\mathrm{gcd}(A,B)=1$. Let $\delta \ = \ 1$ if $\mathrm{max} A \ = \ \mathrm{max} B$ and 0 otherwise. If $|A+B| \ = \ |A| + |B| - 2 \delta$, what structure can $A+B$ have? We can also ask the same question if $\delta$ is replaced by $\delta (A,B)$, where $(A,B) \ = \ 1$ if $A \subseteq B$ and 0 otherwise.


\subsection{Additional Problems}

\noindent Proposed by Vsevolod F. Lev.\\ \

For integer $m\ge 2$, let ${\mathcal F}_m$ denote the class of all real-valued
functions $f$, defined on the interval $[0,1]$ and satisfying the boundary
condition $\max\{f(0),f(1)\}\le 0$ and the ``relaxed convexity'' condition
 $$ f\left(\frac{x_1+\ldots+x_m}m\right)
                          \le \frac{f(x_1)+\dotsb+f(x_m)}m + (x_m-x_1), $$
$0\le x_1\le\dotsb\le x_m\le 1$. Now, let
$F_m:=\sup\{f{\rm col}on f\in{\mathcal F}_m\}$. It is easy to prove that
$F_m\in C[0,1]$, $0=F_m(0)=F_m(1)<F_m(x)$ for all $x\in(0,1)$,
$F_m(1-x)=F_m(x)$ for all $x\in[0,1]$, and, somewhat surprisingly,
$F_m\in{\mathcal F}_m$ (meaning that $F_m$ is the maximal function of the
class ${\mathcal F}_m$). What is $F_m$, explicitly? We have
 $$ F_m(x) = \sum_{k=0}^\infty m^{1-k} \min \{ \|m^kx\| ,1/m \},
                                                      \ m\in\{2,3,4\} $$
(where $\|x\|$ denotes the distance from $x$ to the nearest integer),
but for $m\ge 5$ this fails to hold.

\ \\ \ \\

It is easy to see that for any 0-1 matrix, say $M$, one always can find a
system of simple curves and a system of points in the plane so that their
incidence matrix is exactly the matrix $M$. Suppose now that any pair of
curves is allowed to intersect in at most one point (belonging or not to our
system of points), and let's say that $M$ is \emph{realizable} if such curves
and points can be found. Clearly, a necessary condition for this is that the
scalar product of any two rows of $M$ be at most $1$, but this condition is
insufficient: say, for $q$ large enough, by the Trotter-Szemeredi theorem,
the point-line incidence matrix of $PG(q,2)$ has two many incidences to be
realizable. What are other reasonable necessary / sufficient conditions for
$M$ to be realizable? What are "small" examples of non-realizable 0-1
matrices?


\subsection{Speaker List}

\begin{itemize}

\item Paul Baginski, Smith College

\item Arnab Bhattacharyya, DIMACS, Rutgers University

\item Gautami Bhowmik, Universite de Lille, France

\item Thomas Bloom, University of Bristol, UK

\item Tomas Boothby, Simon Fraser University, Canada

\item Amanda Bower, University of Michigan-Dearborn

\item Jeff Breeding II, Fordham University

\item Javier Cilleruelo, University of Madrid, Spain

\item David Covert, University of Missouri - St. Louis

\item Matthew Devos, Simon Fraser University, Canada

\item Mauro Di Nasso, University of Pisa, Italy

\item Mohamed El Bachraoui, United Arab Emirates University, UAE

\item Leopold Flatto, City College (CUNY) and Bell Labs

\item George Grossman, Central Michigan University

\item Christopher R. H. Hanusa, Queens College (CUNY)

\item Derrick Hart, Kansas State University

\item Kevin Henriot, Universite de Montreal

\item Ginny Hogan, Stanford University

\item Jerry Hu, University of Houston - Victoria

\item Renling Jin, Colege of Charleston

\item Delaram Kahrobaei, New York City Tech (CUNY)

\item Nathan Kaplan, Harvard University

\item Omar Kihel, Brock University, Canada

\item Sandra Kingan, Brooklyn College (CUNY)

\item Thai Hoang Le, University of Texas

\item Seva Lev, University of Haifa, Israel

\item Neil Lyall, University of Georgia

\item Richard Magner, Eastern Connecticut State University

\item Steven J. Miller, Williams College

\item Rishi Nath, York College (CUNY)

\item Mel Nathanson, Lehman College (CUNY)

\item Kevin O'Bryant, College of Staten Island (CUNY)

\item Brooke Orosz, Essex County College

\item Giorgis Petridis, University of Rochester

\item Alex Rice, Bucknell University

\item Tom Sanders, Oxford University, UK

\item Steven Senger, University of Delaware

\item Satyanand Singh, New York Tech (CUNY)

\item Jonathan Sondow, New York

\item Dmitry Zhelezov, Chalmers Institute of Technology, Sweden

\end{itemize}

\newpage
\section{CANT Problem Sessions: 2014}

\subsection{Problem Session I: Wednesday, May 28th (Chair Steven J Miller)}

\subsubsection{Steve Senger}

From last year from a talk of Dmitry Zhelezov.

Let $A \subset \mathbb{R}$, $|A| < \infty$, let $P$ be the longest arithmetic progression in $A A = \{ab: a, b \in A\}$. We have $|P| \le c n^{1+\epsilon} \ll n^2$.

Dmintry (possibly from Hegarty): What if instead of $\mathbb{R}$ we have  $\mathbb{F}_q$, the finite field with $q$ elements? Due to Grosu we get the longest progression is at most $c n^{1+\epsilon}$ if $n \le c \log\log\log p$ where $q=p$ a prime. Question: What bounds can I get on the size of $P$ if we replace $\mathbb{R}$  with $\mathbb{F}_q$? Here $q$ can be anything.

\subsubsection{Steven J Miller}

The following builds on my talk from earlier.

\begin{itemize}

\item How does the structure of number affect the answer or the rate of convergence?

\item How does the answer depend on $c$?

\item What is the best way to compute all the $k$-symmetric means for a given $n$? What if we want just a certain one (such as $k = n/2$?).

\item Find other sequences and compute these means -- is there an interesting phase transition?

\end{itemize}

\subsubsection{Steven J Miller}

Consider the 196 game (or problem). Take an $n$-digit number; if it is not a palindrome reverse the digits and add. If the sum is not a palindrome continue, else stop. Lather, rinse and repeat. It's called the 196 problem as 196 is the first number where we don't know if the process terminates (in a palindrome) or goes off to infinity. We know numbers that do not terminate in base 2 (as well as powers of 2, base 11, base 17 and base 26).

What can you say about this problem? What about other bases than 10? What about other decomposition schemes? See  \textcolor{blue}{\url{http://www.math.niu.edu/~rusin/known-math/96/palindrome}.}

\subsubsection{Nathan Kaplan}

Let $C$ be a cubic curve in $\mathbb{P}^2(\mathbb{F}_q)$. Want a large subset so that there are no three points on a line.

The set of $\mathbb{F}_q$-points form an abelian group $G$. Three points sum to zero if and only if they lie on a line.

Given group $G$ what's the largest subset $H$ s.t. $x+y+z = 0$ with $x,y,z$ distinct has no solutions in $H$?

For example takee $G = \mathbb{Z}/2\mathbb{Z} \times G'$ and take $(1,g)$.

What if $G = \mathbb{Z}/p\mathbb{Z}$, consider $\{0, 1, 2, \dots, \lfloor p/3\rfloor\}$.

What if $\mathbb{Z}/5\mathbb{Z} \times G'$, take things of the form $(1,g)$ and $(4,g)$. For each group ask such a question.

One thing you can do is look at a greedy construction. What is the best percentage you can get?



\subsection{Problem Session II: Thursday, May 29th (Chair Kevin O'Bryant)}

\subsubsection{Kevin O'Bryant ?}

A $k$-GP cover of $[N]$ is a family of ${\mathcal F}$ of $k$-term geometric progressions with
\[[N] \subseteq \bigcup_{F \in {\mathcal F}} F.\]

Set
\[\gamma_k\ =\ \lim_{N\to\infty} \inf_{{\mathcal F}} \frac{|{\mathcal F}|}{N}, \]
the infimum being over all $k$-GP covers of $[N]$.

It is easy to see that $\gamma_3\geq \gamma_4 \geq \cdots$, and by basic counting $\gamma_k\geq 1/k$. The cover
\[
{\mathcal F}\ =\ \left\{ b \cdot 2^{ki}\cdot \{1,2,\dots,2^{k-1}\}: 1\leq b \cdot 2^{ki} \leq N, i\geq 0, b\text{ odd}\right\}
\]
shows that
\[ \gamma_k\ \leq\ \frac{2^k}{2(2^k-1)},\]
so that $\lim_{k\to\infty} \gamma_k$ exists and is in $[0,1/2]$. I conjecture that the limit is positive.

Comment from Bloom (in audience): The GP $\{1,2,\dots,2^{N-1}\}$ is covered by the APs $\{1,2,3,4\}, \{8,16,24,32\},
\dots$, so the analog of $\gamma_k$ does go to 0.

Comment from Xiaoyu (in audience): Each GP has at most 2 squarefree numbers, so $\gamma_k \geq 3/\pi^2$. This proves the
conjecture, but leaves open the precise limit.

\subsubsection{David Newman}

About partitions. Finding two sets of partitions which are equal. Finding a type of partition which can divide into two different sets.
$
\prod(1 + x^n)
$
counts partition into distinct parts. Change plus sign into minus signs:
\[
\prod (1-x^n)
=
1 - x - x^2 + x^5 + x^7 + \cdots.
\]

Most of the time, the partitions (?) are equal, but for 1, 2, 5, 7, ... , they differs by one
\[
\prod \frac{1}{1-x^n}
=
(1+x+x^2 + \cdots) (1 + x^2 + x^4 + \cdots) \cdots.
\]

Question: if you change some of the plus signs on the RHS to minus signs, is is possible to get all the coefficients (when you expand) to be ${-1, 0, 1}$.

``I have an example where it works up to $x^{101}$.''

Kevin O'Bryant: throw in lots of $\epsilon$ (taking value in $\pm 1$), becomes a SAT problem.



\subsection{Problem Session III: Friday, May 24th (Chair Mel Nathanson)}

\subsubsection{Mel Nathanson}

Think of $\frac{a}{b}$ as being a parent of 2 children. Left child is $\frac{a}{a+b}$. Right child is $\frac{a+b}{b}$.
Starting with $1$ as root, this gives tree with rows
\begin{gather}
1
\\
1/2, 2
\\
1/3, 3/2, 2/3, 3/1
\\
1/4, 4/3, 3/5, 5/2, 2/5, 5/3, 3/4, 4
\\
\vdots
\end{gather}

``Calkin-Wilf tree''

Every fraction occurs exactly once in this tree.

Start with $z$ (variable) at root instead of $1$. Apply Calkin-Wilf: $z \mapsto (\frac{z}{z+1}, z+1)$. Get linear fractional transformations. $f(z) = \frac{az+b}{cz+d}$.

Rule from parent to children: apply matrices $\begin{pmatrix} 1 & 1 \\ 0 & 1\end{pmatrix}$ and $\begin{pmatrix} 1 & 0 \\ 1 & 1\end{pmatrix}$

Depth formula (involving continued fraction) holds in the Calkin-Wilf tree for $z$.
\[
f(z) = [q_0, q_1, \ldots, q_{k-1}, q_{k} + z]
\text{ if $k$ is even}
\]

\[
f(z) = [q_0, q_1, \ldots, q_{k-1}, q_{k}, z]
\text{ if $k$ is odd}
\]

The form is different for $k$ even/odd. However, if you use that formula above for $k$ even when $k$ is odd, you get a fractional linear transformation with det $-1$. (i.e., like starting a tree with $1/z$, gives tree with det $-1$.) For a given determinant, only finitely many orphans (i.e., no parent) of that determinant.

Question from Nathanson: For a given determinant, how many orphans?

Question from Harald Helfgott: What if the parent to children rules use the matrices $\begin{pmatrix} 1 & 2 \\ 0 & 1\end{pmatrix}$ and $\begin{pmatrix} 1 & 0 \\ 2 & 1\end{pmatrix}$ instead?

Question from Thao Do: what happens if you start with $i$ instead of $1$ as the root? (i.e., let $z = i$) you get tree with elements of $\mathbb{Q}[i]$.

Nathanson: if $z = -1$, then $z \mapsto \frac{z}{z+1}$ gives you $-1$ again... if you're looking at complex numbers.

Thomas Bloom: use fields of characteristic $p$? (Nathanson: ``I don't know anything about char $p$.'')

\subsubsection{Thomas Bloom}

You have the sum set, different set, product set, ratio set.

Because of commutativity, you expect difference set to be larger than sum set.

Question: Is there some subset $A \subset \mathbb{N}$ that is both MSTD and MPTR?

MSTD ``more sums than differences'':  $|A + A| > |A - A|$

MPTR ``more products than ratios'': as well as $|A \cdot A| > |A / A|$.

Nathanson: what is the probability measure? We've seen MSTD sets before but not the multiplicative.

Note: from an MSTD set, can exponentiate to get an MPTR set.

Comments from Thao Do: in order to have MSTD, must have ``almost symmetric form'' i.e., smallest + largest = 2nd smallest + 2nd largest = 3rd smallest + 3rd largest, etc. (then change around a little and be clever)
$a_1 + b_1 = a_2 + b_2$ implies $a_1 - b_2 = a_2 - b_1$.

$I = \{ (a_1, b_1), (a_2, b_2) : a_1 + b_1 = a_2 + b_2\}$.

$J = \{ (a_1, b_10, (a_2, b_2) : a_1 - b_1 = a_2 - b_2\}$.

Nathanson: tell Miller to have students working on MSTD/MPTR over summer. \textbf{Note from Miller: done!}

Bloom: I don't think these sets exist.



\subsection{Speaker List}

Talks here: \textcolor{blue}{\url{http://www.theoryofnumbers.com/CANT2014-program.pdf}.}

\begin{itemize}

\item Sukumar Das Adhikari, Harish-Chandra Research Institute, India
\item Paul Baginski, Fair?eld University
\item Thomas Bloom, University of Bristol
\item Bren Cavallo, CUNY Graduate Center
\item Alan Chang, Princeton University
\item Jean-Marc Deshouillers, IPB-IMB Bordeaux, France
\item Charles Helou, Penn State Brandywine
\item Nathan Kaplan, Yale University
\item Sandra Kingan, Brooklyn College (CUNY)
\item Angel Kumchev, Towson State University
\item Thai Hoang Le, University of Texas
\item Eshita Mazumdar, Harish-Chandra Research Institute, Allahabad, India
\item Nathan McNew, Dartmouth College
\item Steven J. Miller, Williams College
\item Mel Nathanson, Lehman College, CUNY
\item Lan Nguyen, University of Wisconsin-Parkside
\item Kevin O'Bryant, College of Staten Island, CUNY
\item Alberto Perelli, University of Genova, Italy
\item Giorgis Petridis, University of Rochester
\item Luciane Quoos, Instituto de Matem\'atica, UFRJ, Rio de Janeiro, Brasil
\item Steven Senger, University of Delaware
\item Satyanand Singh, New York City Tech (CUNY)
\item Jonathan Sondow, New York
\item Yonutz V. Stancescu, Afeka College, Tel Aviv, Israel
\item Tim Susse, CUNY Graduate Center
\item Johann Thiel, New York City Tech (CUNY)

\end{itemize}

The abstracts are here: \textcolor{blue}{\url{http://www.theoryofnumbers.com/CANT2014-abstracts.pdf}.}

\section{CANT Problem Sessions: 2015}

\subsection{Phase Transitions in MSTD sets: Steven J Miller}

In previous years I talked about phase transitions in the behavior of $|A+A|$ and $|A-A|$ when each element in $\{0, \dots, N\}$ is chosen independently with probability $p(N) = N^{-\delta}$ as $\delta$ hits 1/2. What happens with three summands? Four?

What happens if we restrict $A$ to special types of sets? How does the additional structure affect the answer?

Also, is there an `explicit' construction of an infinite family of MSTD sets? The word `explicit' is deliberately not being defined; I would like some nice, concrete procedure that does not involve randomness.

Finally, is this the 14\textsuperscript{th} or the 13\textsuperscript{th} CANT?

\subsection{An accidental sequence: Satyanand Singh}

The two outer graphs which form an envelope around $\gamma_{3}((6j+2)^5)$ illustrate that:
\[
{\left(\frac{\ln(6j+2)^5}{3\ln3}\right)} <\gamma_{3}((6j+2)^5)<{\left(\frac{\ln(6j+2)^5}{\ln3}\right)}.
\]

The upper bound is easily seen by finding the power of $3$ that is closest to $(6j+2)^5 $ but does not exceed it. We can also say for certain that $\gamma_{3}((6j+2)^5)\geq3,$ since Bennett dispensed of the two term case in \cite{MB} and equality occurs when $2^5=3^3+3^1+2.$ We were not able to prove the lower bound suggested by the experimental results, i.e., $ \gamma_{3}((6j+2)^5)>{\left(\frac{\ln(6j+2)^5}{3\ln3}\right)}$ for $j\geq{1}$. This would completely resolve the case for $q=5.$

\ \\

\noindent {\bf Problem 1.}
{\it For both $a$ and $b$ odd, where $a>b>0$, find all solutions to the diophantine equation $3^{a}+3^{b}+2=(6j+2)^5$ or show that the only solution is $(a,b,j)=(3,1,0).$}\\
\[
\]
\noindent {\bf Problem 2.}
{\it {For any positive integer $n$, with $(n,3)=1,$ find all solutions to $\gamma_{3}(n^{q})\leq3$ for $q$ a prime number where $q >1000$}}? \\
\[
\]
\noindent {\bf Problem 3.}
{\it {For any positive integer $n$, with $(n,3)=1$, we conjecture that $\gamma_{3}((6j+2)^5)>c\ln{(6j+2)^5}$ where $c$ is a constant such that $0<c<1/(3\ln{(3)}).$}}\\

\subsection{Kevin O'Bryant}


Let $b_1,b_2, \dots$ be an infinite binary sequence, and let $A$ be the set of real numbers of the form
$a_i:=\sum_{n=1}^\infty b_{n+i} \cdot 2^{-n}$. If $A$ has no infinite decreasing subsequence, that is, if $A$ is an
ordinal, what are the possibilities for the order type of $A$? In particular, can the order type be $\epsilon_0$?

Blair, Hamkins, and O'Bryant [forthcoming] have shown that the order type, if infinite, must be at least $\omega^2$,
and can be as large as $\omega \uparrow\uparrow n$ for any $n$.

\subsection{Steven Senger (repeat from previous years):} Given a large finite subset, $A$, of real numbers, and any non-degenerate, generalized geometric progression, $G$, with $|G|\approx |AA|$, can we get a nontrivial bound on $|(AA+1)|\cap G|$?

\subsection{Nathan Kaplan:} We say a set of points, $P\subset \mathbb \mathbb{R}^2$, is in general position if no curve has more points than it ``ought to''. That is there are no three points on a line, no six points on a conic, etc.... The original problem posed by Jordan Ellenberg is ``How does the minimum height of a set of completely generic points grow with the number of points?'' It is available at:
\begin{center}
\bburl{https://quomodocumque.wordpress.com/2014/04/05/puzzle-low-height-points-in-general-position/}
\end{center}

Can anything be said, even if we choose a greedy construction for our points.

\subsection{Nathan Kaplan (repeat from previous years):} In $\mathbb{F}_3^n$, define the function $f(n,m)$ to be the maximum number of lines completely contained in any set of $m$ points. Is it true that the simplest greedy construction (filling in lower dimensional subspaces), is the best possible? That is, is $f(n,m)=f(\lceil \log_3(m)\rceil,m)$?

\subsection{Kevin O'Bryant:} Given a large finite subset, $A$, of real numbers, is it true that $|AA+A|\geq|A+A|$? \\


Oliver Roche-Newton has asked on Math Overflow
\begin{center}
(\bburl{http://mathoverflow.net/questions/204020/is-the-set-aaa-always-at-least-as-large-as-aa/})
\end{center}  if it is possible for
$$| A\cdot A + A| < |A+A|$$ with $A$ being a set of real numbers. Some observations.

\begin{itemize}

\item For a random set $A$ of $k$ real numbers, $A\cdot A + A$ has $\sim k^3/2$ elements while $A+A$ has only $\sim
k^2/2$, so any example needs to have some special structure.

\item Modulo 13, the set $A=\{2,5,6,7,8,11\}$ (the set of positive quadratic non-residues, with a modulus $p\equiv 1
\bmod 4$), is an example.

\item if $| \cdot |$ means Lebesgue measure, then $A=[0,1/2]$ is an example, as $A\cdot A + A = [0,3/4]$ but
$A+A=[0,1]$.

\item if $A$ is a set of 3 or more positive integers, then it cannot be an example, as $$a_1+a_n A, a_2 + a_n
A, \dots, a_n+a_nA$$, where $a_n=\max A$, are necessarily disjoint (reduce modulo $a_n$) and contain at least $|A|^2$
elements altogether, while $|A+A| \leq \binom{|A|+1}{2}$.
\end{itemize}

The audience asked what was known for Hausdorff measure, and suggested considering the problem over the integers,
positive rationals, and complexes.

\newpage

\subsection{Speaker List}\label{sec:participantlist}

\begin{itemize}

\item Paul Baginski, Fairfield University
\item Bela Bajnok, Gettysburg College
\item Dakota Blair, CUNY Graduate Center
\item Lisa Bromberg, CUNY Graduate Center
\item Mei-Chu Chang, University of California
\item Scott Chapman, Sam Houston State University
\item David John Covert, University of Missouri St. Louis
\item Robert Donley, Queensborough Community College (CUNY)
\item Leonid Gurvits, City College (CUNY)
\item Sandie Han, New York City Tech (CUNY)
\item Charles Helou, Penn State Brandywine
\item Alex Iosevich, University of Rochester
\item Renling Jin, College of Charleston
\item Nathan Kaplan, Yale University
\item Mizan Khan, Eastern Connecticut State University
\item Sandra Kingan, Brooklyn College (CUNY)
\item Diego Marques, University of Brasilia
\item Ariane Masuda, New York City Tech (CUNY)
\item Nathan McNew, Dartmouth College
\item Steven J. Miller, Williams College
\item Mel Nathanson, Lehman College (CUNY)
\item Kevin O'Bryant, College of Staten Island (CUNY)
\item Cormac O'Sullivan, Bronx Community College (CUNY)
\item Jasmine Powell, Northwestern University
\item Alex Rice, University of Rochester
\item Steven Senger, Missouri State University
\item Satyanand Singh, New York City Tech (CUNY)
\item Jonathan Sondow, New York
\item Johann Thiel, New York City Tech (CUNY)
\item Yuri Tschinkel, NYU
\item Bart Van Steirteghem, Medgar Evers College (CUNY)
\item Madeleine Weinstein, Harvey Mudd College

\end{itemize}

A list of talks  and abstracts is available online here: \begin{center} \textcolor{blue}{\url{http://www.theoryofnumbers.com/CANT2015-abstracts.pdf}.}\end{center}

\ \\


\newpage

\section{CANT Problem Sessions: 2016}

\subsection{Problem Session I: Tuesday, May 24th (Nathanson Chair)}

\ \\
\noindent \textbf{Mel Nathanson:} Alexander Borisov in his 2005 arXiv paper, ``Quotient singularities, integer ratios of factorials and the Riemann Hypothesis,'' discussed integer valued ratios of factorials and their relation to problems in number theory and algebrai geometry.  Historically the application of such ratios to number theory goes back at least to Pafnuty Chebyshev, who used them to obtain the order of magnitude of $\pi(x)$.
Since $\binom{2n}{n} = \frac{(2n)!}{n! n!}$ is an integer, it is not hard to see
that its prime decomposition must include all primes $p$ such that $n < p < 2n$, and so
$$4^n > \frac{(2n)!}{n! n!}\geq\prod_{n<p<2n}p$$
where $p$ is a prime number. Chebyshev used this fact to show that
$$\frac{x}{\log(x)}\ll\pi(x)\ll\frac{x}{\log(x)}$$

It is a theorem of Eug\'ene Charles Catalan
that for any $k\in\mathbb{N}$ one has
$$\frac{(2n)!(2k)!}{n! k!(n+k)!}\in\mathbb{Z}$$
There are combinatorial proofs of this identity for $k=0,1,2$.

\noindent \textbf{Question:} For  $3 \leq k \leq n$, prove
$$\frac{(2n)!(2k)!}{n! k!(n+k)!}\in\mathbb{Z}$$
by a  counting argument.\\ \

One can show that both $$\frac{(9n)! n!}{(5n)!(3n)!(2n)!}$$
and $$\frac{(14n)! (3n)!}{(9n)!(7n)!(n)!}$$ are integers for all positive
integers $n$.

\noindent \textbf{Question:}
Find all quintuples $(a,b,c,d,e)\in\mathbb{N}^5$  such that
$$ a + b = c + d + e $$
$$ \gcd(a,b,c,d,e) = 1 $$
and
$$\frac{(an)!(bn)!}{(cn)!(dn)!(en)!}\in\mathbb{Z}\text{?}$$
Show that there are  only 29 such quintuples.  \\

\noindent \textbf{Question:} Can one deduce something  interesting about the
distibution of the primes from
an integer identity of the form $\frac{(an)!(bn)!}{(cn)!(dn)!(en)!}$?\\
\ \\

\noindent \textbf{Thomas Blume:} Does there exist $A\subset\mathbb{Z}$ such that
\begin{enumerate}
    \item $\lvert A+A\rvert>\lvert A-A\rvert$,\\
    \item $\lvert A\times A\rvert>\lvert A/A\rvert$.
\end{enumerate}

\ \\
\noindent \textbf{Mel Nathanson:}
J. A. Haight showed in a paper from 1973 that for all $h$ and $l$, there exists a modulus $m$ and $A$ a subset of $\mathbb{Z}/m\mathbb{Z}$ such that $A-A=\mathbb{Z}/m\mathbb{Z}$ but $hA$ omits $l$ consecutive residues in $\mathbb{Z}/m\mathbb{Z}$.
He then used this algebraic result to show the following.

\begin{thm}
There exists $E\subset\mathbb{R}$ such that $E-E=\mathbb{R}$ but $\mu(hE)=0$ for all $h \geq 1$, where $\mu$ is  Lebesgue measure.
\end{thm}

This is a kind of reverse MSTD result for the reals.  \\

\noindent \textbf{Question:}
Let $\epsilon>0$.  Does there exist $A\subset\mathbb{R}$ such that $A-A=\mathbb{R}$ and $\mu(2*A-A)<\epsilon$.\\

Using Haight's results, Ruzsa was able to show that for any fixed $h$, there exists an $A\subset\mathbb{N}$ such that $\lvert A-A\rvert$, but $\lvert hA\rvert$ is small.\\ \

\noindent \textbf{Question:}
Does there exist $A\subset\mathbb{R}$ such that $A-A=\mathbb{R}$ and $\mu(2*A-A)<\epsilon$, $\forall\epsilon>0$.\\ \

Specifically, Ruzsa showed the following.

\begin{thm}
For any $h>1$, and any $\epsilon>0$, there exists a modulus $m$ and an $A\subset\mathbb{Z}/m\mathbb{Z}$ such that $A-A=\mathbb{Z}/m\mathbb{Z}$, and $\lvert h A\rvert<\epsilon\cdot m$.
\end{thm}
If one defines $\Phi(t_1,t_2,\ldots,t_h)=\sum_{i=1}^h t_i$ and $\Upsilon(t_1,t_2)=t_1-t_2$, then a slight reinterpretation of the above theorem says that $\lvert\Phi(A)\rvert<\epsilon\cdot m$ and $\Upsilon(A)=\mathbb{Z}/m\mathbb{Z}$.\\

\ \\

\noindent \textbf{Question:} What, if anything, is special about these linear forms? If one takes $\Phi(t_1,t_2,\ldots,t_h)=\sum_{i=1}^h \phi_it_i$, for some function $\phi$ on the index set of $\Phi$ and similarly a $\psi$ for $\Upsilon(t_1,t_2,\ldots,t_g)=\sum_{i=1}^g \psi_it_i$, for what functions $\psi$ and $\phi$ does one get a Haight like result?\\ \

Consider $(\phi_1,\phi_2,\ldots,\phi_h)$ and $I=\{1,2,\ldots,h\}$. Then let $S_I^{(\Phi)}=\sum_{i\in I}\phi_1$ and $S_I^{(\Upsilon)}$, then\\

\noindent \textbf{Exercise:} Show that when $\Upsilon(t_1,t_2)=t_1-t_2$ and $\Phi(t_1,t_2)=2t_1-t_2$ one gets a Haight like result, (i.e., $\exists A\subset\mathbb{R}$ such that $A-A=\mathbb{R}$ and $\mu(2*A-A)<\epsilon$.\\

\noindent \textbf{Question:} What if the measure $\mu$ in the above statements is replaced with Hausdorff dimension? Are there Haight like results that one can describe where the difference set is of dimension 1, and the $h$ fold sum set is of fractional dimension?\\

\ \\
\noindent \textbf{William Keith:}
Let $P=\prod_{i=1}^\infty(1+q^{2i})$ and $Q=\prod_{i=1}^\infty(1+q^{2\cdot4^i})$ for $q$ a prime. Note that $(1-q)^2\equiv(1+q^2)$ and $P\cdot Q\equiv1+q+q^2+\ldots$. \ \\

\noindent \textbf{Question:} When is it that $F=\sum f(m)q^m$, where $f(m)$ is odd with positive density, that $(F)^k$ has zero density for the odd coefficients?

\ \\
\noindent \textbf{Larsen Urban:} Let $A\subset\mathbb{N}^d$, if $A+A=A^c\setminus T_A$, where $T_A=\{z\in\mathbb{N}^d\ \mid\ z\text{ is unrelated to anything in }A\}$, then what can one say about $A$? Let $A,B\subset\mathbb{N}^d$, when is it that $\varnothing=A\cap B$, and $A+B=(A\cup B)^c$.

\subsection{Problem Session II: Wednesday, May 25th (Iosevich Chair)}

\ \\
\noindent \textbf{Urban Larsson:} Suppose $X$  is a sum-free set on $\Z_{\ge 0}$; such sets are known, but now require $\min X =  \in \Z_{\ge 0}$. What is the maximum density of such a set $X$ if it is sum-free and has smallest element $k$?  What if we vary $k$? Has this been studied? Is it $$X \ = \ \{i p_k + k, \dots, i p_k + 2k - 1: i \in \Z_{\ge 0}, \ \ p_k \ = \ 3k-1.\}$$

\ \\
\noindent \textbf{Brendan Murphy:} Inspired by Alex and Tom's talks, say you color $\F_q^\ast$, $E \subset \F_q^d$, $d \ge 2$, $x$ and $y$ connected by an edge if $||x-y|| = t \in \F_q^\ast$, with $||x|| = x_1^2 + \cdots + x_d^2$. How large must $E$ be so that we have monochromatic ($||x-y||$, $x \cdot y$). Comment from audience: hope. Alex: for large sets might be ok (maybe $|E| > C q^{(d+1)/2}$), small sets....

\ \\
\noindent \textbf{Sarfraz Ahmad:} Goal is to prove that for all $i > 0$ we have $$\sum_{j=0}^i \frac{(-1)^j (i-j + 1/2)^i}{j! (i-j)!} \ = \ 1.$$ Comment from audience: is this a difference order problem? $$\Delta_h^i f(x) \ = \ \sum_{j=0}^i \ncr{i}{j} (-1)^i f(x+j),$$ where $$\Delta_h^{i_1+i_2} \ = \ \Delta_h^{i_1} \left(\Delta_h^{i_2}\right).$$

\subsection{Problem Session III: Thursday, May 26th (Miller Chair)}

\ \\
\noindent \textbf{Steven Miller:} Prove unconditionally that there are infinitely many subsets of the primes that are MSTD sets. What about other special sets? Answered at lunch and in talk: Use  Green-Tao, follows  immediately.

\ \\
\ \\
\noindent \textbf{Mel Nathanson:} For all $A, B \subseteq \R^n$, we have $(A \cap \Z^n) + (B\cap \Z^n) \subseteq (A+B) \cap \Z^n$. The special case $A=B$ is interesting.
For example,  in $\R^3$, the Reeve polytope $A$ is the convex hull of the set
$(0,0,0), (1,0,0), (0,1,0), (1,1,2)$. The lattice points in $A\subset \Z^3$
are $\{(0,0,0), (1,0,0), (0,1,0), (1,1,2)\}$. The sumset $2A$ contains the lattice  point
$(1,1,1) = (1/2, 1/2, 0) + (1/2, 1/2, 1)$, but it is not the sum of two lattice points in $A$.
Thus,  $(A \subset \Z^3) + (A \subset \Z^3)$ is a proper subset of $2A \cap \Z^3$.
Describe the lattice polytopes $A$ such that $2(A \subset \Z^n) = (2A) \cap\Z^n$.

\begin{figure}
\begin{center}
\scalebox{.5}{\includegraphics{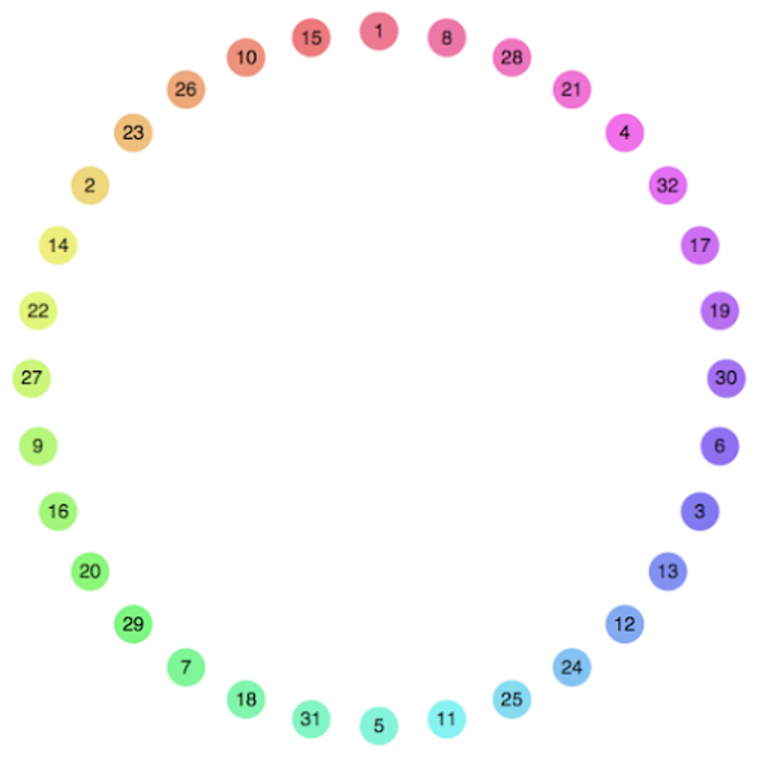}}
\caption{\label{fig:colorwheel} B\"ogeholz's arrangement of numbers such that sums of adjacent elements are squares.}
\end{center}\end{figure}

In the plane there exist lattice triangles  where the sum of the  triangles contains
a lattice point that is not the sum of lattice points in the triangle.
For example, the triangles with vertices at (0,0), (1,0), (0,1) and at (0,1), (1,3) and (2,4) have this property.

\ \\
\noindent \textbf{Kevin O'Bryant:} A couple of months ago on Google+, Harald B\"ogeholz found an arrangement of the integers from 1 to 32 such that any adjacent pair adds
up to a square. See his graphic in Figure \ref{fig:colorwheel}, or \bburl{https://plus.google.com/u/0/106537406819187054768/posts/Y9qaWEwiLuv?cfem=1}.

For $N \le 31$ it is \emph{impossible} to arrange the numbers up to $N$ on a circle so that each adjacent pair sums to a perfect square. For example, if $N=19$ what goes
next to 16? Could have 9 but then in trouble as need two neighbors. Considering the graph with vertices $\{1,2,\dots,N\}$ and edges connecting numbers that sum to a
square, we are asking for a Hamiltonian cycle. Unique ways for $N=32, 33$; number of ways of doing is not monotonic.

Implied question: is it possible for every $N \ge 32$? For infinitely many $N \ge 32$?

Could add all the sums, that gets each number twice, so that would be $N(N+1)$, has to be a linear combination of the squares, gives a Diophantine condition which maybe
could be easily checked. For some $N$, this condition can satisfied even though there is no Hamiltonian cycle.

Generalizations: What about three in a row added? What about a cube? A sphere?

\ \\
\noindent \textbf{Sergei Konyagin:} What is the cardinality of the maximal subset of $\{1, \dots, N\}$ such that $A$ does not contain an MSTD subset? Comment: if have an arithmetic progression of length 15 fail.

\ \\
\noindent \textbf{Colin Defant:} Define  $\sigma_c(n) = \sum_{d|n} d^c$, look at $\sigma_c(\N)$, gives a set of complex numbers, look at closure, for which complex numbers $c$ is $\overline{\sigma_c(\N)}$ connected? If ${\rm Re}(c) < -3.02$ (approximately) not connected, this bound is probably not optimal.

\subsection{Problem Session IV: Friday, May 27th (O'Bryant Chair)}

\ \\
\noindent \textbf{Sergei Konyagin:} Given a natural number, $r\geq 2$, consider the equation over natural numbers, $1\leq x_i,y_i\leq N.$
\begin{equation}\label{recips}
\frac{1}{x_1}+\frac{1}{x_2}+\cdots+\frac{1}{x_r}\ = \ \frac{1}{y_1}+\frac{1}{y_2}+\cdots+\frac{1}{y_r}.
\end{equation}
There are roughly $r!N^r$ trivial solutions of Equation \ref{recips}, where the $x_i$ are permutations of the $y_i$. Define $\mathcal F_{r,n}$ to be the number of nontrivial solutions to Equation \ref{recips}. Konyagin and Korolev have shown that as $N\rightarrow \infty$,
$$\mathcal F_{r,n}\leq O\left(N^{r-\frac{1}{4}}\right).$$
They conjecture that the upper bound should really be $O\left(N^{r-1+o(1)}\right).$
Consider also the following example. For any choice of $z$ with $1\leq z\leq N.$
$$\frac{1}{2z}+\frac{1}{4z}+\frac{1}{x_3}+\cdots+\frac{1}{x_r}\ =
\ \frac{1}{3z}+\frac{1}{3z}+\frac{1}{x_3}+\dots+\frac{1}{y_r}.$$
To bar such examples, we could consider the assumption that $x_i\neq y_i$ for all $i$ in Equation \ref{recips}. For this variant, we could have $x_1=x_2=2z$, $y_1=z$, $x_3=z'$, $y_2=y_3=2z'$, etc\dots The conjecture for the variant is that the number of nontrivial solutions to Equation \ref{recips} with this additional restriction should be no more than
$$O\left(N^{\left\lfloor\frac{2r}{3} \right\rfloor+o(1)} \right).$$
We know that, for some positive constants $c$ and $C$, $\mathcal F_{r,n}\approx r!N^r$ when $r\leq c\left( \frac{\log N}{\log \log N}\right)^\frac{1}{3},$ but that $\mathcal F_{r,n}> r!N^r$ when $r > C\left( \frac{\log N}{\log \log N}\right)^\frac{1}{3}.$ This is known when one considers the choices of $y_j$ to be the smooth numbers.

Question 1 (Shparlinski): Consider the solutions to Equation \ref{recips} with $M+1\leq x_i,y_i\leq M+N$.

Question 2 (Konyagin): Consider the solutions to the following equation
$$\frac{a_1}{x_1}+\frac{a_2}{x_2}+\dots+\frac{a_r}{x_r}=\frac{1}{y_1}+\frac{1}{y_2}+\dots+\frac{1}{y_r},$$
where $a_j$ are nonzero rational numbers, and the $x_j$ and $y_j$ are as before.

Question 3 (Senger): Consider the solutions to Equation \ref{recips}, except with different numbers of terms on each side.\\ \

\ \\
\noindent \textbf{Kamil Bulinski:} Let the group $G=\mathbb Z / N\mathbb Z \oplus \mathbb Z / M\mathbb Z$, and
$$G\ =\ \bigsqcup_{i=1}^m(a_i+H_i),$$
where the $H_i$ are subgroups of $G$. Must $H_i=H_j$ for some $i\neq j$? Note that this is false for the case that $G=\mathbb Z / 2\mathbb Z \oplus \mathbb Z / 2\mathbb Z \oplus \mathbb Z / 2\mathbb Z,$ as $$G\ = \ \{(0,0,0),(1,0,0)\}\cup\{(0,0,1),(0,1,1)\}\cup\{(1,1,0),(1,1,1)\}\cup\{(0,1,0),(1,0,1)\},$$ a union of four disjoint lines.

\ \\
\noindent \textbf{Kevin O'Bryant:} Sun's Conjecture: If $a_1+H_1,a_2+H_2, \dots a_m+H_m$ are disjoint, then there exist $i<j$ such that gcd$([G:H_i],[G:H_j])\geq m$.

\ \\
\noindent \textbf{Brian Hopkins:} Let $p(n,3)$ denote the number of partitions of a natural number, $n$, into exactly three parts. It is known that $p(n,3)$ is the nearest integer to $\frac{n^2}{12}.$ This tells us that for Pythagorean triples, $a,b$, and $c$, where $a^2+b^2=c^2$, we have that
\begin{equation}\label{star}
p(a,3)+p(b,3)\ =\ p(c,3).
\end{equation}

Question 1: Is there a direct bijective proof of Equation \ref{star}?

Question 2: If a triple, $(a,b,c)$ satisfies Equation \ref{star}, it may not be a Pythagorean triple. Characterize the triples for which Equation \ref{star} holds.

\ \\
\noindent \textbf{Steven Senger:} Tom Sanders spoke on colorings of the natural numbers where there exists a quadruple, $(x,y,x+y,xy)$, whose entries are all the same color.

Question 1: Can one show that any four-coloring (with equally dense sets of colors) of the natural numbers will guarantee the existence of a quadruple $(x,y,x+y,xy)$ are all different colors? Note that if each residue class modulo 4 is given a different color, then restricting possible colors to fixed places may render the answer as no, so we must allow any color to be in any entry.

Answer (Ryan Alweiss): NOPE! Actually, even this restriction is irrelevant, as giving each residue class modulo 4 will render the answer negative. To see this, note that if $x$ and $y$ are both even, $x+y$ will also be even, so they cannot come from three distinct equivalence classes modulo 4. To see this, note that if $x$ and $y$ are both odd, then $xy$ will also be odd, so they cannot come from three distinct equivalence classes modulo 4. Therefore, $x$ and $y$ must have different parity. This will also not work, as can be seen by checking each case.

Question 2: Perhaps more colors could work?

\ \\
\noindent \textbf{Tom Bloom:} Let $A\subset \{1, 2, \dots, 2N \}$ is $N$-circular (as discussed by Kevin O'Bryant earlier, there exists a permutation of $\{1, 2, \dots, N\}$ such that pairwise, $a+b \in A$).

Question 1: For which $p(n)$ is it true that a random subset of $[2N]$ is $N$-circular with high probability.

Question 2: If $\{1, 2, \dots, N\}=X\sqcup Y$ such that $(X+Y)\cap A=\varnothing,$ then $A$ is not $N$-circular. Are there any other natural obstructions?

\ \\

\subsection{Speaker List}

\begin{center}
\fbox{\Huge {\rm col}or{red} CANT 2016 Speakers}  \\
 \vspace{0.5cm} {\LARGE  Fourteenth Annual Workshop on \\
Combinatorial and Additive Number Theory}  \\ \vspace{0.3cm}
{\Large CUNY Graduate Center  \\
May 24--27, 2016}
\end{center}

\vspace{1cm}

\begin{enumerate}
\item
{\bf Sarfraz Ahmad}, Comsats Institute of Information Technology,
Lahore, Pakistan

\item
{\bf Paul Baginski}, Fairfield University

\item

{\bf Gautami Bhowmik}, Universit of Lille, France

\item
{\bf Pierre Bienvenu}, University of Bristol, UK
\item

{\bf Arnab Bose}, University of Lethbridge, Canada

\item
{\bf Kamil Bulinski},  University of Sydney, Australia

\item
{\bf Sam Cole}, University of Illinois at Chicago

\item
{\bf Colin Defant}, University of Florida

\item
{\bf Mohamed El Bachraoui}, United Arab Emirates University, United Arab Emirates

\item
{\bf George Grossman},   Central Michigan University

\item
{\bf Sandie Han}, New York City Tech (CUNY)

\item
{\bf Brian Hopkins}, St. Peter's University

\item
{\bf Alex Iosevich}, University of Rochester

\item
{\bf William J. Keith}, Michigan Technological University

\item
{\bf Mizan Khan}, Eastern Connecticut State University

\item
{\bf S. V. Konyagin}, Steklov Mathematical Institute of the  Russian Academy of Sciences,
Russia

\item
{\bf Ben Krause}, University of British Columbia, Canada

\item
{\bf Urban Larsson}, Dalhousie University, Halifax, Canada

\item
{\bf Jiange Li}, University of Delaware

\item
{\bf Ray Li}, Carnegie-Mellon University

\item
{\bf Neil Lyall}, University of Georgia

\item
{\bf Akos Magyar}, University of Georgia

\item
{\bf Ariane Masuda}, New York City Tech (CUNY)

\item
{\bf Nathan McNew}, Towson State University

\item
{\bf Steven J. Miller},  Williams College

\item
{\bf Amanda Montejano},  UMDI-Facultad de Ciencias,
Universidad Nacional Aut\'onoma de M\'exico, Quer\'etaro, M\'exico

\item
{\bf Brendan Murphy}, University of Rochester

\item
{\bf Rishi Nath}, York College (CUNY)

\item
{\bf Mel Nathanson},
Lehman College (CUNY)

\item
{\bf P\' eter P\' al P\' ach},   Budapest University of Technology and Economics, \\
Hungary

\item
{\bf Zhao Pan}, Carnegie-Mellon University

\item
{\bf Giorgis Petridis}, University of Rochester

\item
{\bf Bradley Rodgers}, University of Michigan

\item
{\bf Ryan Ronan}, CUNY Graduate Center

\item
{\bf Tom Sanders}, Oxford University, UK

\item
{\bf James Sellers}, Pennsylvania State University

\item
{\bf Steven Senger}, Missouri State University

\item
{\bf Satyanand Singh}, New York City Tech (CUNY)
\item
{\bf Jonathan Sondow}, New York

\item
{\bf Yoni Stancescu}, Afeka College, Israel
\item
 {\bf Johann Thiel}, New York City Tech (CUNY)

\item
{\bf Andrew Treglown}, University of Birmingham, UK

\item {\bf Yuri Tschinkel}, New York University

\item
{\bf Van Vu}, Yale University

\item
 {\bf Huanzhong Xu}, Carnegie-Mellon  University

\item
{\bf Victor Xu}, Carnegie-Mellon University


\item
{\bf Xiaorong Zhang}, Carnegie-Mellon University

\end{enumerate}

\vspace{2cm}
\subsection{Countries represented}
\begin{enumerate}

\item
Australia

\item
Brazil

\item
Canada

\item
France

\item
India

\item
Israel

\item
Mexico

\item
Pakistan

\item
Russia


\item
Sweden

\item
United Arab Emirates

\item
United Kingdom

\item
United States

\end{enumerate}



\section{CANT Problem Sessions: 2017}

\noindent Fifteenth Annual Workshop on Combinatorial and Additive Number Theory

\subsection{Problem Session I: Wednesday, May 24th (Chair Kevin O'Bryant)}

\subsubsection{Sam Chow: \textcolor{blue}{\href{mailto:sam.chow@york.ac.uk}{sam.chow@york.ac.uk}}}

Problem on behalf of Oleksiy Klurman.

Theorem of Shao (2013): Let $\mathcal{A}$ be a subset of the primes $\mathcal{P}$ with $$\underline{\delta} \ = \ \underline{\lim} \frac{|\mathcal{A} \cap [N]|}{|\mathcal{P} \cap [N]} \ > \ \frac58.$$ Then for $N$ large and odd $$p_1 + p_2 + p_3 \ = \ N$$ has a solution with the $p_i \in \mathcal{A}$.

This is sharp. For example, $\mathcal{A} = \{p \in \mathcal{P}: p \equiv 1, 2, 4, 7, 13 \bmod 15\}$ and $N \equiv 14 \bmod 15$ has no solution.

Question: Let $\mathcal{A} \in \N$, $\underline{\delta} = \underline{\lim}|\mathcal{A} \cap [N]| / N > \delta^\ast$. Then for $N$ large $$x_1^2 + \cdots + x_5^2 \ = \ N$$ has a solution with the $x_i \in \mathcal{A}$. Local problem: $x_1^2 + \cdots + x_5^2 \equiv N \bmod q$.

From the audience: is it easier with more variables? Answer: Claim it doesn't get easier as go from 5 variables to 100.  What is the smallest value of $\delta^\ast$ one can take?

\subsubsection{Colin Defant: \textcolor{blue}{\href{mailto:cdefant@ufl.edu}{cdefant@ufl.edu}}}

Definition: A $k$-antipower is a word $w = w_1 w_2 \cdots w_k$ where $w_1, w_2, \dots, w_k$ are distinct words which all have the same length (so $|w_1| = |w_2| = \cdots = |w_k|$).

Definition: Thue-Morse Word: Start with $T = 0$. The complement of 0 is 1, append to end, have 01. Take complement 10 and append to end, have 0110, keep going, get infinite word 0110100110010110.... Given a positive integer $k$ let $\gamma(k)$ be the smallest odd positive integer $m$ such that the prefix of $T$ of length $km$ is a $k$-antipower. I proved with Shyam Narayanan d that $$\frac34 \ \le \ \liminf_{k\to\infty}  \frac{\gamma(k)}{k} \ \le \ \frac9{10};$$ what is the value of the liminf? Expect it is 9/10. We know the limsup is 3/2.

See for example \bburl{http://www.combinatorics.org/ojs/index.php/eljc/article/view/v24i1p32/pdf} and  \bburl{https://arxiv.org/pdf/1705.06310.pdf} and \bburl{https://arxiv.org/pdf/1606.02868}.

\subsubsection{Salvatore Tringali: \textcolor{blue}{\href{mailto:salvo.tringali@gmail.com}{salvo.tringali@gmail.com}} (preferred) or \textcolor{blue}{\href{mailto:salvatore.tringali@uni-graz.at}{salvatore.tringali@uni-graz.at}}}

Some preliminaries: Let $H$ be a multiplicatively written monoid with identity $1_H$. We denote by $H^\times$ the group of units of $H$, and by $\mathcal A(H)$ the set of all $a \in H \setminus H^\times$ such that there do not exist $x, y \in H \setminus H^\times$ for which $a = xy$. We refer to the elements of $\mathcal A(H)$ as the \textit{atoms} of $H$: If you think of the special case where $H$ is the multiplicative monoid of a unital ring, you will find that atoms are no different from the usual notion of an irreducible element in commutative algebra.

Given $x \in H$, we set
$
\mathsf L_H(x) := \{k \in \mathbf N^+: x = a_1 \cdots a_k \text{ for some }a_1, \ldots, a_k \in \mathcal A(H)\}
$ if $x \ne 1_H$, and $\mathsf L_H(x) := \{0\} \subseteq \mathbf N$ otherwise: We call $\mathsf L_H(x)$ the \textit{set of lengths} of $x$ (relative to the atoms of $H$), while the family $\mathscr L(H) := \{\mathsf L_H(x): x \in H\}$ is termed the \textit{system of sets of lengths} of $H$. It is straightforward that
$$
\bigl\{\{0\}\bigr\} \subseteq \mathscr L(H) \subseteq \bigl\{\{0\}, \{1\}\bigr\} \cup \mathcal P(\mathbf N_{\ge 2}),
$$
where $P(\mathbf N_{\ge 2})$ is the power set of $\mathbf N_{\ge 2}$. Moreover, if $H$ is a reduced BF-monoid, then
$$
\mathscr L(H) \subseteq \bigl\{\{0\}, \{1\}\bigr\} \cup \mathcal P_{\rm fin}^\ast(\mathbf N_{\ge 2}),
$$
where $\mathcal P_{\rm fin}^\ast(\mathbf N_{\ge 2})$ is the collection of all \textit{non-empty, finite} subsets of $\mathbf N_{\ge 2}$.
Lastly, $\{1\} \in \mathscr L(H)$ if and only if $\mathcal A(H) \ne \emptyset$, in which case we also have that $\bigl\{\{0\}, \{1\}\bigr\}$ is properly contained in $\mathscr L(H)$, because $k \in \mathsf L_H(a^k)$ for all $k \in \mathbf N^+$ and $a \in \mathcal A(H)$.

We say that $H$ is: \textit{reduced} if $H^\times = \{1_H\}$; \textit{BF} (short for ``bounded factorization'') if $1 \le |\mathsf L_H(x)| < \infty$ for every $x \in H \setminus H^\times$; \textit{unit-cancellative} if $xy = x$ or $yx=x$, for some $x, y \in H$, implies $y \in H^\times$ (this is a generalization of cancellativity). We call a function $\lambda: H \to \mathbf N$ a \textit{length function} if $\lambda(y) < \lambda(x)$ for all $x, y \in H$ such that $x = uyv$ for some $u, v \in H$ with $u \notin H^\times$ or $v \notin H^\times$.
\begin{thm}
Assume $H$ is a unit-cancellative monoid and admits a length function. Then $H$ is \textup{BF}.
\end{thm}
With these definitions in place (which are largely unnecessary, but put things in a certain perspective), take $\mathcal P_{{\rm fin}, 0}(\mathbf N)$ to be the set of all non-empty, finite subsets of $\mathbf N$ containing $0$ endowed with the (binary) operation of set addition
$$
H \times H \to H: (X, Y) \mapsto X + Y := \{x+y: x \in X,\,y \in Y\}.
$$
It is seen that $H$ is a reduced, unit-cancellative (but highly non-cancellative!), commutative BF-monoid, referred to as the \textit{restricted power monoid} of $(\mathbf N, +)$: The identity is the singleton $\{0\}$, and a length function is given by the map $H \to \mathbf N: X \mapsto |X| - 1$. Notably, ``most'' $X \in \mathcal P_{{\rm fin}, 0}(\mathbf N)$ are atoms, in the sense that, if $\alpha_n$ is the number of atoms contained in the discrete interval $[\fixed[-0.6]{\text{ }}[ 0, n ]\fixed[-0.6]{\text{ }}]$, then $\alpha_n/2^n \to 1$ as $n \to \infty$. We get from here and the above that
\begin{equation}
\label{equ:tringali(1)}
\mathscr L(\mathcal P_{{\rm fin},0}(\mathbf N)) \subseteq \bigl\{\{0\}, \{1\}\bigr\} \cup \mathcal P_{\rm fin}^\ast(\mathbf N_{\ge 2}),
\end{equation}
and we have the following.

\begin{conj}
The inclusion in \eqref{equ:tringali(1)} holds as an equality.
\end{conj}

In more explicit terms, we are asking whether, for every non-empty finite set $L \subseteq \mathbf N_{\ge 2}$, there exists $X \in \mathcal P_{{\rm fin},0}(\mathbf N)$ with the property that $L$ is equal to the set of all $k \in \mathbf N^+$ such that $X = A_1 + \cdots + A_k$ for some atoms $A_1, \ldots, A_k \in \mathcal P_{{\rm fin},0}(\mathbf N)$.

What we know so far:
\begin{enumerate}
\item $\mathscr L(\mathcal P_{{\rm fin},0}(\mathbf N))$ contains every one-element subset of $\mathbf N$ and every two-element subset of $\mathbf N_{\ge 2}$.
\item If $L$ is in $\mathscr L(\mathcal P_{{\rm fin},0}(\mathbf N))$, then so is $L + h$ for every $h \in \mathbf N$.
\item $[\fixed[-0.6]{\text{ }}[ 2, n ]\fixed[-0.6]{\text{ }}] \in \mathscr L(\mathcal P_{{\rm fin},0}(\mathbf N))$ for all $n \ge 2$.
\end{enumerate}
See \bburl{https://arxiv.org/abs/1701.09152} (in particular, Corollary 2.23 and Sections 4 and 5 therein) for further details.
%

\subsection{Problem Session II: Friday, May 26th (Chair Steven J. Miller)}

\subsubsection{Steven J. Miller: \textcolor{blue}{\href{mailto:sjm1@williams.edu}{sjm1@williams.edu}}}

Consider the sequence of papers of trying to generalize the works on Rankin on sequences avoiding three term arithmetic progressions.

\begin{itemize}

\item Joint work with Andrew Best, Karen Huan, Nathan McNew, Jasmine Powell, Kimsy Tor, Madeleine Weinstein: {\it Geometric-Progression-Free Sets over Quadratic Number Fields}. To appear in the  Proceedings of the Royal Society of Edinburgh, Section A: Mathematics. \bburl{https://arxiv.org/pdf/1412.0999v1.pdf}.

\item Joint with Megumi Asada, Eva Fourakis, Sarah Manski, Gwyneth Moreland and Nathan McNew: {\it Subsets of $\mathbb{F}_q[x]$ free of 3-term geometric progressions}. To appear in Finite Fields and their Applications. \bburl{https://arxiv.org/pdf/1512.01932.pdf}.

\item Joint with Megumi Asada, Eva Fourakis, Eli Goldstein, Sarah Manski, Gwyneth Moreland and Nathan McNew: {\it Avoiding 3-Term Geometric Progressions in Non-Commutative Settings}, preprint. Available at: \bburl{https://web.williams.edu/Mathematics/sjmiller/public_html/math/papers/Ramsey_NonComm2015SMALLv10.pdf}.

\end{itemize}

The speaker would love to work with others on various generalizations; if you are interested please contact him as he is possibly pursuing some of these with his REU students.

\begin{itemize}

\item What about avoiding 4 terms? Avoiding patterns?

\item Look at cubic or higher number fields; how do the answers depend on the properties of the field?

\item Look at matrix analogues. Could look at $M, MR, MR^2$. Maybe these matrices live in a group? Maybe we fix a point $\overrightarrow{v}$ and want to make sure there aren't three matrices $M_i$ such that $M_1\overrightarrow{v}, M_2 \overrightarrow{v}$ and $M_3\overrightarrow{v}$ do not form a geometric progression. Might have to be carefully as depending on the matrix size might be mapping vectors to vectors and not scalars.

\item Generalize the quaternion arguments to octonions. To sedenions?

\end{itemize}

\subsubsection{Mel Nathanson: \textcolor{blue}{\href{mailto:MELVYN.NATHANSON@lehman.cuny.edu}{MELVYN.NATHANSON@lehman.cuny.edu}}}

Suppose you have a finite set $A$ of non-negative integers. Let $n_2(A)$ be the largest integer $n$ such that $\{0,1,2,\dots,n\} \subset A+A$. Fix $A$ with $|A| = k$, and let $n_2(k)$ be the maximum of all $n_2(A)$ with $|A| = k$. Easy to check that $k^2/4 \le n_2(k) \le k^2/2$. Does the limit of $n_2(k)/k^2$ exist?

\subsubsection{Mel Nathanson: \textcolor{blue}{\href{mailto:MELVYN.NATHANSON@lehman.cuny.edu}{MELVYN.NATHANSON@lehman.cuny.edu}}}

Take a finite set of lattice points $A \subset \Z^n$. Let $K$ be the convex hull of $A$. Dilate: $h\ast K$. Count the number of lattice points inside: $$E_h(A) \ = \ |h\ast K \cap \Z^n|.$$ We know this equals ${\rm vol}(K)h^n$ plus lower order terms.

This is a continuous operation. Assume the subgroup generated by $A$ is all of $\Z^n$. As $K$ is convex, $h \ast K = hK$. As $A \subset K$ we have $hA \subset hK = h \ast K$.

Thus $hA \subset (h \ast K) \cap \Z^n$. It is not everything, but it is a lot. We have $|hA| = {\rm vol}(K)h^n$ plus lower order terms. We know almost nothing about this polynomial except that the leading coefficient is the volume. The difference between the two polynomials are the lower order terms. What can one say about the points missing or in the boundary layer? Have to have some geometrical description.

\subsubsection{Steven Senger: \textcolor{blue}{\href{mailto:StevenSenger@missouristate.edu}{StevenSenger@missouristate.edu}}}

Let $A$ be a subset of $[0,1]$, say $[.01, .02] \cup [.90, .91]$. Have .1 not in $A-A$. Can come up with a subset of the reals that avoid certain differences.

Try to take a step further. Can you find an $A \subset \R^2$ whose distance set $$\Delta(A) \ = \ \{|x-y|: x, y \in A\}$$ avoids certain values?

What is the behavior of a set $A \subset \mathbb{F}_q^2$ \emph{avoiding} some set $B \subset \mathbb{F}_q^2$ of ``distances''. Beliefs:

\begin{itemize}

\item You can't have a ``big'' set avoiding ``many'' distances.

\item If missing one distance you are missing ``many'' distances.

\end{itemize}

\subsubsection{Steven Senger: \textcolor{blue}{\href{mailto:StevenSenger@missouristate.edu}{StevenSenger@missouristate.edu}}}

Let $A \subset \R$ (or any field), assume $|A| < \infty$. Can there be be a ``geometric progression'' $G$, with $|G| \approx |GG|$, such that $$|(AA+1) \cap G|\ \approx\ |AA|\ \approx\ |G|.$$

Any energy techniques used have not been helpful.

Comment from Miller: replace the plus 1 with plus a constant, and see if a pigeon-hole argument can give that there must be at least one constant where this is true....

This is an \$8000 question. Note: typos here should not be attributed to the typist.

\subsubsection{Salvatore Tringali: \textcolor{blue}{\href{mailto:salvo.tringali@gmail.com}{salvo.tringali@gmail.com}} (preferred) or \textcolor{blue}{\href{mailto:salvatore.tringali@uni-graz.at}{salvatore.tringali@uni-graz.at}}}

Let $f: \mathcal{P}(\mathbf{N}) \to \mathbf{R}$ such that, for all $X, Y \subseteq \mathbf{N}$, $h \in \mathbf{N}$, and $k \in \mathbf{N}^+$, the following hold:
\begin{itemize}
\item $f^\ast(X) \le f^\ast(\mathbf{N}) = 1$.
\item $f^\ast(X \cup Y) \le f^\ast(X) + f^\ast(Y)$.
\item $f^\ast(k \cdot X + h) = f^\ast(X)/k$, where $k \cdot X + h = \{kx + h: x \in X\}$.
\end{itemize}

The class $\mathscr F$ of functions satisfying these conditions is ``large'' and includes various upper densities that are commonly encountered in Analysis and Number Theory. Most notably, the following are in $\mathscr F$:

\begin{itemize}

\item the upper $\alpha$-density (with $\alpha$ a fixed parameter $\ge -1$), given by
    $$
    f^\ast(X) \ := \ \limsup_n \frac{\sum_{i \in X \cap [\fixed[-0.6]{\text{ }}[1, n ]\fixed[-0.6]{\text{ }}]} i^\alpha}{\sum_{i \in [\fixed[-0.6]{\text{ }}[1, n ]\fixed[-0.6]{\text{ }}]}  i^\alpha}.
    $$
    The upper logarithmic ($\alpha = -1$) and upper asymptotic ($\alpha = 0$) densities are special cases.

\item the upper Banach density.

\item the upper P\'olya density.

\item the upper analytic density.

\item the upper Buck density.

\item any ``convex combination'' of the form $\sum_{i=1}^\infty a_i f_i$, where $f_i$ are functions satisfying the conditions and the coefficients $a_i$ are non-negative real numbers adding to $1$.

\end{itemize}

Now, set $f_\ast(X) := 1 - f^\ast(\mathbf{N}\setminus X)$ for every $X \subseteq \mathbf{N}$, and let $$\mathcal{D} \ := \ \{X \subseteq \mathbf{N}: f^\ast(X) = f_\ast(X)\}.$$

Denote by $f$ the restriction of $f^\ast$ to $\mathcal D$. Fix $A, B \in \mathcal D$, and let $\alpha \in [f(A), f(B)]$ (note that the interval may be empty). Does there exist an $X \in \mathcal D$ such that $A \subseteq X \subseteq B$ and $f(X) = \alpha$?

See \bburl{https://arxiv.org/abs/1506.04664} (in particular, Sections 2, 4, and 5) and \bburl{https://arxiv.org/abs/1510.07473} for further details.


\subsection{Speakers and Participants}

\bi

\item Ali Armandnejad, Vali-e-Asr University of Rafsanjan, Iran
\item Abdul Basit, Rutgers - New Brunswick
\item Sam Chow, University of York, England
\item David Chudnovsky, NYU
\item Gregory Chudnovsky, NYU
\item Colin Defant, University of Florida
\item  Robert Donley, Queensborough Community College (CUNY)
\item  Joseph Gunther, CUNY Graduate Center
\item  Brandon Hanson, Pennsylvania State University
\item  Charles Helou, Penn State
\item  Brian Hopkins, St. Peter's University
\item  Robert Hough, SUNY at Stony Brook
\item  Alex Iosevich, University of Rochester
\item  William J. Keith, Michigan Technological University
\item  Mizan Khan, Eastern Connecticut State University
\item  Byungchan Kim, SeoulTech, Republic of Korea
\item  Hershy Kisilevsky, Concordia University, Canada
\item  Sandor Kiss, Budapest University of Technology and Economics, Hungary
\item  Nana Li, Bard College at Simon's Rock
\item  Jared Lichtman, Dartmouth College
\item  Neil Lyall, University of Georgia
\item  Akos Magyar, University of Georgia
\item  Michael Maltenfort, Northwestern University
\item  Azita Mayeli, Queensborough Community College CUNY
\item  Nathan McNew, Towson State University
\item  Steven J. Miller, Williams College
\item  Mel Nathanson, Lehman College CUNY
\item  Mengquing Qin, Missouri State University
\item  Hans Parshall, University of Georgia
\item  Giorgis Petridis, University of Georgia
\item  Sinai Robins, University of Sao Paulo, Brazil
\item  Ryan Ronan, CUNY Graduate Center
\item  Csaba Sandor, Budapest Univ. of Technology and Economics, Hungary
\item  James Sellers, Pennsylvania State University
\item  Steve Senger, Missouri State University
\item  Satyanand Singh, New York City Tech CUNY
\item  Jonathan Sondow, New York
\item  Jack Sonn, Technion, Israel
\item  Yoni Stancescu, Afeka College, Israel
\item  Stefan Steinerberger, Yale University
\item  Salvatore Tringali, University of Graz, Austria
\item  Yuri Tschinkel, Courant Institute, NYU
\item  Ajmain Yamin, Bronx High School of Science
\item  Yifan Zhang, Central Michigan University
\ei


\newpage

\section{CANT Problem Sessions: 2018}

\subsection{Problem Session I: Tuesday, May 22th (Chair Steven J. Miller)}

\subsubsection{Steven J Miller: I (\email{\textcolor{blue}{\href{mailto:sjm1@williams.edu}{sjm1@williams.edu}}})}

In the talk today by Kevin Kwan and Steven J. Miller the asymptotic growth rate of $k$-near perfect numbers was given. These problems were inspired by the study of perfect numbers, which sadly have many intractable problems. Explicitly, do odd perfect numbers exist?

\subsubsection{Steven J Miller: II (\email{\textcolor{blue}{\href{mailto:sjm1@williams.edu}{sjm1@williams.edu}}})}

Consider the Fibonacci game, defined as follows. We have bins $B_i$, where the label of bin $i$ is $F_i$ ($\{F_i\}$ are the Fibonacci numbers, with $F_1 = 1, F_2 = 2$ and $F_{n+2} = F_{n+1} + F_n$). This two player game goes as follows. Start with a positive integer $N$, and have $N$ beans in bin $1$. The two players alternative moving. A legal move is one of the following.

\begin{itemize}

\item Combine one bean each from two consecutive bins, removing them from their bins and putting one in the next (thus if we have 4 beans in bin 6 and 9 bins in bin 7, we can combine and have 3 bins in bin 6, 8 bins in bin 7, and now 1 bin in bin 8). Note we also allow combining two beans from bin 1, removing them and adding one bean to bin 2.

\item Take two beans in bin $i$ and if $i \ge 3$ remove and add one bean to bin $i+1$ and 1 bean to bin $i-2$, or if $i=2$ remove and add one to bin 3 and one to bin 1. (We could also view combining two beans from bin 1 and getting one in bin 2 as this type of move.)

\end{itemize}

Kristen Flint and Steven J. Miller proved that if you make the game deterministic, you always end in the Zeckendorf decomposition and the number of moves is related to the number of summands in the Zeckendorf decomposition.

Alyssa Epstein and Steven J. Miller proved all games terminate, no matter how people move. Moreover, interestingly, unless $N=2$ then Player Two \emph{always} has a winning strategy. Unfortunately the proof is non-constructive, leading to....

\noindent \textbf{Question: Is there a constructive winning strategy for Player Two?}

Comment: could it depend on parity?

\ \\

\subsubsection{Jared Duker Lichtman: I (\email{\textcolor{blue}{\href{mailto:Jared.D.Lichtman@gmail.com}{Jared.D.Lichtman@gmail.com}}})}

Let $\sigma(n) = \sum_{d|n} d$. If $\sigma(n) = 2n$ then $n$ is perfect (if $\sigma(n)$ exceeds $2n$ then $n$ is called abundant,
else $n$ is called deficient).

An nondeficient number is called \textbf{primitive} if $d$ is deficient for all $d$ dividing $n$.

Are there any odd primitive nondeficient numbers? Up to 1000, there are 17 such numbers, only one of which (i.e. 945) is odd. So how
many are there up to $x$?

If $n$ is an even perfect number we know $n = 2^e p$ where $p = 2^{e+1}-1$ is a Mersenne prime. If we relax our condition on $p$ to
be a prime in $[2^e, 2^{e+1}]$ then $n=2^ep$ is primitive abundant. Can generalize: $n = 2^e p_1 \cdots p_k$ for distinct primes
$p_i \in [(k-1)2^{e+1}, k2^{e+1}]$.

To count even primitive nondeficient numbers up to $x$, Erd\H{o}s' idea was to choose $e$ and $k$ carefully so that $\binom{P}{k}$ is
large, where $P$ is the number of primes in the interval $[(k-1)2^{e+1}, k2^{e+1}]$. Specifically, taking $e \approx \sqrt{\log x
\log\log x}$ and $k \approx \sqrt{\log x/ \log \log x}$ gives a lower bound asymptotic to $x / \exp(c \sqrt{\log x \log\log x})$.
Note this is of the correct form, since there is an upper bound of the same form (though a different constant $c$).

\textbf{Can we do the same for odd primitive nondeficient numbers?} Note that for $n = 2^e p_1 \cdots p_k$,
\begin{eqnarray}\sigma(n) & \ = \ & \sigma(2^e)\prod \sigma(p_i) \nonumber\\ &=& (2^{e+1}-1) \prod (p_i + 1). \nonumber
\end{eqnarray}
Erd\H{o}'s argument leverages the fact that $\sigma(2^e)/2^e = 2 - 2^(-e) \to 2$ as $e\to\infty$. It seems unclear how to replicate
in the odd case.

\ \\

Another question is the following. We say a set $A$ of natural numbers is \textbf{primitive} no element of $A$ divides another. Many
familiar sets are primitive (for example, the primes). Consider the set of integers up to $n$, denoted $[1,n]$. How can you
partition this into primitive subsets? Of course, we can pick singleton sets, but that requires a lot of sets.

\noindent \textbf{How many disjoint primitive sets does it take to cover $[1,n]$?}

Note that the second-half of the interval $(n/2, n]$ is primitive, so we may iteratively take $(n/2, n]$, $(n/4, n/2], \ldots$, to
get a partition using $\approx\log_2 n$ primitive sets.

Can we do better? \textbf{Answer from audience: no!} Consider powers of 2, cannot have two in the same set.

Maybe we can generalize? What about restricting to the odd numbers in $[1,n]$, or square-free?

\ \\

\subsubsection{Arseniy (Senia) Sheydvasser: I  (\email{\textcolor{blue}{\href{mailto:sheydvasser@gmail.com}{sheydvasser@gmail.com}}})}

Ulam sequences $U(a,b)$: $a, b, u_2, u_3, \dots$ where $a, b$ are the two initial terms, and $u_k$ is the smallest integer that can be written as the sum of two distinct prior terms in exactly one way.

For example: $U(1,2)$: 1, 2, 3, 4, not 5 but 6, not 7 but 8, ....

Lots of open questions. We found numerical evidence supporting the following conjecture. Look at the family $U(1,n)$.

\noindent \textbf{Rigidity Conjecture: There exists a positive integer $N > 0$ and coefficients $a_i, b_i, c_i, d_i$ such that for all $n > N$, $U(1,n)$ is a disjoint union of $[a_i n + n_i, c_i n + d_i]$. We think $N=4$.}

For all $C>0$ there exists an $N>0$ and coefficients $a_i, b_i, c_i, d_i$ such that $U(1,n) \cap [1,Cn]$ is a disjoint union of $[a_i n + n_i, c_i n + d_i] \cap [1,cn]$. Sadly $C$ can depend on $N$.


\subsection{Problem Session II: Friday, May 24th (Chair Kevin O'Bryant)}

\subsubsection{Joel Moreira I: (\email{\textcolor{blue}{\href{mailto:joel.moreira@northwestern.edu}{joel.moreira@northwestern.edu}}})}

A set $S \subset \mathbb{N}$ is \emph{syndetic} if there exists $L \subset \mathbb{N}$ such that for any $n \in \mathbb{N}$, we have
\[
S \cup \{n, n + 1, \cdots, n + L \} \ \neq\ \emptyset.
\]
\textbf{(Beiglboeck, Bergelson, Hindman, Strauss) Do syndetic sets contain arbitrarily long geometric progressions? Do syndetic sets contain a subset of the form $\{a, an^2\}$?}

\ \\

\noindent \textbf{Problem: If $S$ is syndetic, does it contain \[A \cdot B\ =\ \{ab: a \in A, b \in B \}, \]
where $A$ and $B$ are infinite?}

\ \\

\subsubsection{Melvyn Nathanson: I (\email{\textcolor{blue}{\href{mailto:MELVYN.NATHANSON@lehman.cuny.edu}{MELVYN.NATHANSON@lehman.cuny.edu}}})}

Let  $A = (a_{i,j})$ be an $n \times n$ matrix such that
$a_{i,j}>0$ for all $i,j \in \{1,\ldots, n\}$.
The $i$th \emph{row sum} of $A$ is
\[
{\rm row}_i(A)\ =\ \sum_{j=1}^n a_{i,j}\ >\ 0.
\]
The matrix $A$ is \emph{row stochastic} if ${\rm row}_i(A) = 1$
for all $i \in \{1,\ldots, n\}$.

The $j$th \emph{column sum} of $A$ is
\[
{\rm col}_j(A)\  =\ \sum_{i=1}^n a_{i,j}\ >\ 0.
\]
The matrix $A$ is \emph{column stochastic} if ${\rm col}_j(A) = 1$
for all $j \in \{1,\ldots, n\}$.

The matrix $A$ is \emph{doubly stochastic} if it is
both row stochastic and column stochastic.

Define the $n \times n$ positive diagonal matrix
\[
X(A)\ =\ \diag \left( \frac{1}{{\rm row}_1(A)}, \frac{1}{{\rm row}_2(A)},\ldots, \frac{1}{{\rm row}_m(A)} \right).
\]
Multiplying $A$ on the left by $X(A)$ multiplies each coordinate in the $i$th
row of $A$ by $ 1/{\rm row}_i(A)$, and so
\begin{align*}
{\rm row}_i\left( X(A) A\right)
&\ =\ \sum_{j=1}^n (X(A) A)_{i,j}
 \ =\ \sum_{j=1}^n \frac{a_{i,j}}{{\rm row}_i(A)}
 \ =\ \frac{ {\rm row}_i(A)}{{\rm row}_i(A)}  = 1
\end{align*}
for all $i \in \{1,2,\ldots, n\}$.
The process of multiplying $A$ on the left by $X(A)$ to obtain the
row stochastic matrix $X(A) A$ is called
\emph{row scaling} or \emph{row normalization}.
Note that $X(A) A = A$ if and only if  $X(A) = I$ if and only if $A$ is row stochastic.
Also, the row stochastic matrix $X(A)A$ is not necessarily column stochastic.

Similarly, we define the $n \times n$ positive diagonal matrix
\[
Y(A)\ =\ \diag \left( \frac{1}{{\rm col}_1(A)}, \frac{1}{{\rm col}_2(A)},\ldots,
\frac{1}{{\rm col}_n(A)} \right).
\]
Multiplying $A$ on the right by $Y(A)$ multiplies
each coordinate in the $j$th column
of $A$ by $1/{\rm col}_j(A)$, and so
\[
{\rm col}_j(AY(A))\ =\ \sum_{i=1}^n (A Y(A))_{i,j}\ =\ \sum_{i=1}^n  \frac{a_{i,j}}{{\rm col}_j(A)}
\ = \ \frac{{\rm col}_j(A)}{{\rm col}_j(A)} \ =\  1
\]
for all $j \in \{1,2,\ldots, n\}$.
The process of multiplying $A$ on the right by $Y(A)$ to obtain a
column stochastic matrix $A Y(A)$ is called
\emph{column scaling} or \emph{column normalization}.
Note that $AY(A) = A$ if and only if  $Y(A) = I$ if and only if $A$ is column stochastic.
Also, the column stochastic matrix
$A Y(A)$ is not necessarily row stochastic.

The following algorithm is often called  ``alternate minimization,''
but it should more appropriately be called ``alternate normalization''
or ``alternate scaling.''

Let $A = (a_{i,j})$ be a positive $n \times n$ matrix.
Construct inductively an infinite sequence of positive $n \times n$
matrices by alternate operations of row scaling and column scaling:
\begin{align*}
A^{(0)}&\ = \  A \\
A^{(1)} &\ = \  X( A^{(0)} )   \cdot  A^{(0)} \\
A^{(2)} &\ = \  A^{(1)}  \cdot  Y( A^{(1)} ) \\
A^{(3)} &\ = \  X( A^{(2)} )   \cdot A^{(2)} \\
A^{(4)} &\ = \  A^{(3)}  \cdot  Y( A^{(3)} ) \\
A^{(5)} &\ = \  X( A^{(4)} )   \cdot A^{(4)} \\
A^{(6)} &\ = \  A^{(5)}  \cdot  Y( A^{(5)} ) \\
& \vdots
\end{align*}
Sinkhorn's theorem states that this sequence of matrices converges to a doubly
stochastic matrix.

If $A^{(L)}$ is doubly stochastic for some $L$, then $A^{(\ell)} = A^{(L)}$
for all $\ell \geq L$. In this presumably exceptional case, we say that the algorithm
terminates in at most $L$ steps.  For $2\times 2$ matrices, if the algorithm terminates in at most $L$ steps, then the algorithm terminates in at most two steps.

\ \\

\noindent \textbf{Problem:  Let $n \geq 3$. Does there exist an integer $L^{\ast}(n)$ such that, if $A$ is a positive $n\times n$ matrix for which the alternate minimization algorithm terminates in $L$ steps, then the alternate minimization algorithm
terminates in $L^{\ast}(n)$ steps?}

\ \\

\subsubsection{Steve Senger: I: (\email{\textcolor{blue}{\href{mailto:StevenSenger@MissouriState.edu}{StevenSenger@MissouriState.edu}}})}

\noindent \textbf{Problem: Given a large finite set $A \subset \mathbb{R}$ with $|A| = n$, can $(A \cdot A + 1)$ have large intersection with a generalized geometric progression of length about size $|A \cdot A + 1|$ with finite number of generators? Can $(A \cdot A + 1)$ have large intersection with $B \cdot C$, where $B$ and $C$ are both large?}

\ \\

\subsubsection{Brad Isaacson: I: (\email{\textcolor{blue}{\href{mailto:BIsaacson @citytech.cuny.edu}{BIsaacson @citytech.cuny.edu}}})}

Eisenstein first come with the observation that for
\[
C_k(x) \ =\ \sum_{n \in \mathbb{Z}} \frac{1}{(n + x)^k},
\]
we have
\[
c_1(x) \ = \  \pi \cot(\pi x),
\]
\[
c'_1(x) \ = \  - c_2(x) \ = \  - \pi^2 (1 + \cot^2(\pi x)),
\]
and thus
\[
c_2(x) \ = \  \pi^2 +(c_1(x))^2,
\]
\[
c_3(x) \ = \  c_1(x)c_2(x) \ = \  c_1^3(x) + \pi^2 c_1(x).
\]
We have a pattern that $c_k(x) = v_k(c_1(x))$ for a polynomial $v_k$,
where
\[ v_0(t) \ = \  1, \]
\[v_1(t) \ = \  t, \]
\[k v_{k + 1}(t) \ = \  (t^2 + \pi^2) v_k' \]
for $k \geq 1$.

\ \\

\noindent \textbf{Problem: Is there a closed form expression for $v_k$? We can use either Bernulli numbers, or Sterling numbers, or anything else. Note that the Bernulli polynomials satisfy a similar equation $k B_k = B_{k + 1}'$.}

\ \\

\subsection{Speakers and Participants}

See \bburl{http://www.theoryofnumbers.com/CANT2018-speakers.pdf}.

\ \\


\newpage

\section{CANT Problem Sessions: 2019}

\subsection{Problem Session I: Tuesday, May 21th (Chair Kevin O'Bryant)}

\subsubsection{Colin Defant}

For $c \in \mathbb{C}$, let \[ \sigma_c(n)\ :=\ \sum\limits_{d\vert n} d^c.\]

Now let \[\sigma_c(\mathbb{N})\ :=\ \{ \sigma_c(n): n \in \mathbb{N} \}\] and \[ \mathcal{N}(c)\ =\ \#\text{ of connected components of } \overline{\sigma_c(\mathbb{N})},\] where the overline denotes topological closure. One should note that $\overline{\sigma_c(\mathbb{N})}$ is a bounded set if and only if $\Re(c) < - 1$. We have the following theorem.
\ \\

\noindent \textbf{Theorem [Defant]: If $\Re(c) < -1$ then $\mathcal{N}(c) < \infty$.} \\ \

\noindent \textbf{Question:} Does there exist $c \in \mathbb{C} \setminus \mathbb{R}$ with $\Re(c) < -1$ and $\mathcal{N}(c) = 1$? \\ \

\noindent \textbf{Question:} For which bounded multiplicative functions $f : \mathbb{N} \longrightarrow \mathbb{C}$ is it the case that $\overline{f(\mathbb{N})}$ has finitely many connected components?\\

\subsubsection{Colin Defant}
With the same notation as above, let \[E_m \ :=\ \{ c \in \mathbb{R}: \mathcal{N}(c) = m\}.\]

\ \\

\noindent \textbf{Theorem [Zubrilina 2017]: $E_4 = \emptyset$.}

\ \\

\noindent \textbf{Theorem [Achenjang, Berger 2018]: $E_6 = \emptyset$}

\ \\

We say that $m \in \mathbb{N}$ is a Zubrilina number if $E_m = \emptyset$.

\ \\

\noindent \textbf{Conjecture (Defant): There are infinitely many Zubrilina numbers.} \ \\

The methods for $E_4$ and $E_6$ reduce to a somewhat algorithmic check of various inequalities, which motivates the following: \\
\textbf{Question:} Can we find an algorithm to compute the Zubrilina numbers?\\

\subsubsection{Satyanand Singh}
$\lambda_{2,3}(h)$ is the smallest positive integer that can be represented as the sum of $h$, and no less than $h$, elements from the set \[ A_{2,3} \ :=\ \{0\} \cup \{ \pm 2^j: j = 0,1,2, \dots\} \cup \{ \pm 3^k: k = 0,1,2,\dots\}.\]

\ \\

\noindent \textbf{Question:} Find $\lambda_{2,3}(h)$ for $h \ge 5$.\\ \

For context, it is known that $\lambda_{2,3}(1) = 1$, $\lambda_{2,3}(2) = 5$, $\lambda_{2,3}(3) = 21$ and $\lambda_{2,3}(4) = 150$. Certain Maple calculations also suggest that $\lambda_{2,3}(5)$ might be $2581$, as this number seems to be always left out when one case-checks various combinations of powers. A useful tool to attack this problem might be Baker's method for linear forms in logarithms of primes. Tijdemann-Hadouche showed that $\lambda_{2,3}(h) \rightarrow \infty$ and $h\rightarrow \infty$.

\subsection{Problem Session II: Wednesday, May 22nd (Chair  David Grynkiewicz)}

\subsubsection{Andrew Odesky}

\subsubsection{Andrew Odesky}

Let $H$ be a finite group and $k$ a number field. A problem I'm interested in leads to a Diophantine problem of
finding $k$-points on a certain variety $X=X_H$, known to be a homogenous space for some
$G$: $X = G /H$. The variety $X$ is canonically associated to $H$. The group $G$ is the unit group of
the group algebra of $H$ over $k$, which is a connected, linear algebraic group over $k$.

There is an explicit description of $X$ as an intersection of quadrics:
$X = \cap_i Q_i \subset \mathbb{A}^d$. Here
$d$ depends on the size of $n = |H|$ by $d_n = \frac12 n^2 + \frac32 n$,
and $\dim X = n$.

It is provably the case that $X(k) \neq \varnothing$. \\ \

\textbf{Question:} What can be said about $X(k)$ in general when $k$ is a number field? \\ \

It might be impossible to say anything substantial with current technology.
Maybe the Circle Method applies, I don't know.

This problem is motivated by the Inverse Galois Problem. It is already well-known
that inverse Galois problems can be reformulated into Diophantine problems,
and this has been a reasonably successful approach in general. Our Diophantine
problem here is another instance of this. The hope is that the description, which is
fairly concrete and specific, will enable more tools to be applied.

In a precise sense if $X$ has enough rational
points (if $X$ has the "Hilbert property over $k$") then
the inverse Galois problem holds for $H$ over $k$.
For example: When $H$ is abelian then $X$ has lots
of rational points (it is a torus), and has the Hilbert property over $k$.

\subsubsection{Fred Schneider}

For a more detailed description, see
\bburl{https://docs.google.com/spreadsheets/d/1FQm8LunkZz2oGDZnl_soZCJOVWqv4H_rFMIe-Cyet30/edit?ts=5ce8070d\#gid=1034210868}.
Trying to determine set counts over ranges of consecutive positive
integers.

Inspired by puzzle \#670 on \bburl{https://www.primepuzzles.net/}, try to
find three numbers in a row that have the same sum of divisors. When
trying to tackle this problem, I thought of it in terms of finding all the solutions $x$ to $f(x) = n$.
Imagine that there are a large number of solutions and it's more efficient to divide them into different subsets.
After that, we must find valid sequences of these subsets for possible solutions over a consecutive range.
When allowing for duplicates, this is a relatively simple problem but the much more interesting problem
is finding the number of distinct sets.

Let's start with considering how to generate all of the sets.  What sort of behavior must they follow?

For two consecutive numbers, the first is even, and the second odd (or vice
versa). So, there is one unique set.  We can label this (1,2).

For 3 consecutive numbers:

Powers of 2: Consider string $4^\ast, 1, 2$, shorthand for
a multiple of 4 or more, an odd number,  a multiple of 2 but not 4.
You could also have $1, 4, 1$ and $2, 1, 4$ and $1, 2, 1$.
For powers of 3, the possibilities are  $3, 1, 1$
and $1, 3, 1$ and $1, 1, 3$. "3" means any multiple of 3.  ``1'' means a non-multiple of 3.

If we overlay (aka take the cross-product of) the ``2'' and ``3'' rows, we will get
seven sets out of 4*3=12 combinations.  Sorting each set's elements, we get:
(1,1,6), (1,1,12), (1,2,3), (1,2,12), (1,3,4), (1,4,6), (2,3,4)

One thing that jumps out is that if we overlay the 3's
on top of the $4^\ast, 1, 2$ get three unique sets, but over $1, 4, 2$ get
only 2, over $2, 1, 4$ get zero, and over $1, 2, 1$ get 2 again.
When considering symmetry, we can prune redundant sets.

As you go further out, you have to consider the effect of new primes and larger prime powers.  For instance, for length 5,
you must start considering $8^\ast$ in the ``2'' rows.

I am trying to come up with an efficient way to enumerate all the distinct sets. It's easy to find the total number
of distinct tuples (like a Cartesian product) but the distinct set count is a challenge.
I wrote a C++ program to do this.  Several optimizations are explained in the Google sheet and code.
But, as the range length increases, I run into memory and storage issues.

\subsubsection{Noah Kravitz}

Related to Ryan Matzke's talk from earlier. Ryan asked what is the biggest $(k,\ell)$-sum-free set....

Let $G$ be a compact abelian group with Haar measure $\mu$. Let $$\lambda_{k,\ell} \ = \ \sup\left\{\mu(A): A \subset G\ {\rm satisfies} \ kA \cap \ell A = \emptyset\right\}.$$

If $G$ is the direct product of an identity component and a discrete part, say $\mathbb{I} \times M$, then in a certain sense can break up and take the maximum of the two constants: $$\lambda_{k,\ell}(G) \ = \ \max\left\{\lambda_{k,\ell}(\mathbb{I}), \lambda_{k,\ell}(M)\right\}.$$ Moreover if $\mathbb{I}$ is non-trivial then $$\lambda_{k,\ell} \ = \ \frac1{k + \ell}.$$

Question: What about $\lambda_{k,\ell}(M)$? What if $M$ is profinite, say $M$ is the direct sum of a bunch of cyclic groups.

Can always do a generalization of a middle third argument.... Can lift results back....

From Ryan Matzke: perhaps a lower bound is $$\left\lfloor \frac{|G|- \gcd(|G|,k-\ell) - 1}{k+\ell} \right\rfloor.$$

\subsubsection{Paolo Leonetti}

Take any function: $f: \mathcal{P}(\mathbb{N}) \to [0,1]$ such that

\begin{itemize}
\item $f(\mathbb{N}) = 1$.

\item $f(X) \le f(Y)$ if $X \subset Y$.

\item $f(X \cup Y) \le f(X) + f(Y)$.

\item $f(k \cdot X) = \frac1{k} f(X)$ for all $k \ge 1$.

\item $f(X+1) = f(X)$.

\end{itemize}

For all $x \in [0,1]$ does there exist an $A \subset \mathbb{N}$ such that $f(A+A) = x$? In other words, is this true for all such $f$? It is known to hold for certain special $f$. True for upper asymptotic and upper logarithmic density.

\subsubsection{Paul Baginski (Steve asked Paul to ask so he can type)}

Inspired by David Grynkiewicz's talk. Take a finite abelian group $G$, take a sequence of elements $S = g_1 \cdot g_2 \cdots g_n$ of elements in  $G$. Instead of just adding allow addition and subtraction. So consider $$\left\{\sum_{|I| \subset [1,n]} \lambda_i g_i: \lambda_i \in \{-1,1\}\right\}.$$ Considering all subsequences $T$ of $S$, have a choice of adding or subtracting each term in $T$. When does zero belong to this set? If $S$ is long enough then 0 will be in. Shortest length ensuring 0 is in is called the plus-minus Davenport constant. It is the least $\ell$ such that any sequence of length at least $\ell$ has a plus-minus zero sum subsequence.

Regular Davenport constant: certain values are known, not known in general.

The plus minus Davenport constant is more recent and less studied. Some general bounds. Adhikari-Grynkiewicz-Sun have some really basic upper bounds that are close to the lower bounds (much closer than the Davenport constant, where there is more of a difference). Exact values in a couple of cases. If $G$ is a 2-group then $D_{\pm}(G) = \log_2 |G|$. If $G = \mathbb{Z}_3^r$ then $D_\pm(G) = r+1$, while if $G$ is cyclic then $D_\pm(\mathbb{Z}_n)$ is either the floor or the ceiling of $\log_2 n$ (think ceiling).

Question: What is the structure of the longest sequence that avoid a plus minus sum of zero? What is the structure of the shortest sequences with a plus minus zero sum but no proper plus minus zero sum? For $G = \mathbb{Z}_2^r$ the answer to the first question is just a basis, and for question 2 is a basis union the sum of the basis. If $G = \mathbb{Z}_3^r$ believe for the first it is still a basis, and for the second it is a basis with union of the sum \emph{or} could have a basis with union of the negative of the sum. If $G = \mathbb{Z}_n$ it is open.

\ \\

\subsubsection{David Grynkiewicz}

Let $G$ be a (finite) abelian group and take a bunch of cardinality 2 subsets. Consider cardinality 2 subsets $A_1, \dots, A_n$. Suppose $$\sum_{i=1}^n A_i \ \ \ {\rm is\ aperiodic};$$ this means $$\{x \in G: x + \sum A_i \ = \ \sum A_i\} \ = \ {\rm trivial}.$$

Want a unique element in $$A_1 + \cdots + A_{n-1} + A_n.$$ Is there a unique expression element? Does there exist a unique expression element in $A_1 + \cdots + A_n$ so that $a_i \in A_i$ such that there is no other way to represent $a_1 + \cdots + a_n$ in $A_1 + \cdots + A_n$? This came up once in a paper, I found a way around and didn't need the result, but are there clever constructions / ways to prove?

Just did 2 subsets as that is what needed at the time, with larger sets more ways to have it.

What if try just cyclic groups? What if do $\mathbb{Z}_p$?

\ \\

\emph{Comment from Noah Kravitz:  For $\mathbb{Z}/12\mathbb{Z}$, for example, take 3 copies of $\{0,2\}$ and 2 copies of $\{0,3\}$ for a counter-example}.

\subsection{Problem Session III: Thursday, May 23rd (Chair Steven J Miller)}

\subsubsection{Steven J. Miller and Tudor Popescu}

Inspired by trying to walking to infinity along the primes, what happens if you try to walk to infinity along the square-free numbers? Looking at the number of steps you can walk and stay square-free converges fairly quickly to close to what a random model predicts, assuming all numbers are equally likely to be square-free.

    \begin{itemize}
    \item Is it possible to walk to infinity by appending a bounded number of digits to a prime at each stage while staying prime?
    \item Is it possible to walk to infinity by appending a bounded number of digits to a square-free number at each stage while staying square-free? \end{itemize}

    We believe that the answers are no for primes, and yes for square-free numbers, since the density of square-free integers is positive. This is suggested by looking at the probabilistic models for the longest prime walk and square-free walk respectively.

    We ran a $10000$ sample simulation, starting at a random $9$ digit square free number, and randomly appending a digit and seeing how long we could continue. The results are as follows:
\begin{center}
 \begin{tabular}{||c c c||}
 \hline
 Length & Expected & Actual \\
 \hline\hline
 0 & 3921 & 3878 \\
 \hline
 1 & 2384 & 2233 \\
 \hline
 2 & 1449 & 1426 \\
 \hline
 3 & 881 & 918 \\
 \hline
 4 & 536 & 585 \\
 \hline
 5 & 326 & 362 \\
 \hline
 $\ge$ 6 & 503 & 503 \\
 \hline
\end{tabular}
\end{center}

We get that $\chi^2 = 38.352,$ and the null hypothesis will be rejected. Thus perhaps we are seeing the main term, and there are some dependencies in being square-free....

\subsubsection{Yunied Puig de Dios}

Generally speaking, if we want to produce a subset $A$ of $\Z_+$ with some property and we are lucky enough that $A$ can be
realized as a recurrence set of some bounded linear operator on some separable infinite-dimensional Banach space $X$, in the
sense that for some $x\in X$ and some non-empty open subset $U$ of $X$, the set $A=\{n\geq 0: T^nx\in U\}$, then we are in a
very advantageous situation. Indeed, this means that $A$ might enjoy much richer properties inherited from $T$. Then we might
prefer to deal with a more sophisticated object like $T$, but much more easy to deal with. Linear dynamics is a relatively young
area on the intersection of functional analysis and operator theory considered to be born in 1982. By now, the area has
developed enough in such a way that we have at our disposal a whole machinery which makes desirable to be in the above-mentioned
situation. We showed recently that any subset satisfying the conclusion of K\v r\'{i}\v z's Theorem is not that exotic as it can
be realized as the recurrence set of some bounded linear operator on $c_0(\Z_+)$. Is this the case of any other subset of the
integers known as the solution of some well-known problem in combinatorics?




\subsubsection{Arseniy (Senia) Sheydvasser}

      The Ulam sequence starts with 1, 2, and subsequent terms are the smallest that can be written as a sum of two distinct prior terms in exactly one way. So the
      first few terms are $$U_{1,2} \ = \  \{1, 2, 3, 4, 6, 8, 11, 13, 16, \dots \}.$$ How quickly does this sequence grow? Numerical evidence suggests linearly. Best
      known bound is exponential, so very far off from what we believe is the truth.

      More generally, can consider $$U_{1,n} \ = \ \{1, n, n+1, n+2, \dots \}.$$

      A remarkable theorem (Hinman, Kuca, Schlesinger, S.: 2018): there are $a_i, b_i, c_i, d_i$ such that for all $C > 0$ there is an integer $N$ such that if $n \ge
      N$ then $$U_{1,n} \ = \  \bigsqcup_{i=1}^\infty [a_i n + b_i, c_i n + d_i] \ \bigcap \ [1, Cn].$$ The proof is non-constructive and uses non-standard methods.

      Argues by assuming $n$ is a hyper-integer....

      Suggests a more general setting for questions like this. Suppose that we have an integer sequence $S_n$. Let $\mathcal{A}$ be an algorithm such that if I feed
      into $\mathcal{A}$ the index $n$ and cut-off $k$ then we get the first $k$ terms of the sequence. Suppose that $\mathcal{A}$ is still an algorithm if we allow
      $(n,k)$ to be non-standard integers. Can use an ultrafilter, but what you get might not be an algorithm (might not run in finite time). What can we say about
      $S_n$? What does the existence of an algorithm tell us in general? There are some examples known -- Hofstadter $Q$-sequence (Nathan Fox in 2018). Also
      Sidon-like sequences (Kevin O'Bryant). Are there other interesting examples?

\subsubsection{Steven Senger}

I've been asking this for ten years: Does there exist a large finite $A \subset \mathbb{R}$ such that $$|(AA+1)(AA+1)| \ \le \ |AA|$$ (the dumbest thing you can do with energy doesn't work...)? (See for example the 2011 problems.)

\subsection{Problem Session IV: Friday, May 24th (Chair Misha Rudnev)}

\subsubsection{Misha Rudnev}

Let $A = \{(a,b): a \neq 0\}$, $x \mapsto ax + b$ be a finite set of affine transformations of the real plane. Denote $\ell = (a,b) \in
A$. We want to find the energy of $A$: $$E(A) \ = \ \#\{\ell_1^{-1} \ell_2
= \ell_3^{-3} \ell_4: \ell_i \in A\}.$$ This question has a nice geometric
interpretation (unless the four points are collinear): if $l_i$ are viewed as points in the $(a,b)$-plane, then the lines $(l_1l_2)$ and $(l_3l_4)$ are parallel, while the lines $(l_1l_3)$ and $(l_2l_4)$  intersect the forbidden vertical axis $a=0$ at the same point.

Many ``abelian'' problems of additive combinatorics deal with energy estimates. In this case, we are interested in non-commutative energy.

Clearly $E(A)$ can be as large as $|A|^3$ if the whole set $A$ lies on a line. Let us forbid this, by restricting that, say at most $\sqrt{|A|}$ points can be collinear or just count the contribution into $E(A)$ of the above-described, non-collinear quadruples $(l_1,\dots,l_4)$, call it $E^*(A)$ and seek a better upper bound than the trivial $|A|^3$.

If so, then it follows from the famous Guth-Katz incidence theorem (adapted as Theorem 5, \bburl{https://arxiv.org/pdf/1412.2909.pdf}) that each individual one of the two above geometric  conditions allows for  $O(|A|^3\log|A|)$ quadruples $(l_1,...,l_4)$. If the two conditions were in some sense ``independent'', then one would expect $E^*(A) \lesssim |A|^2$ (up to a log factor). But this is too good to be true. Considering $A=S\times S$, where $S$ is either an arithmetic or geometric progression shows that in this case $E(A)\gg |A|^{5/2}$. Furthermore, in the general special Cartesian product case $A=S\times S$, one can prove that $E(A)\ll |A|^{11/4 - c},$ for some small $c>0$  (\bburl{https://arxiv.org/pdf/1812.01671.pdf}, Lemma 21).

The question is whether (i) $E^*(A)\lesssim |A|^{5/2}$ (possibly with some log factors) is the correct general upper bound, and less ambitiously, whether any {\em quantitative} non-trivial upper bound can be established for a general $A$ (not a Cartesian product). The only thing we know so far is that if at most $|A|^{1-\delta}$ points of $A$ are allowed to be collinear, $E^*(A)\ll |A|^{3-c}$, for some small $c>0$,  (\bburl{https://arxiv.org/pdf/1812.01671.pdf}, Corollary 13).

\subsubsection{Paolo Leonetti}

Is it true that for all positive sequences $x_n$ such that
\begin{itemize}

\item $\lim_{n\to\infty} x_n = 0$,

\item $x_n$ is decreasing,

\item $\sum_{n\ge 1} x_n = \infty$

\end{itemize}
then there is a $\theta \in (0,1)$ such that for all $\ell > 0$ there is a sub-sequence $\{x_{n_k}\}$ such that $\sum_{k \ge 1} x_{n_k} = \ell$ and $x_{n_k} \ll \theta^k$? The answer is true if also assume $\liminf_{n\to\infty} x_{n+1}/x_n > 1/2$.

\subsection{Speakers and Participants Lists}

\ \\

\noindent CANT 2019: Jean Bourgain Memorial Speakers

\begin{itemize}
\item  Alex Iosevich, University of Rochester
\textcolor{blue}{\href{mailto:iosevich@gmail.com}{iosevich@gmail.com}}

\item Alex Kontorovich, Rutgers University, New Brunswick
\textcolor{blue}{\href{mailto:alex.kontorovich@rutgers.edu}{alex.kontorovich@rutgers.edu}}

\item Ben Krause, California Institute of Technology
\textcolor{blue}{\href{mailto:benkrause2323@gmail.com}{benkrause2323@gmail.com}}

\item Neil Lyall, University of Georgia
\textcolor{blue}{\href{mailto:lyall@math.uga.edu}{lyall@math.uga.edu}}

\item Akos Magyar, University of Georgia
\textcolor{blue}{\href{mailto:amagyar@uga.edu}{amagyar@uga.edu}}

\item Mariusz Mirek, Rutgers University, New Brunswick
\textcolor{blue}{\href{mailto:m.a.mirek@gmail.com}{m.a.mirek@gmail.com}}

\item  Yumeng Ou, Baruch College (CUNY)
\textcolor{blue}{\href{mailto:yumeng.ou@baruch.cuny.edu}{yumeng.ou@baruch.cuny.edu}}

\item  Misha Rudnev, University of Bristol, Bristol, UK
\textcolor{blue}{\href{mailto:m.rudnev@bristol.ac.uk}{m.rudnev@bristol.ac.uk}}

\item Fred Schneider \textcolor{blue}{\href{mailto:fws.nyc@gmail.com}{fws.nyc@gmail.com}}

\item  Adam Sheffer, Baruch College (CUNY)
\textcolor{blue}{\href{mailto:adam.sheffer@baruch.cuny.edu}{adam.sheffer@baruch.cuny.edu}}

\item  Ilya Shkredov, Steklov Institute of mathematics, Moscow, Russia
\textcolor{blue}{\href{mailto:ilya.shkredov@gmail.com}{ilya.shkredov@gmail.com}}

\item  Jozsef Solymosi, University of British Columbia, Canada
\textcolor{blue}{\href{mailto:solymosi@math.ubc.edu}{solymosi@math.ubc.edu}}

\item Van Vu, Yale University
\textcolor{blue}{\href{mailto:van.vu@yale.edu}{van.vu@yale.edu}}

\item  Hong Wang, MIT
\textcolor{blue}{\href{mailto:hong.wang1991@gmail.com}{hong.wang1991@gmail.com}}

\end{itemize}

\ \\

\noindent CANT 2019 Speakers

\begin{itemize}

\item Paul Baginski, Fairfield University
\textcolor{blue}{\href{mailto:baginski@gmail.com}{baginski@gmail.com}}

\item Amanda Burcroff, University of Michigan
\textcolor{blue}{\href{mailto:burcroff@umich.edu}{burcroff@umich.edu}}

\item  Colin Defant, Princeton University
\textcolor{blue}{\href{mailto:cdefant@math.princeton.edu}{cdefant@math.princeton.edu}}

\item  Robert W. Donley, Jr., Queensborough Community College (CUNY)
\textcolor{blue}{\href{mailto:RDonley@qcc.cuny.edu}{RDonley@qcc.cuny.edu}}

\item  Heidi Goodson, Brooklyn College (CUNY)
\textcolor{blue}{\href{mailto:Heidi.Goodson@brooklyn.cuny.edu}{Heidi.Goodson@brooklyn.cuny.edu}}

\item  David Grynkiewicz, University of Memphis
\textcolor{blue}{\href{mailto:diambri@hotmail.com}{diambri@hotmail.com}}

\item  Matthew Hase-Liu, Harvard University
\textcolor{blue}{\href{mailto:matthewhaseliu@college.harvard.edu}{matthewhaseliu@college.harvard.edu}}

\item Charles Helou, Penn State, Brandywine
\textcolor{blue}{\href{mailto:cxh22@psu.edu}{cxh22@psu.edu}}

\item  Robert Hough, SUNY - Stony Brook
\textcolor{blue}{\href{mailto:robert.hough@stonybrook.edu}{robert.hough@stonybrook.edu}}

\item  Jing-Jing Huang, University of Nevada, Reno
\textcolor{blue}{\href{mailto:jingjingh@unr.edu}{jingjingh@unr.edu}}

\item Trevor Hyde, University of Michigan
\textcolor{blue}{\href{mailto:tghyde@umich.edu}{tghyde@umich.edu}}

\item  Mizan R. Khan, Eastern Connecticut State University
\textcolor{blue}{\href{mailto:KHANM@easternct.edu}{KHANM@easternct.edu}}

\item Sandor Kiss, Budapest University of Technology and Economics, Hungary
\textcolor{blue}{\href{mailto:ksandor@math.mme.hu}{ksandor@math.mme.hu}}

\item Noah Kravitz, Yale University
\textcolor{blue}{\href{mailto:noah.kravitz@yale.edu}{noah.kravitz@yale.edu}}

\item  Paolo Leonetti, Graz University of Technology, Austria
\textcolor{blue}{\href{mailto:leonetti.paolo@gmail.com}{leonetti.paolo@gmail.com}}

\item  Huixi Li, University of Nevada, Reno
\textcolor{blue}{\href{mailto:huixil@g.clemson.edu}{huixil@g.clemson.edu}}

\item  Ariane Masuda, New York City Tech (CUNY)
\textcolor{blue}{\href{mailto:AMasuda@citytech.cuny.edu}{AMasuda@citytech.cuny.edu}}

\item  Ryan W. Matzke, University of Minnesota - Twin Cities
\textcolor{blue}{\href{mailto:matzk053@umn.edu}{matzk053@umn.edu}}

\item  Steve Miller, Williams College
\textcolor{blue}{\href{mailto:sjm1@williams.edu}{sjm1@williams.edu}}

\item Akshat Mudgal, University of Bristol
\textcolor{blue}{\href{mailto:am16393@bristol.ac.uk}{am16393@bristol.ac.uk}}

\item  Rishi Nath, York College (CUNY)
\textcolor{blue}{\href{mailto:rnath@york.cuny.edu}{rnath@york.cuny.edu}}

\item  Mel Nathanson, Lehman College (CUNY)
\textcolor{blue}{\href{mailto:melvyn.nathanso@lehman.cuny.edu}{melvyn.nathanso@lehman.cuny.edu}}

\item  Timothy Newlin, United States Military Academy, West Point
\textcolor{blue}{\href{mailto:timothy.newlin@westpoint.edu}{timothy.newlin@westpoint.edu}}

\item  Kevin O'Bryant, College of Staten Island (CUNY)
\textcolor{blue}{\href{mailto:obryant@gmail.com}{obryant@gmail.com}}

\item  Andrew Odesky, University of Michigan
\textcolor{blue}{\href{mailto:aodesky@umich.edu}{aodesky@umich.edu}}

\item  Jun Seok Oh, University of Graz, Austria
\textcolor{blue}{\href{mailto:junseok.oh@uni-graz.at}{junseok.oh@uni-graz.at}}

\item  Nikita Pereverzin, United States Military Academy, West Point
\textcolor{blue}{\href{mailto:nikita.pereverzin@westpoint.edu}{nikita.pereverzin@westpoint.edu}}

\item  Cosmin Pohoata, California Institute of Technology
\textcolor{blue}{\href{mailto:apohoata@caltech.edu}{apohoata@caltech.edu}}

\item  Yunied Puig de Dios, University of California, Riverside
\textcolor{blue}{\href{mailto:puigdedios@gmail.com}{puigdedios@gmail.com}}

\item  Alex Rice, Millsaps College
\textcolor{blue}{\href{mailto:arice2386@gmail.com}{arice2386@gmail.com}}

\item Alisa Sedunova, MPIM Bonn, Germany
\textcolor{blue}{\href{mailto:alisa.sedunova@phystech.edu}{alisa.sedunova@phystech.edu}}

\item  Steven Senger, Missouri State University
\textcolor{blue}{\href{mailto:StevenSenger@MissouriState.edu}{StevenSenger@MissouriState.edu}}

\item  Satyanand Singh, New York City Tech (CUNY)
\textcolor{blue}{\href{mailto:SSingh@citytech.cuny.edu}{SSingh@citytech.cuny.edu}}

\item  Jonathan Sondow, New York
\textcolor{blue}{\href{mailto:jsondow@alumni.princeton.edu}{jsondow@alumni.princeton.edu}}

\item  Wenbo Sun, Ohio State University
\textcolor{blue}{\href{mailto:sun.1991@osu.edu}{sun.1991@osu.edu}}

\item  Aled Walker, Trinity College, Cambridge, UK
\textcolor{blue}{\href{mailto:aledwalker@gmail.com}{aledwalker@gmail.com}}

\end{itemize}

\section{CANT Problem Sessions: 2020}

\section{CANT Problem Sessions: 2021}

\section{CANT Problem Sessions: 2022}


\newpage

\section{CANT Problem Sessions: 2023}

\noindent \large \textbf{CANT problem session, May 24, 2023} \normalsize \\

\noindent Communicated by Steve Senger \\ \

\textbf{Problem 1a, communicated by Noah Kravitz:} Given a set $A\subset \mathbb R / \mathbb Z,$ with the property that $A + A - 2A \neq \mathbb R / \mathbb Z,$ where $2A$ denotes dilation by 2, not a sum set, what bounds can we have on the Lebesgue measure of $A$? It is known that $\mu(A)\leq \frac{1}{3}-\epsilon,$ for some small $\epsilon>0,$ but it is conjectured to be $\leq \frac{1}{4}.$ \ \\

\textbf{Problem 1b, communicated by Noah Kravitz:} The discrete version of the previous question is as follows. Given a set $A\subset \mathbb Z / (p\mathbb Z),$ with the property that $A + A - 2A \neq \mathbb Z / (p\mathbb Z),$ where $2A$ denotes dilation by 2, not a sum set, what bounds can we have on the size of $A$? It is known that $|A|\leq \frac{p}{3},$ but it is conjectured to be $\leq \frac{p}{4}.$ Finite field/ring variants also appear to be open. \ \\

\textbf{Problem 2, communicated by Noah Kravitz:} Given a natural number $n$ with MANY divisors ,(for example, suppose it is divisible by $m!$ for some large $m$), consider a set $A \subset \mathbb Z / (n\mathbb Z)$ with the property that it avoids these two configurations:
\begin{itemize}
\item $\{x, 2x\}$
\item $\{a, b, a+b\},$ whenever $a$ and $b$ are coprime, and $\gcd(a,b,n)=1.$
\end{itemize}
Note that the second condition is weaker than ``sum free," as the triple of numbers $\{a,b,a+b\}$ is allowed in $A$ provided that the $\gcd$ relations do NOT hold. The question is, how large can such $A$ be? The conjecture is that $|A|\leq \frac{n}{2}$ with sharpness achieved by $A$ being the odd numbers between 1 and $n,$ when $n$ is even. \ \\

\textbf{Problem 3, communicated by Marc Technau:} Marc heard this from Christoph Glasser and Titus Dose. We first define what it means for a subset of $\mathbb Z$ to be happy. All singleton sets are happy. Complements of happy sets are happy. Given two happy sets, $A$ and $B$, the following sets are also happy: $A \cup B, A \cap B, A+B,$ and $A\cdot B,$ where the last two denote sum set and product set, respectively. The question is, is $\mathbb N$ happy? Some notes:
\begin{itemize}
\item It is straightforward to show that $\mathbb Z$ and $\varnothing$ are happy.
\item $\mathbb N$ is not $\mathbb Z,$ or the problem would be easy.
\item Someone in the room conjectured this would be solved by Friday.
\item The positive and negative powers of two are happy, as $\overline{(2\mathbb Z+ 1)\mathbb Z}$ is happy.
\item ``Most" sets must be unhappy, as there are only countably many happy sets.
\end{itemize} \ \\

\textbf{Problem 4, communicated by Steven Senger (repeat from previous years):} Given a large finite subset, $A\subset \mathbb R,$ can we show that $|AA|=o(|(AA+1)(AA+1)|)$?


\newpage

\section{CANT Problem Sessions: 2024}

\subsection{Problem Session I: Wednesday, May 22nd, 2024}

\large
\noindent Chaired by Steven J. Miller (Wednesday, May 22nd, 2024). \\ \
\normalsize

\noindent \textbf{Problem 1:} Proposed by Steven Senger: StevenSenger@gmail.com: Find an asymptotic upper bound for the number of quadruples $(a, b, c, d)$ of points in any large finite point set so that $a, b, c, d \in \mathbb{R}^3$ the angle of $\angle abc = \angle bcd$ equals 90 degrees. It is known that there are $O(n^{10/3})$ points, conjectured that the truth is of size $n^3$. \\ \

\noindent \textbf{Problem 2:} Asked by Steven J. Miller to Daniel Flores: Can you extend to $K$-multimagic cubes or higher problems? Answer: Gets very technical, and in $d$-dimensions the bound is independent of $d$. Also: I believe that a nontrivial 4-multimagic square of order 20 should exist. Maybe someone better at computation can find it. \\ \

\noindent \textbf{Problem 3:} Proposed by Steven J. Miller (to be done by his SMALL 2024 students): sjm1@williams.edu: Can you bypass the calculations for other problems and use the Theory of Normalization Constants to solve problems related to Zeckendorf decompositions?\\ \

\noindent \textbf{Problem 4:} Asked by Steven to Steven: Generalize Problem 1 to look at other configurations. How critical is 90 degrees? Can look at other configurations. \\ \

\subsection{Problem Session II: Thursday, May 23rd, 2024}

\large
\noindent Chaired by Kevin O'Bryant (Thursday, May 23rd, 2024). \\ \
\normalsize

\noindent \textbf{Problem 1:} Proposed by Leonid Fel: lfel@technion.ac.il: Partitions with constraints Gaussian polynomials (posed before, but now with more specificity): Consider
\[\sum_{k=1}^m kx_k\ =\ s,\]
where the $x_k$ are nonnegative. This gives partitions of $s$. For example, $$5 \ = \ 5+0 \ = \ 4+1 \ = \ 2+2+1 \ = \ \cdots.$$ This is called $W(s)$. We now consider this with some constraints. Namely, suppose $\sum x_k \leq m$. The basic partitions come from the function
\[M \ = \ \frac{1}{\prod_{k=1}^\infty(1-z^k)} \ =\ \sum W(s)z^s.\]
Following Shu, we write the generating function
$$G(m,n,z)\ =\ \frac{ \prod_{k=1}^{m+n} (1-z^k) }{ \prod_{k=1}^n (1-z^k)\prod_{k=1}^m (1-z^k) } \ = \ \sum_{s=0}^{mn} P_k^m(s) z^s. $$
For example, $$G(2,3,z) \ = \ 1+z+2z^2+2z^3+2z^4+z^5+z^6.$$
Famously (1980s), this distribution is unimodal, but not log-concave.

\begin{itemize}
\item $1<P(1)\leq \dots P \dots > P_{mn-1}>1$.
\item $0<s<mn$.
\end{itemize}

Naturally motivated by mathematical physics (referenced in Leonid's previous talks): Following Wilf, not conjecturing, just asking.
Let $m_1 n = m_2 n = 0 $ mod 2, and $n\geq 3.$ Then is it true that
\[P_n^{m_1}\left(\frac{m_1 n}{2}\right) + P_n^{m_2}\left(\frac{m_2 n}{2}\right)\ \leq\ P_n^{m_1+m_2}\left(\frac{(m_1+m_2)n}{2}\right)?\]
An example of this inequality holding: $P_3^2(3) = 2, P_3^4(6)=5,$ and $P_3^6(9)=8.$

Addendum/motivation:
For $x^3=(x_1,x_2,x_3)\in \mathbb R^3$, $I_1,I_2, I_3$ are invariants. Let $$S(\lambda, x^3) \ =\ \lambda^3+AI_1\lambda^2-AI_2\lambda-I_3.$$ This is invariant when we map $\lambda$ to $\lambda^{-1}$. Let $$S_2(\lambda,x^3) \ = \ \lambda^6+AI_1\lambda^5+(BI_1^2+CI_2)\lambda^4-(BI_2^2+CI_1I_3)\lambda^2-AI_2I_3\lambda-I_3^2.$$
These are invariant polynomials. $P_n^m(s)$ is the number of Young diagrams partitioning $s$ into at most $n$ parts of size at most $m$. \\ \

\noindent \textbf{Problem 2:} Proposed by Steven Senger: stevensenger@gmail.com: HRT conjecture in finite fields: The classical HRT Conjecture (Heil-Ramanathan-Topiwala) claims that for any reasonable ($L^2$) function $f(x)$ on the reals any set of time-frequency shifts $f(x)$ should be linearly independent. Here a time-frequency shift by $(a,b)$ is $\chi(bx)f(x-a)$, where $\chi(z) = e^{-2\pi 1 z}$. My question is, to what extent can this hold in finite fields? Reference: Okoudjou and Oussa: \bburl{https://arxiv.org/pdf/2110.04053}. \\ \

\noindent \textbf{Problem 3:} Proposed by Noah Lebowitz-Lockard: nlebowi@gmail.com: Start with the game Nim: Given a partition $a_1,a_2,\dots>a_k$ you take a turn by decreasing any one of the $a_i$. The game ends when all are zero, and the last person to play wins. Define $\ell(n)$ to be the number of partitions of $n$ corresponding to a loss for Player 1 (under optimal play).  Example: $(1,1)\rightarrow (1,0)\rightarrow (0,0)$, so Player 2 wins.

Question: Can we bound $\ell(n)$ for $n$ even? Note: Player 1 can always win for $n$ odd. Noah can show $\ell(n) \gg \frac{p(n)}{n}$, where $p(n)$ is the number of partitions of $n$. Conjecture: $\ell(n) \gg \frac{p(n)}{\sqrt{n}}$. Noah's approach is by writing the numbers in binary and looking at the sum of each column. Each column is even if and only if this is a losing position for Player 1. \\ \

\noindent \textbf{Problem 4:} Proposed by George Yuan: qluan21@g.ucla.edu: For relatively prime $a,b\in\mathbb{N}$, consider
    \begin{equation}\label{eq1}
    ax+by\ =\ \frac{(a-1)(b-1)}{2},
    \end{equation}
    \begin{equation}\label{eq2}
    ax+by+1\ =\ \frac{(a-1)(b-1)}{2}.
    \end{equation}
M. Beiter and H. V. Chu have previously shown that for each pair of relatively prime $a,b\in\mathbb{N}$, exactly one of these equations will have a solution, and the solution will be unique.
Let $$\Gamma(a, b) \ := \ \begin{cases}
        1 & \text{if $(a/d, b/d)$ uses Equation $\eqref{eq1}$}\\
        2 & \text{if $(a/d, b/d)$ uses Equation $\eqref{eq2}$}
    \end{cases}.$$
Let
    \[G(x) \ :=\ \frac{\#\{(a, b) \in \mathbb{N}^2: 1 \leq a \leq b\leq x, \gcd(a, b) = 1, \Gamma(a, b) = 1\}}{\#\{(a, b) \in \mathbb{N}^2: 1 \leq a \leq b\leq x\}}\]
    and
    \[H(x)\ :=\ \frac{\#\{(a, b) \in \mathbb{N}^2: 1 \leq a \leq b\leq x, \Gamma(a, b) = 1\}}{\#\{(a, b) \in \mathbb{N}^2: 1 \leq a \leq b\leq x\}}.\]
What will $\lim_{x \to \infty} G(x)$ and $\lim_{x \to \infty} H(x)$ be (if they indeed exist)?

\subsection{Problem Session II: Thursday, May 23rd, 2024}

\large
\noindent Chaired by Steven Senger (Friday, May 24th, 2024). \\ \
\normalsize


\newpage

\section{CANT Problem Sessions: 2025}

\noindent \large \textbf{CANT Problem Session, May 21, 2025} \normalsize \\

\noindent Communicated by Steven J. Miller: sjm1@williams.edu \\ \

\noindent \textbf{Problem 1, communicated by Steven J. Miller:}  From Miller's talk: try to generalize the observed behavior in the different regimes to higher dimensions (this may be done by Miller's Polymath Jr group in Summer 2025).

 \ \\

\noindent \textbf{Problem 2, communicated by Steven J. Miller:} Come up with an explicit construction for a family of MSTD sets which has positive density, or at least better than current best results (on the order of $2^n / n^2$).

 \ \\

\noindent \textbf{Problem 3, communicated by Mel Nathanson:} Given a set $A$ with prescribed values for $|A|$ and $|A+A+A|$, what are the options for $|A+A|$? What if we prescribe $|k_iA|$ for $0 < k_1 < k_2 < \cdots < k_ell$, what can we say about the other $|kA|$? Note that for $h$ sufficiently large, if $A$ is normalized to have smallest element 0 or 1, largest element $a$ and no $d>1$ divides all elements then $|hA| = h a + {\rm const}$ for some constant.

 \ \\

\noindent \textbf{Problem 4, communicated by Samuel Allen Alexander:} Find a sequence of sets with $A_{i-1} \subset A_i$ that alternate being sum and difference dominated. Observations by Stevens (Miller, Senger): should be trivial if sets are beyond exponentially growing so the addition can dwarf the previous. Also trivial if can fill in the $A_{i-1}$ and make that a middle.

 \ \\

\noindent \large \textbf{CANT Problem Session, May 23, 2025} \normalsize \\

\noindent Communicated by Kevin O'Kevin obryant@gmail.com  \\ \

\noindent \textbf{Problem 5, communicated by Brian Hopkins (bhopkins@saintpeters.edu):} Let  $p(n,3)$ be the number of partitions of $n$ with 3 parts. We have 
$$p(n,3) \ =\ \left[\frac{n^2}{12}\right],$$ the nearest integer to $n^2/12$. Given a Pythagorean Triple $a^2+b^2=c^2$, we know that $$p(a,3) + p(b,3)\ = \ p(c,3).$$ Give a combinatorial proof, please! \\ \

\noindent \textbf{Problem 6, communicated by Kevin O'Bryant (kevin.obryant@csi.cuny.edu):} Recall that a $B_h$-set is a set of integers $\mathcal{A}$ with no solutions to
\[x_1+\cdots+x_h \ =\ y_1+\cdots+y_h,\qquad x_1\le \cdots \le x_h,y_1\le \cdots y_h, \ \ (x_1,\ldots,x_h)\ \neq\ (y_1,\ldots,y_h),\ \ x_i\in \mathcal{A}, \ \ y_i
\in \mathcal{A}.\]

For a positive integer $h$, set $\gamma_0=0,\gamma_1=1$, and for $k\ge 1$, let $\gamma_k=\gamma_k(h)$ be the least integer $x$ so that
$x>\gamma_{k-1}$ and $\{\gamma_0,\ldots,\gamma_{k-1},x\}$ is a $B_h$-set. For $h=2$, this is the Mian-Chowla set, which is famously
intractable---literally nothing is known that is not computational beyond the straightforward combinatorial $\gamma_k(2) = O(k^3)$. Other fixed
$h$ seems similar, although there has been little effort, even computationally. Instead of fixing $h$, we ask you to consider the problem for
fixed $k$. It is easy to see that $\gamma_2(h)=h+1$, Nathanson [1] has shown that $\gamma_3(h)=h^2+h+1$, and Nathanson \& O'Bryant [2] have
shown that $\gamma_4(h) = \frac12 (h^3+2h^2+3h+2)$ if $h$ is even and $\gamma_4(h) = \frac12 (h+1)^3$ if $h$ is odd.

Our question is whether $\gamma_5(h)$ is also polynomial in $h$ when restricted to $h\equiv0\pmod{6}$. If so, it must be quartic [3]. More
generally, is there a modulus $m_k$ so that $\gamma_k(h)$ is given by a degree $k-1$ polynomial in $h$ for sufficiently large $h$ when
restricted to congruence classes modulo $m_k$? That is, is $\gamma_k(h)$ a quasipolynomial for each $k$? \\ \

\noindent [1] M. B. Nathanson, The third positive element in the greedy $B_h$-set, \bburl{https://arxiv.org/abs/2310.14426}. \\ \

\noindent [2] M. B. Nathanson and K. O'Bryant, The fourth positive element in the greedy $B_h$-set, \bburl{https://arxiv.org/abs/2311.14021}. \\ \

\noindent [3] K. O'Bryant, Bounds for Greedy $B_h$-sets, \bburl{https://arxiv.org/abs/2312.10910}. \\ \

\ \\

\noindent \textbf{Problem 7, communicated by Taylor S. Daniels (daniel84@purdue.edu):}  For each prime $p \equiv 1 \mod{4}$, let $\gamma = \gamma(p) \in \{1,\ldots,p-1\}$ be the \emph{smallest} primitive root $(\mathrm{mod}\,p)$, and let $i = i(p)$ be the integer satisfying
    \[
        i \equiv \gamma^{\frac{p-1}{4}} \mod{p} \qquad\text{and}\qquad 0 < i < p.
    \]
By way of Wilson's theorem we observe that
$$ \left(\prod_{0 < \mu < p \atop \mu {\rm\ even}} \mu \right)^{2}\ =\ \left(\prod_{0 < \mu < p \atop \mu {\rm\ even}} \mu \right)\left(\prod_{0 < \mu < p \atop \mu {\rm \ odd}} (p-\mu) \right) \ \equiv\ (-1)^{\frac{p-1}{2}} \prod_{0 < \mu < p} \mu \equiv -1 \mod p, $$
whereby
    \[
        \prod_{0 < \mu < p \atop \mu {\rm\ even}} \mu\ \equiv\ \pm i \mod{p}.
    \]
On the other hand,
    \[
        \prod_{0 < \mu < p \atop \mu {\rm\ even}} \mu \ =\ 2^{\frac{p-1}{2}}(\tfrac{p-1}{2})!\ \equiv\ \pm (\tfrac{p-1}{2})! \mod{p},
    \]
and it follows then that
    \(
    \label{eq:pHalfFact(modp)}
        (\tfrac{p-1}{2})!\ \equiv\ \pm i \mod{p}.
    \)

\ \\
\noindent \textbf{Definition:}   For all $p \equiv 1 \mod{4}$, let $i=i(p)$ and $\gamma=\gamma(p)$ be defined as previously, and with these let $\epsilon(p) \in \{0,1\}$ be defined so that
        \[
            (\tfrac{p-1}{2})!\ \equiv\ (-1)^{\epsilon(p)} i \mod{p}.
        \] \ \\

It is easy to have a computer check values of $\epsilon(p)$ for small $p \equiv 1 \mod{4}$; in particular we find that
    \(
    \label{eq:eps(p)-0}
        \epsilon(p)\ =\ 0 \quad \text{for $p = 5$, $17$, $29$, $37$, $61$, $73$, $89$, $\cdots$},
    \)
and
    \(
    \label{eq:eps(p)-1}
        \epsilon(p) \ =\ 1 \quad \text{for $p = 13$, $41$, $53$, $101$, $113$, $137$, $\cdots$}.
    \)
Unfortunately, we are unaware of a general pattern for $\epsilon(p)$, and neither sequence of primes above appear in the OEIS. \\ \

\noindent \textbf{Problem:}    Determine characterizations of the sequences in \eqref{eq:eps(p)-0} and \eqref{eq:eps(p)-1}. \\ \

\noindent \textbf{Remark:}     We originally specified that $\gamma$ should be the \emph{smallest} primitive root (mod $p$). When $p \equiv 1 \mod{4}$, if $0 < g < p$ is any primitive root (mod $p$), then so is $p-g$. Thus, one could alternatively start by defining $\gamma$ to be, for instance, the smallest \emph{even} or \emph{odd} primitive root (mod $p$).
    The derivation of \eqref{eq:pHalfFact(modp)} is clearly unaffected by this, but this \emph{does} change the resulting sequences of primes in \eqref{eq:eps(p)-0} and \eqref{eq:eps(p)-1}. However, as with \eqref{eq:eps(p)-0} and \eqref{eq:eps(p)-1}, the sequences resulting from choosing $\gamma$ to be the smallest even (or smallest odd) primitive root(s) (mod $p$) do not appear in the  OEIS either. \\ \

\end{document}